\documentclass[opre,nonblindrev]{informs3} % current default for manuscript submission

\DoubleSpacedXI % Made default 4/4/2014 at request
\usepackage{endnotes}
\let\footnote=\endnote
\let\enotesize=\normalsize
\def\notesname{Endnotes}%
\def\makeenmark{$^{\theenmark}$}
\def\enoteformat{\rightskip0pt\leftskip0pt\parindent=1.75em
  \leavevmode\llap{\theenmark.\enskip}}

\usepackage{natbib}
 \bibpunct[, ]{(}{)}{,}{a}{}{,}%
 \def\bibfont{\small}%
\TheoremsNumberedThrough     % Preferred (Theorem 1, Lemma 1, Theorem 2)
\ECRepeatTheorems

\EquationsNumberedThrough    % Default: (1), (2), ...
\usepackage{amsmath,amssymb}
\usepackage{float}
\usepackage{graphicx}
\usepackage{enumerate}
\usepackage[toc,page]{appendix}
\usepackage{url} % not crucial - just used below for the URL 
\usepackage[ruled,vlined,linesnumbered,resetcount]{algorithm2e}

\usepackage{dsfont}
\usepackage{float}
\usepackage{amsmath,amssymb}
\usepackage{graphicx}
\usepackage{enumitem}
\usepackage{chngcntr}
\usepackage{apptools}
\AtAppendix{\counterwithin{lemma}{section}}

\allowdisplaybreaks

\newcommand{\epc}{\hspace{1pc}}
\def\g{g^*_\pi}
\def\h{h^*_\pi}
\def\gn{\widehat g_{N}^{\pi, \beta}}
\def\hn{\widehat h_{N}^{\pi}}
\def\VC{\operatorname{VC}}
\def\EB{\operatorname{EB}}

\def\S{\mathcal{S}}

\def\EE{\mathbb{E}}
\def\TV{\operatorname{TV}}
\def\PP{\mathbb{P}}
\def\U{\mathcal{U}}

\def\G{\boldsymbol{G}}
\def\F{\mathcal{F}}
\def\Hn{\widehat H^{\pi}_N}

\def\Hpi{H^\pi}
\def\T{\mathcal{T}}

\def\goes{\rightarrow}

\def\E{\mathcal{E}}

\def\W{\mathcal{W}}

\def\f{\boldsymbol{f}}
\def\A{\mathcal{A}}
\def\D{\mathcal{D}}

\def\R{\mathbb{R}}
\def\Pn{\mathbb{P}_n}
\def\PN{\mathbb{P}_N}
\def\Mn{\mathbb{M}_n}

\def\Upin{\widehat {U}^{\pi, \beta}_n}
\def\UpiN{\widehat {U}^{\pi, \beta}_N}
\def\Upi{{U}^{\pi, \beta}}

\def\P{\mathcal{P}}

\def\RR{\mathcal{R}}

\def\leqconst{\lesssim}
\def\geqconst{\gtrsim}

\newcommand{\pushright}[1]{\ifmeasuring@#1\else\omit\hfill$\displaystyle#1$\fi\ignorespaces}
\newcommand{\innerprod}[2]{\langle#1, #2 \rangle}

\newcommand{\norm}[1]{\|#1\|}

\newcommand{\abs}[1]{|#1|}
\newcommand{\dkh}[1]{\{#1\}}
\newcommand{\bigdkh}[1]{\big\{#1 \big\}}
\newcommand{\Bigdkh}[1]{\Big\{#1\Big\}}
\newcommand{\bigabs}[1]{\big|#1\big|}
\newcommand{\Bigabs}[1]{\Big|#1\Big|}

\newcommand{\peng}[1]{{\color{blue}{#1}}}

\newcommand{\indicator}[1]{ \mathds{1}_{\{#1\}}}

\newcommand{\mynu}{\nu}

\def\EE{\mathbb{E}}

\def\f{\boldsymbol{f}}

\def\R{\mathbb{R}}
\def\Regret{\operatorname{Regret}}

\def\Rem{\operatorname{Rem}}

\def\I{\,\mathcal{I}}
\def\ones{\mathbf{1}}

\def\B{\mathcal{B}}
\def\argmin{\operatorname{argmin}}
\def\argminb{\operatorname*{argmin}}
\def\argmax{\operatorname{argmax}}
\def\N{\mathcal{N}}

\def\X{\mathcal{X}}

\def\jiao{\cap}

\newcommand{\samfixed}[1]{}

\def\Q{\mathcal{Q}}
\def\given{\, | \,}
\def\Given{\, \Big| \,}
\def\GG{\mathbb{G}}

\def\transpose{\top}

\renewcommand{\G}{\mathcal{G}}

\def\Mn{\widehat{M}_N}

\usepackage{tcolorbox}

\usepackage{mathrsfs}

\begin{document}
\RUNAUTHOR{Qi and Liao}
\RUNTITLE{Robust Batch Policy Learning}

\TITLE{Robust Batch Policy Learning in Markov Decision Processes}

\ARTICLEAUTHORS{%
\AUTHOR{Zhengling Qi}
\AFF{Department of Decision Sciences, George Washington University,
\EMAIL{qizhengling@gwu.edu}} %, \URL{}}
\AUTHOR{Peng Liao}
\AFF{Department of Statistics, Harvard University, \EMAIL{pengliao@g.harvard.edu}}
}

\ABSTRACT{%
We study the offline data-driven sequential decision making problem in the framework of Markov decision process (MDP). In order to enhance the generalizability and adaptivity of the learned policy, we propose to evaluate each policy by a set of the average rewards with respect to distributions centered at the policy induced stationary distribution. Given a pre-collected dataset of multiple trajectories generated by some behavior policy, our goal is to learn a robust policy in a pre-specified policy class that can maximize the smallest value of this set. Leveraging the theory of semi-parametric statistics, we develop a statistically efficient policy learning method for estimating the defined robust optimal policy. A rate-optimal regret bound up to a logarithmic factor is established in terms of total decision points in the dataset. 
%Our regret theoretical guarantee subsumes the long-term average reward MDP setting as a special case and can be possibly extended to the discounted setting. (PL: I don't think this worth mentioning in the abstract)
}%


\maketitle
\newpage

%\peng{Check the use of the phrase ``iid setting'' and dependent data.}
%\peng{Decide if we want to use the work ``indefinite''}
%\begingroup
\section{Introduction}

% guarantee robust performance
%
%be applied for an indefinite time horizon and protect against different initial distribution with uniformly guaranteed robust performance.

An essential goal in data-driven sequential decision making problems is to construct a policy that maximizes the average reward over a certain amount of the time. Depending on the applications, the duration of the policy for use in the future is often unknown and is likely to be different from what we consider at the stage of constructing or optimizing the policy. See two motivating examples below. %\peng{For example, ... (more specific? this paragraph is too vague for nonRL people)} 
Furthermore, the performance measure used in learning the optimal policy often depends on the choice of initial state's distribution. For example, one widely used performance measure is based on the value function, which is the conditional expectation of discounted sum of rewards starting at a given state. To form the objective in optimizing the policy within a specified policy class, one often averages over initial state's distribution, which may change when implementing the learned policy in the environment. Therefore,  given the uncertainty of deploying policies in practice, it is critical to learn a policy with strong generalizability and adaptivity.  Motivated by this, our goal in this paper is 
to learn a robust policy in the sense that it can guarantee the uniform performance over the unknown planning horizon and the distributional change in the initial state. 

%In this work, we consider the batch policy learning problem in indefinite-horizon Markov decision processes (MDPs). Given the pre-collected data of multiple trajectories consisting of states, actions and rewards, our goal is to learn a policy with guaranteed performance that can be used in the future without knowing the duration of the use of the policy and the reference distribution at the future initial state.

%the length of horizons over which a agent wants to maximize average rewards and the reference distribution at the starting state may be different. Given the uncertainty of time duration of a policy being used in the future and reference distribution at the starting state, one may want to learn a robust policy from the pre-collected data that can guarantee the uniform performance against these uncertainties. 

We consider a batch reinforcement learning (batch RL) problem under the framework of Markov decision process (MDP), where data are pre-collected in the form of multiple trajectories consisting of states, actions and rewards. Recently there is an increasing interest in studying batch RL (e.g., \cite{ernst_tree-based_2005,antos2008fitted,mnih2015human,farahmand2016regularized,dabney2018distributional,le2019batch,kallus2020statistically,kumar2020conservative,jin2021pessimism} among many others) as effective solutions for finding optimal decision rules by leveraging the rich observational data in many applications (e.g., \citep{komorowski2018artificial,luckett2019estimating,levine2020offline,shi2020statistical} among many others). 

Our work is also inspired by the following two real-world examples. The first one is the recently emerging mobile health (mHealth) applications. An essential goal of mHealth is to deliver a customized intervention via notification or text message at the right time and the right location for helping individuals make healthy decisions \citep{nahum2016just}. Prior to the actual implementation of interventions, pilot studies are often first conducted to test the software and evaluate multiple intervention components using randomization \citep{klasnja2015microrandomized,Liaoetal2015}. The data collected from these studies can also be used to estimate a good ``warm-start'' policy for the use in the future. It is thus important for the learned policy to ensure decent performance across different individuals and the length of time that the policy is used. The second motivating example comes from the inventory management \citep{powell2007approximate}. A retailer, facing the uncertainty of consumers' daily demands, needs to decide how much inventory to purchase every day. Due to the unknown information of stochastic lead time \citep{kaplan1970dynamic}, a robust ordering policy that can protect against uncertainty in delivery is desirable in order to maintain inventory level so that stockout and holding costs are minimized.

In order to learn a desired robust policy,
we propose to evaluate a policy by the average rewards with respect to a set of distributions centered at the policy induced stationary distribution. Under standard mixing conditions, we show that this set contains the average rewards over different lengths of time-horizons and reference distributions of initial state. Such appealing property of this set motivates us to perform policy optimization that improves the minimal value in this set. As a result, we can guarantee the robust performance of the learned policy when
implemented in the future against the uncertainty characterized by the proposed set. To the best of our knowledge, such criterion for policy learning in MDPs has not been studied before. Thanks to the celebrated convex duality theory, policy optimization under our proposed novel criterion can be formulated as an M-estimation problem in statistics (see, for example, \cite{van2000asymptotic}) and hence the original max-min problem becomes more computationally tractable (see Theorem \ref{duality} for details). More importantly, based on this equivalent representation we develop a statistically efficient policy learning method to estimate the optimal policy under the proposed robust criterion over a parametrized policy class. In particular, we show that our proposed algorithm can achieve the rate-optimal regret bound up to a logarithm factor in terms of the number of trajectories and the number of decision points at each trajectory, thus efficiently using the pre-collected data and also breaking the curse of the time-horizon \citep{kallus2019efficiently}. Our theoretical result generalizes the previous work by \cite{liao2020batch}, which studied the policy learning under the long-term average reward, and can be extended to the discounted sum of rewards setting. To the best of our knowledge, this is the first in-class near-optimal regret bound established in the literature of batch RL in terms of the total number of decision points, which itself may be of independent interest.

 Our approach can also be viewed as an example of distributionally robust optimization (DRO). DRO has recently attracted a lot of interests in the community of machine learning and statistics due to its superior performance in terms of generalization. See some recent literature such as \citep{pflug2014multistage,wozabal2014robustifying,gao2016distributionally,blanchet2019quantifying,esfahani2018data}   and two recent review papers by \citep{kuhn2019wasserstein,rahimian2019distributionally}. In the MDPs, DRO has been mainly studied in the setting of discounted sum of rewards. The major discussion is focused on the uncertainty of the temporal difference and the corresponding parameter estimation. See for example \citep{xu2010distributionally, smirnova2019distributionally,abdullah2019wasserstein,derman2020distributional} for more details. In particular, \citep{smirnova2019distributionally} established a (sub-optimal) sample complexity result for their distributionally robust modified policy iteration method in the setting of finite state and action spaces.
It is also known that there is a strong connection between DRO and risk measure \citep{ben2007old}. In the risk-sensitive sequential decision making, one line of research is to modify the criterion of searching a policy by taking risky scenarios into consideration. See the early papers by \citep{sobel1982variance,filar1989variance}.  Another line of research is to control the uncertainty of the exploration process such as temporal differences (e.g., \cite{mannor2011mean,gehring2013smart}). See some recent developments in risk-sensitive reinforcement learning such as \citep{prashanth2013actor,shen2014risk,chow2015risk,tamar2015policy,prashanth2018risk,qi2019estimation, qi2019estimating,zhong2020risk}. Our proposed criterion can be regarded as using DRO or risk measure to robustify the policy optimization objective. Thus our method inherits the nice property of DRO in improving the generalizability of the learned policy to the new and unseen data. Compared with existing literature on DRO and risk-sensitive RL, we focus on improving the average rewards over varying time horizons with the unknown initial state distribution. To the best of our knowledge, this has not been studied before. 
{More notably, few existing works in the literature considered the statistical efficiency of algorithms (i.e., how to efficiently use the data), which is essential in batch RL. As the amount of available training data is often limited, in contrast with the online setting, it is necessary to develop a data-efficient learning method to perform policy optimization. }

\begin{comment}

Our main contributions can be summarized into three folds. First, we propose a robust criterion that evaluates a policy by the average rewards with respect to a set of distributions centered at the stationary distribution. Based on this criterion, by performing policy optimization that improves the minimal value of this uncertain set, we can guarantee the robust performance of the learned policy when implemented in the future; Second, by leveraging the convex duality theory, policy optimization under our proposed novel criterion can be formulated as an M-estimation problem in statistics and hence the original max-min problem becomes more computationally tractable. Based on this equivalent representation, we develop a statistically efficient policy learning method to estimate the optimal policy under the proposed robust criterion over a parametrized policy class; Lastly, we show that our proposed algorithm can achieve the rate-optimal regret bound up to a logarithm factor in terms of the number of trajectories and the number of decision points at each trajectory. To the best of our knowledge, this is the first in-class near-optimal regret bound established in the literature of offline reinforcement learning in terms of the total number of decision points, which is an important theoretical contribution.

\end{comment}

The rest of the paper is organized as follows. In Section 2, we introduce the framework of the time-homogeneous MDP, related concepts and notations. In Section 3,  we formally introduce a robust average reward criterion that can be used to improve the generalizability of the learned policy.  We then discuss our statistically efficient learning method to estimate the optimal policy under our proposed robust criterion in Section 4. In section 5, we provide strong theoretical guarantees for our proposed method including the uniformly finite sample error bounds for nuisance functions estimation, the statistical efficiency of our proposed estimator in evaluating a policy and the strong finite-sample regret bound of our learned policy. All these results are seemingly new in the current literature. In Section 6, we use a simulation study to demonstrate the promising performance of our proposed method. We provide some discussions and point out some interesting future research directions in Section 7. All proofs of technical results and details of computation can be found in the Supplementary Material.

\section{Framework}
\subsection{Time-homogeneous Markov Decision Processes}
In this section, we briefly introduce discrete time homogeneous MDPs and the necessary notations. For a comprehensive description, we refer to \citep{puterman1994markov} and \citep{hernandez2012further}.  Denote $\S$ as the state space, and $\A$ as a finite action space. Let $\B(\S)$ and $\B(\S \times \A)$ be the family of Borel subsets on $\S$ and $\S \times \A$ respectively. We assume $\B(\S \times \A)$ contains all pairs of $(s, a)$ for every $(s, a) \in (\S \times \A)$. %\peng{I dont understand this sentence \QZL{This is correct to avoid some technical issues. See Meyn's book}}.
We further define the stochastic kernel $P$ on $\S$ given a measurable subset of $\S \times \A$. This means $P(\bullet | s, a)$ is a probability measure on $\B(\S)$ for every $(s, a) \in \S \times \A$ and $P(B | \bullet, \bullet)$ is a non-negative measurable on $\S \times \A$ for every $B \in \B(\S)$. We denote $t = 1, 2, 3,\cdots$, as a series of discrete time steps. The time-homogeneous MDP process begins as $(S_1, A_1, S_2,\cdots, S_t, A_t, \cdots)$ on $\Omega = \Pi_{t = 1}^\infty (\S_t \times \A_t)$, and measurable with respect to $\mathbb{F} = \otimes_{t = 1}^\infty \B(\S_t \times \A_t)$, with some probability measure $\PP$, where $(\S_t \times \A_t)$ is a copy of $(\S \times \A)$. Denote the history up to $k$-th time as $H_k =S_1 \times \Pi_{t = 1}^{k-1}(S_t \times A_t)$ for $k \geq 2$ and $H_1 = \S_1$. The distribution $\PP$ satisfies that for $t \geq 2$, $\PP(S_{t+1}\in B | A_t = a_t, H_t = h_t) = P(B \, | \, s_t, a_t)$ for every $B \in \B(\S)$ and $h_t = (s_1, a_1, s_2, \cdots, s_t) \in H_t$, thus satisfying Markovian and time-homogeneous properties. We assume the reward only depends on the current state, that is, $R_t = \mathcal{R}(S_t)$, where $\mathcal{R}$ is a known measurable function defined over $\S$. In addition, we assume $\mathcal{R}$ is uniformly bounded by a positive constant $R_{\max}$.  Such assumption on the reward was commonly used in the literature, such as \cite{baxter2001infinite}. Other forms of reward will be discussed in Section 6.

%and the transition probability measure $P$ has a density $p$.

The tuple $(\S, \A, P)$ is usually called an MDP. In this work, we focus on the time-invariant Markovian policy $\pi$, which is a function mapping from the state space $\cal S$ into a probability distribution over the action space $\cal A$. More specifically, $\pi(a \, | \, s)$ denotes the probability of selecting the action $a$ given the state $s$. Together, an MDP $(\S, \A, P)$, a policy $\pi$ and an initial state distribution $\nu$ define a joint probability measure $\PP^\pi$ over $(S_1, A_1, \cdots, S_t, A_t, \cdots)$ such that  (1) $\PP^\pi(H_1 \in B) = \nu(S_1 \in B)$ for every $B \in \B(\S)$; (2) for $t \geq 2$, $\PP^\pi(S_{t+1}\in B | A_t = a_t, H_t = h_t) = P(B \, | \, s_t, a_t)$ for every $B \in \B(\S)$ and (3) $\PP^\pi(A_t = a_t | H_t) = \pi(a_t | s_t)$. We use $\EE_\pi$ to denote the expectation with respect to $\PP^\pi$. For simplicity, throughout this paper, we assume all probability measures have probility densities with respect to the Lebesgue measure.

\subsection{Batch RL}
%\peng{we use $\nu$ to refer the initial distribution in the training data? if so we need to say this below}
 
In the batch setting, we are given a training dataset $\D_n$ collected from previous studies that consists of sample size $n$ independent and identically distributed (i.i.d.) trajectories of length $T_0$, i.e., %(Without loss of generality, we assume lengths of trajectories in $\D_n$ are equal) 
%peng{Why we say without loss of generality here? I don't think this is a ``without generality'' statement. }\QZL{I delete this}: 
\[
\D_n = \left\{ D_i \right\}_{i = 1}^n = \left\{S_1^i, A_1^{i}, S_2^{i}, \cdots, S_{T_0}^i, A_{T_0}^i, \cdots, S_{T_0+1}^i \right\}_{i = 1}^n.
\]
Each trajectory as the form of $D = \left\{S_1, A_1, S_2, \cdots, S_{T_0}, A_{T_0}, S_{T_0+1} \right\}$ is assumed to be generated by some behavior policy $\left\{\pi_{bt}(\bullet \, |\, H_t)\right\}_{t = 1}^{T_0}$, where $\pi_{bt}(\bullet \, |\, H_t)$ maps the history $H_t$ to a probability mass function defined on $\A$. The distribution of the initial state in $D$ is denoted by $\nu$. In our theoretical analysis given in Section 4, we assume the behavior policy being time-stationary. But implementing our method introduced below does not need this assumption, so we let the behavior policy be history-dependent to keep its generalization. 
%\peng{We assume the behavior policy is time-stationary. use $\pi_b$ throughout}

%\QZL{Be lazy, so I delete this}\peng{(The initial distirbution of state is denoted by $\nu$)?} \peng{In this work we assume the trajectory length is identical, but this can be easily generalized to incorporate different length (maybe explain why?)} 

A primary goal of batch RL is to learn a policy in a policy class $\Pi$
% in a pre-specified class of policies $\Pi$ 
that maximizes the average reward over some time horizon $T$ (i.e., planning horizon) and with respect to some state distribution $\GG$ (i.e., a reference distribution) \citep{puterman1994markov}. More specifically, for a given policy $\pi$ and an initial state $s$, we define its average reward as
%\peng{What about we define}
\begin{align*}
	\eta^\pi_{T}(s) = \EE_\pi \left[\frac{1}{T}\sum_{t=1}^{T}R_{t} \, \,\middle\vert\, \, S_1=s \right].
\end{align*}
Then the integrated average reward with respect to a reference distribution $\GG$ is defined as
\begin{align}\label{finite average}
    \eta^\pi_{T}(\GG) = \int \eta^{\pi}_{T}(s) d\GG(s),
\end{align}
%\begin{align}\label{finite average}
%	\eta^\pi_{T_1}(\GG) = \EE_{\mathbb{G}}\left\{\EE_\pi \left[\frac{1}{T_1}\sum_{t=1}^{T_1}R_{t} \, \,\middle\vert\, \, S_1 \right]\right\} = \int \EE_\pi \left[\frac{1}{T_1}\sum_{t=1}^{T_1}R_{t} \, \,\middle\vert\, \, S_1=s \right] d\GG(s),
%\end{align}
where $\GG$ could be different from the initial distribution $\nu$.  Note that an policy that maximizes $\eta^\pi_{T}(\GG)$ over $\Pi$ may not be optimal if the reference/initial distribution or the time horizon $T$ is changed when implementing in the future. By letting $T$ goes to infinity, we have the long-term average reward for each policy, denoted by $\eta^\pi$.  Through this paper, we assume that for any $\pi \in \Pi$, the induced Markov chain by $P^\pi$ is positive Harris and aperiodic. In this case, $\eta^\pi$ always exists and is independent of the reference distribution. See Theorem 13.3.3 of \cite{meyn2012markov} for more details,  and Sections 5 and 9 of \citep{meyn2012markov} for the definition of positive Harris and aperiodic. 

Next we introduce average visitation density, which motivates our proposed robust criterion.  Define the average visitation density induced by the policy $\pi$ and the initial distribution $\GG$ up to the decision time $t$ as
\[
\bar d^\pi_{t;\GG}(s)\triangleq \frac{1}{t}\sum_{j = 1}^{t}d^{\pi}_{j, \GG} \left(s\right), %\peng{\text{this notation is too much I think.. We use $d$ and $p$.}}
\]
where each $d^{\pi}_{j, \GG}$ is the marginal probability density of $S_t$ induced by $\pi$  and $\bar d^\pi_{1;\GG} = \GG$. Similarly we define $\bar d^D_{T_0;{\nu}}$ as the average visitation density across the decision points in the trajectory $D$ of length $T_0$ with the initial distribution $\nu$. In addition, let $d^\pi$ be the stationary density induced by the policy $\pi$.  Through this paper, for every policy $\pi \in \Pi$ and $t\geq1$, we assume $\bar d^\pi_{t;\GG}\ll d^\pi$, i.e., $\bar d^\pi_{t;\GG}$ are absolutely continuous with respect to $d^\pi$, to avoid some technical difficulties. Finally, we remark that we can rewrite $\eta_{T_0}^\pi(\GG)$ and $\eta^\pi$ as $\int_{s \in S} \mathcal{R}(s) \bar d^\pi_{T_0;\GG}(s)ds$ and $\int_{s \in S} \mathcal{R}(s) d^\pi(s)ds$ respectively.

%For example, if we deploy this optimal policy for the future use with length and reference distribution different from $T_1$ and $\GG$ respectively, it may not be able to guarantee the improved performance compared with other policies. 
%see an analytical example in Appendix \peng{ABCD?} 

%\QZL{Good} \peng{I think same holds for model-based method, so I remove  ``using model-free methods''. Plus we did not discuss what model-free means}

\section{A Robust Average Reward Criterion}
In this section, we introduce a new robust average reward criterion, which will be used to learn an optimal policy for improving the performance of decision making in unseen scenarios such as unknown length of horizon and initial state distribution. We first introduce the uncertainty set:
\[
\Lambda^\pi_c \triangleq \left\{u \in \Lambda(\S) \, | \, \left\|u(\bullet) - d^\pi(\bullet) \right\|_{\TV} \leq c, u \ll d^\pi    \right\},
\]
where $\Lambda(\S)$ is the class of probability measures over the state space $\S$, $\| \bullet \|_{\TV}$ denotes the total variation distance between two probability measures, and $c \in [0, 1]$ is a constant that controls the size of $\Lambda^\pi_c$. Next, we consider a set of average rewards:\[
\U^\pi_c \triangleq \left\{\EE_u\left[\RR(S)\right] \, | \, u \in \Lambda^\pi_c   \right\}.
\]
Recall that $\EE_u$ is the expectation with respect to a probability measure $u$ over the state space $\S$. 
Basically $\U^\pi_c$  represents average rewards with respect to a probability ball centered at the stationary distribution $d^\pi$. The key observation is that $\U^\pi_c$ contains average rewards over different lengths of horizons and initial state distribution, which is essential in achieving our aforementioned goal. To see this, the ergodicity implies that for every $\GG \in \Lambda(\S)$, $\left\| \bar d^\pi_{t;\GG}(\bullet) - d^\pi(\bullet) \right\|_{\TV}  \rightarrow 0$ as $t \rightarrow \infty$. As a result, for any $c$, there must exist a $T$ such that for every $T_1 \geq T$, $\left\|\bar d^\pi_{T_1;\GG}(\bullet) - d^\pi(\bullet) \right\|_{\TV} \leq c$.
Therefore $\eta^\pi_{T_1}(\GG) = \int_{s \in S} \mathcal{R}(s) \bar d^\pi_{T_1;\GG}(s)ds$ must belong to $\U^\pi_c$ for $T_1 \geq T$ by the remark at the end of Section 2. Moreover, $\U^\pi_c$ also quantifies the uncertainty where there is some distributional perturbation on the underlying dynamics induced by the policy. %We remark that the requirement of $u \ll d^\pi$ is to ensure elements in $\U^\pi_c$ can be identified by $d^\pi$. 

Based on the appealing properties of $\U^\pi_c$, it is thus desirable to use it to quantify the uncertain performance of each policy $\pi$ deployed in practice. To protect against such uncertainty, we propose to use the smallest value of $\U^\pi_c$ to evaluate a policy $\pi$, i.e.,
\begin{align}\label{robust measure}
\min_{u \in \Lambda^\pi_c} \, \, \EE_u\left[\RR(S)\right].
\end{align}
Then the optimal robust policy with respect to \eqref{robust measure} in the policy class $\Pi$ is defined as
\begin{align}\label{robust policy}
\pi^\ast_c \in \max_{\pi \in \Pi} \min_{u \in \Lambda^\pi_c} \, \, \EE_u\left[\RR(S)\right],
\end{align}
i.e., the policy that maximizes the worst-case average rewards with respect to the probability ball $\Lambda^\pi_c$. Hence, if $\pi^\ast_c$ is deployed in practice, we can guarantee that the worst-case performance (in terms of average reward) against the probability uncertain set $\Lambda^\pi_c$ is the best, which enhances the generalizability of the learned policy.  

The constant $c$ controls the robust level of $\pi^\ast_c$. When $c = 0$, it degenerates to $\pi^\ast_0 \in \argmax_{\pi \in \Pi} \eta^\pi$, which is an in-class optimal policy with respect to the long-term average reward. When $c = 1$, $\Lambda^\pi_c = \Lambda(\S)$, i.e., the class of all probability distributions. Then $\pi^\ast_1$ can be any policy in $\Pi$ since \eqref{robust measure} is the same for every policy. The larger $c$ is, the more near-term rewards are considered for the policy optimization. In contrast, smaller $c$ weighs more on distant rewards. Therefore, the constant $c$ balances the short-term and long-term effect we consider when finding a robust optimal policy. We provide some insights on how to choose $c$ in the Appendix.
\section{Efficient Statistical Estimation}
In this section, we discuss how to estimate $\pi^\ast_c$ in $\eqref{robust policy}$ given a batch data $\D_n$. Specifically, in Section \ref{sec: duality}, we make use of the convex duality theory to formulate problem $\eqref{robust policy}$ as an M-estimation problem. Based on this result, by leveraging semi-parametric statistics, we show how to efficiently estimate the objective function of our policy optimization problem in Section \ref{sec: EIF}. The related nuisance functions estimation is discussed in Section \ref{sec: nuisance function estimation}. Lastly, we present our overall policy optimization procedure in Section \ref{sec: algorithm}.  Throughout this section, we fix the constant $c$ and use the following notations. %For an arbitrary set $\X$, let $\P(\X)$ be the probability distribution on $\X$.  For any function $f(s, a, s')$ of a transition sample $(s, a, s')$ and distribution $w \in \P(\S \times \A \times \S)$, denote the $L_2(w)$ norm by $\norm{f}^2_w = \int f^2(s, a, s') dw(s, a, s')$. If the norm does not have a subscript, then the expectation is taken with respect to the data generating distribution, that is,  $\norm{f}^2 = \EE[(1/T) \sum_{t=1}^{T}f^2(S_t, A_t, S_{t+1})]$.  
For any function of the trajectory $f(D)$, the sample average is denoted by $\Pn f(D) = (1/n)\sum_{i=1}^{n} f(D_i)$. A transition tuple is either denoted by $Z =(S, A, S')$ or $Z_t =(S_t, A_t, S_{t+1})$ at time $t$. We let $N= nT_0$. 
%the average over samples and decision times is denoted by $\PN f(Z) = \frac{1}{N}\sum_{i=1}^n\sum_{t=1}^{T_0}f(S_{it}, A_{it}, S_{i(t+1)})$.  

%\peng{Suggest we don't use the $Z$ notation given we have so many notations already. I don't see we use this?\QZL{We use several times. You will see}} \peng{We use $P_Nf(Z)$?} \QZL{We use $P_N$ in our proof}
\subsection{Dual Reformulation}\label{sec: duality}
We first reformulate problem \eqref{robust policy} by using the convex duality theory. Define a function $\phi(x) \triangleq \frac{1}{2}| x -1|$ for $x \geq 0$, and $\phi(x) := +\infty$ for $x < 0$. Then by the definition of total variation distance, we can rewrite the set $\Lambda^\pi_c$ as
\begin{align}\label{dual feasible set}
\Lambda^\pi_c = \left\{u \in \Lambda(\S)  \, \Big| \, \EE_{d^\pi}\left[\phi\left(\frac{u(S)}{d^\pi(S)}\right)\right] \leq c, ~ u \ll d^\pi    \right\},
\end{align}
where $\EE_{d^\pi}$ denotes the expectation with respect to the stationary distribution $d^\pi$ over $\S$. By the change of variable (i.e., let $W(s) =\frac{u(s)}{d^\pi(s)}$), we consider a set defined as
\begin{align}\label{dual feasible set 2}
\W^\pi_c = \left\{W \in L^1(\S, \B(\S), d^\pi) \, | \, \EE_{d^\pi}\left[\phi\left(W(S)\right)\right] \leq c, W(s) \geq 0,\, \mbox{for every $s \in \S$}, \,  \EE_{d^\pi}[W(S)] = 1    \right\},
\end{align}
where $L^1(\S, \B(\S), d^\pi) $ is $L^1$ space defined on the measure space $(\S, \B(\S), d^\pi)$.
Using $\W^\pi_c$, we can rewrite our problem \eqref{robust policy} as
\begin{align}\label{robust policy (2)}
    \max_{\pi \in \Pi} \min_{W \in \W^\pi_c} \, \, \EE_{d^\pi}\left[W(S)\RR(S)\right],
\end{align}
where $W(s)$ can be interpreted as a likelihood ratio of $\frac{u(s)}{d^\pi(s)}$ for every $u \in \Lambda^\pi_c$. Define $R_{\min} = \inf_{s \in \S}\RR(s)$. Now we present our first key theorem.
\begin{theorem}\label{duality}
Assume that for every $\pi \in \Pi$, the essential infimum of $\RR$ under $d^\pi$ is $R_{\min}$. Then the following holds:
\begin{align}\label{optimal value equivalent}
    \min_{W \in \W^\pi_c} \, \, \EE_{d^\pi}\left[W(S)\RR(S)\right] = cR_{\min} + \left(1-c\right)\max_{\beta \in \mathbb{R}} \left\{\beta -\frac{1}{\left(1-c\right)}\EE_{d^\pi}\left[\left(-\RR(S)+\beta \right)_+\right]\right\},
\end{align}
\begin{align}\label{optimal solution equivalent}
    &\argmax_{\pi \in \Pi} \min_{W \in \W^\pi_c}\EE_{d^\pi}\left[W(S)\RR(S)\right]=\argmax_{\pi \in \Pi} \max_{\beta \in \mathbb{R}} \left\{\beta -\frac{1}{\left(1-c\right)}\EE_{d^\pi}\left[\left(-\RR(S)+\beta \right)_+\right]\right\}.
\end{align}
\end{theorem}

	Theorem \ref{duality} transforms the max-min problem \eqref{robust policy (2)} into an M-estimation problem using the convex duality theory. Such result, adapted from \citep{shapiro2017distributionally}, makes problem \eqref{robust policy} more computational tractable. In the original formulation, the constraint set (e.g., total variation distance) is very difficult to compute/estimate since the stationary distribution $d^\pi$ is not directly observed and needs to be updated along with $\pi$ during the policy optimization procedure. By transforming into an M-estimation problem, we avoid solving a constraint max-min problem. Furthermore, while $d^\pi$ is still not observed, we can leverage semi-parametric statistics to estimate the objective function more directly, under which the computation can be relatively easy to perform. See following sections for more details. 
	
	Interestingly, maximizing the objective function in the RHS of Equation \eqref{optimal solution equivalent} with respect to $\beta$ is equivalent to computing the $(1-c)$-conditional value-at-risk ($(1-c)$-CVaR) of the reward under the stationary distribution induced by the policy $\pi$ \citep{ben1986expected,rockafellar2000optimization}. CVaR is a coherent risk measure \citep{artzner1999coherent}, frequently used in the domain of finance and engineering. The original CVaR is defined as the truncated mean of some loss above a certain quantile \citep{rockafellar2000optimization}. Here we use $(1-c)$-CVaR to represent the truncated mean of the reward lower than a $(1-c)$-quantile to align with the reward instead of the loss in our problem. One maximizer $\beta^\ast$ (the leftmost of the optimal solution set) in \eqref{optimal solution equivalent} is the corresponding $(1-c)$-quantile of the reward with respect to the stationary distribution $d^\pi$. Since rewards are uniformly bounded, we can show that $|\beta^\ast| \leq R_{\max}$. Therefore it is enough to restrict $\beta$ to be between $-R_{\max}$ and $R_{\max}$. That is, we can obtain $\pi^\ast_c$ by jointly solving
\begin{align}\label{equivalent maximization problem}
    \max_{\pi \in \Pi, |\beta| \leq R_{\max}} \left\{M(\beta, \pi) \triangleq \beta -\frac{1}{\left(1-c\right)}\EE_{d^\pi}\left[\left(-\RR(S)+\beta \right)_+\right]\right\}.
\end{align}

When there is no temporal dependence among the trajectory $D$, the overall problem becomes a single-stage decision making problem which has been extensively studied in the literature. See a review paper by \citep{kosorok2019precision}. In this case, problem \eqref{equivalent maximization problem} degenerates to the policy learning under the CVaR criterion, which was recently studied by \citep{qi2019estimating}. Compared with the single-stage problem, one notable challenge in our problem is that data are not generated by $d^\pi$ and thus the objective function in \eqref{equivalent maximization problem} cannot be directly estimated by the sample-average approximation. We need to leverage the Markov and stationarity assumptions to estimate the objective function in \eqref{equivalent maximization problem} so that the long term effect of the policy is captured. To the best of our knowledge, such robust formulation has not been studied in the literature.

\subsection{A Statistically Efficient Evaluation Method}\label{sec: EIF}
%\peng{We need to decide if we want to keep the dependence on $\beta$. Right now it is quite confusing. Sometimes we have $\beta$ but sometimes we don't. I suggest we keep $\beta$ everywhere. I fix a lot of this already. }

To estimate $\pi^\ast_c$, given limited batch data $\D_n$, we need to first develop an efficient estimator to evaluate the objective function in \eqref{equivalent maximization problem} for any given $\beta$ and $\pi$, after which we can optimize the objective function. 
%Note that $M(\beta, \pi)$ can be viewed as the long-term average reward under modified reward function $\beta - \frac{1}{1-c}(\beta - \RR)_+$. (I don't think they need to know here)
%In this subsection, we develop an estimator for $M(\beta, \pi)$ which satisfies the doubly robust property and achieves the statistical efficiency bound; see Section 4.3 for the details. 
The M-estimation formulation of \eqref{equivalent maximization problem} motivates us to leverage semi-parametric statistics (e.g., \citep{tsiatis2007semiparametric}) to construct an efficient estimator.
Before we introduce our estimator of $M(\beta, \pi)$, we take a detour and consider the following two alternative estimators, which motivate ours.

It can be seen that $M(\beta, \pi)$ is the long-term average reward under a modified reward function: $\beta - \frac{1}{1-c}(\beta - \RR)_+$.  Then one can construct an estimator based on the relative value function of the modified reward. For any given policy $\pi$ and $\beta$, the relative value function (e.g., \citep{hernandez2012further}) can be defined as
\begin{align}
Q^{\pi, \beta}(s, a) := 
\operatorname*{lim}_{t^* \goes \infty} \frac{1}{t^*}\sum_{t=1}^{t^*} \EE_\pi\left[\sum_{k=1}^{t} \left\{\beta - \frac{1}{1-c}\left(\beta - R_{k}\right)_+ - M(\beta, \pi)\right\} \Given S_1 = s, A_1 = a\right],
\label{relative value function}
\end{align}
which we assumed is always well defined.
%This definition is different from the state-action value function used in the discounted sum of rewards setting because the undiscounted sum of rewards in general cannot be bounded \citep{mahadevan1996average}. The definition of \eqref{relative value function} relies on the Poisson equation associated with the long-term average reward criteria \citep{hernandez2012further}.
The Bellman equation related to the relative value function is
\begin{align}
\EE\left[\beta - \frac{1}{1-c}\left(\beta - R_{t}\right)_+ + \sum_{a'} \pi(a'|S_{t+1}) Q(S_{t+1}, a)\given S_t = s, A_t = a\right] = Q(s, a) - \eta,
%	\eta = 0. 
\label{Bellman equation Value}
\end{align}
with respect to $\eta$ and $Q$.
As given by Theorem 7.5.7 of \cite{hernandez2012further}, solving the above equation \eqref{Bellman equation Value} with respect to $(\eta, Q)$ gives us the unique solution $M(\beta, \pi)$, and $Q^{\pi, \beta}$ up to some constant respectively. Therefore, based on the estimating equation \eqref{Bellman equation Value}, one can construct estimators for both $Q^{\pi, \beta}$ and $M(\beta, \pi)$ by using the generalized method of moments \citep{hansen1982large}. This method requires to model $Q^{\pi, \beta}$. If we impose some parametric model on $Q^{\pi, \beta}$, we may suffer from model mis-specification, thus causing biases for estimating $M(\beta, \pi)$. Alternatively, if a nonparametric model is used for $Q^{\pi, \beta}$, while it could be consistent, the resulting estimator for $M(\beta, \pi)$ or regret of the learned policy may not be rate-optimal, say $\sqrt{N}$-consistent. Before we discuss the second estimator for $M(\beta, \pi)$,  define the relative value difference function, which will be used later, as
\begin{align}
	U^{\pi,\beta}(s, a, s') := \sum_{a' \in \A} \pi(a'|s') Q^{\pi, \beta}(s', a') - Q^{\pi, \beta}(s, a), \label{Upi}
	\end{align}
	where $(s, a,s')$ is a transition tuple.

The second estimator of $M(\beta, \pi)$ can be constructed by adjusting the mismatch between the data generating mechanism by the behavior policy and the stationary distribution of a given policy $\pi$. This is motivated by recently proposed marginal importance sampling \citep{liu2018breaking}. Note that
\begin{align}
M(\beta, \pi) &= \int_{s \in S, a \in \A} d^\pi(s) \pi(a|s)\left(\beta - \frac{1}{1-c}\left(-\RR(s)+ \beta\right)\right) ds da \notag \\[0.1in]
& = \EE\left[\frac{1}{T_0}\sum_{t=1}^{T_0}\frac{d^\pi(S_t) \pi(A_t|S_t)}{\bar d^D_{T_0;{\nu}}(S_t, A_t)}\left(\beta - \frac{1}{1-c}\left(-\RR(S_t) + \beta\right)\right)\right] \label{average reward using ratio}.
\end{align}
Recall that $\bar d^D_{T_0;{\nu}}$ is the average visitation density in the batch data $D$. Based on this observation, we can first estimate a ratio function defined as
	\begin{align}
	\omega^\pi(s, a) = \frac{d^\pi(s) \pi(a|s)}{\bar d^D_{T_0;{\nu}}(s, a)}, \label{ratio}
	\end{align}
	after which we can use the sample-average approximation of \eqref{average reward using ratio} and plug in the estimator of ratio function to estimate $M(\beta, \pi)$. The sample-average approximation procedure is valid because the expectation in $\eqref{average reward using ratio}$ is with respect to the data generating process. However, using such an estimator has the same issue as the first one. 
	
 %This can be interpreted as importance weights, similar to those in the average treatment effect estimation \citep{precup2000eligibility}. However, in addition to estimating the behavior policy from the data in the classic importance weights, we also need to correct the discrepancy between the state distribution under the behavior and target policies since policies can not only affect the selection of actions but also the transition of states.  

Towards that end, we combine these two estimators together and introduce an estimator of $M(\beta, \pi)$ that enjoys doubly robust property for model mis-specification and meanwhile achieves statistical efficiency bound, which is the best one can hope for; See the discussion of double robustness and statistical efficiency bound in Section 5.3. Our proposed estimator is inspired by \citep{liao2020batch} and relies on two nuisance functions: one is the relative value difference $U^{\pi, \beta}$ and the other is the ratio function $\omega^\pi$ defined above.  Such estimator is derived from the efficient influence function (EIF) \citep{newey1990semiparametric} of $M(\beta, \pi)$ given as
\begin{align}\label{estimating equation}
	%\phi(Z; U, \omega) =  {\omega}(S, A) \left[\beta - \frac{1}{1-c}\left(\beta - \RR(S) \right)_+ +  U(S, A, S')-  M(\beta, \pi) \right].
    (1/T_0)\sum_{t=1}^{T_0}{\omega}^\pi(S, A) \left[\beta - \frac{1}{1-c}\left(\beta - \RR(S) \right)_+ +  U^{\pi, \beta}(S, A, S')-  \eta \right]. 
\end{align}
%One can show that $\EE[\frac{1}{T_0}\sum_{t=1}^{T_0}\phi(Z_t; U^{\pi, \beta}, \omega^\pi)] = 0$.
One can show that the expectation of the above EIF is zero if and only if $\eta = M(\pi, \beta)$ for any $\pi$ and $\beta$, which naturally forms an estimating equation. 
Based on this, we can first construct estimators for two nuisance functions $U^{\pi, \beta}$ and $\omega^\pi$, denoted by $\widehat U_N^{\pi, \beta}$ and $\widehat \omega_N^\pi$, and then estimate $M(\beta, \pi)$ by solving the empirical version of the plug-in estimating equation, or equivalently
%where $\widehat U_N^{\pi, \beta}$ and $\widehat \omega_N^\pi$ are estimators of $U^\pi$ and $\omega^\pi$ respectively. 
%Specifically, we solve
%$$\Pn \dkh{(1/T_0)\sum_{t=1}^{T_0} \widehat \omega^{\pi}_N(S_{t}, A_{t}) [\beta - \frac{1}{1-c}\left(\beta - R_{t}\right)_+ +\widehat U^{\pi, \beta}_N(S_t, A_t, S_{t+1})  - M(\beta, \pi)]} = 0,$$
%which gives the proposed estimator as
	\begin{align}
	\widehat{M}_N(\beta, \pi) = \frac{\Pn \dkh{(1/T_0)\sum_{t=1}^{T_0} \widehat \omega^{\pi}_N(S_{t}, A_{t}) [\beta - \frac{1}{1-c}\left(\beta - R_{t}\right)_+ + \widehat U^{\pi, \beta}_N(S_t, A_t, S_{t+1})]} }{\Pn \dkh{(1/T_0)\sum_{t=1}^{T_0} \widehat \omega^{\pi}_N(S_{t}, A_{t})} }. \label{pol.eval}
	\end{align}
In Section 5.3, we demonstrate that under some technical assumptions, the proposed estimator $\widehat{M}_N(\beta, \pi)$ has the doubly robust property and achieves statistical efficiency bound, i.e., the supermum of Cramer-Rao low bounds for all parametric submodels that contain the true parameter, using the same notion in \citep{kallus2019efficiently}.

%\peng{this section needs to re-order.. We should put the discussion of estimating $M(\beta, \pi)$ in the front and then introduce the nusiance functions. Otherwise the reader have no idea why we need to consider the nuisance functions.  Also did we tell reader how we estimate the optimal policy? This is the purpose of this section..}
\subsection{Nuisance Functions Estimation}\label{sec: nuisance function estimation}
%\peng{Put it into the appendix? Just mention we use the coupled estimation bla bla.}
The doubly robust structure of our estimator has a weak requirement on the convergence rate of each nuisance function estimation for achieving the optimal convergence rate to the targeted parameter $M(\beta, \pi)$. This promotes the use of nonparametric models for estimating these nuisance functions. In the following, we briefly discuss how to nonparametrically estimate the relative value difference function and the ratio function. %The details are given in the Appendix.

\paragraph{Estimation of relative value difference function.}  We use the Bellman equation given in \eqref{Bellman equation Value} to estimate the nuisance function $U^{\pi, \beta}$ via estimating $Q^{\pi, \beta}$. Recall that $Z_t = (S_t, A_t, S_{t+1})$ be the transition sample at time $t$ and define the so-called temporal difference (TD) error as
$$
\delta^{\pi, \beta}(Z_t; \eta, Q) = \beta - \frac{1}{1-c}\left(\beta - R_{t}\right)_+ + \sum_{a'} \pi(a'|S_{t+1}) Q(S_{t+1}, a) - Q(S_t, A_t) - \eta. 
$$
%\peng{where is $\beta$ in LHS?}
%We will omit the subscript $\beta$ in $\delta^{\pi, \beta}$ when there is no confusion. 
As a result of the Bellman equation \eqref{Bellman equation Value}, we can rewrite $(M(\beta, \pi), Q^{\pi, \beta})$ as an optimal solution of the following optimization problem.
	\begin{align}\label{population U}
	(M(\beta, \pi), Q^{\pi, \beta}) \in \argmin_{\eta \in \mathbb{R}, Q} \EE\left[\frac{1}{T_0}\sum_{t=1}^{T_0}\left(\EE[\delta^{\pi, \beta}(Z_t; \eta, Q)|S_t, A_t]\right)^2\right]
	\end{align}
The above Bellman equation can only identify the relative value function $Q^{\pi, \beta}$ up to a constant (\cite{hernandez2012further}).
%That is, the set of solutions of (\ref{Bellman equation Value}) is given by $\{(M(\beta, \pi), Q): Q = Q^{\pi, \beta} + C \ones, C \in \R\}$ where $\ones(s, a) = 1$ for all $(s, a)$;  see \cite{puterman1994markov}, p. 343 for details. 
Fortunately, since our goal is to estimate $U^{\pi, \beta}$, estimating one specific version of $Q^{\pi, \beta}$ is enough. For example, one can impose one  restriction on $Q^{\pi, \beta}$ to make it identifiable.   Define a shifted relative value function by $\tilde Q^{\pi, \beta}(s, a) = Q^{\pi, \beta}(s, a) - Q^{\pi, \beta}(s^*, a^*)$ for an arbitrarily chosen state-action pair $(s^*, a^*) \in \S \times \A$.  By restricting to $Q(s^*, a^*) =0$, the solution of Bellman equations (\ref{Bellman equation Value}) is unique and given as $(M(\beta, \pi), \tilde Q^{\pi, \beta})$. For the ease of notation, we will use $\widehat Q^{\pi, \beta}_N$ to denote the estimator of the shifted value function $\tilde Q^{\pi, \beta}$.

We know that $\tilde Q^{\pi, \beta}$ can be characterized as the minimizer of the above objective function \eqref{population U} which involves the conditional expectation of a function inside. 
%One may try to directly minimize the empirical version of it and obtain the corresponding estimator. However, this approach in general introduces non-negligible bias. See \cite{antos_learning_2008} for more details. 
Borrowing ideas from \citep{farahmand2016regularized,liao2019off}, we first estimate the projection of $\delta^{\pi, \beta}(Z_t; \eta, Q)$ onto the space of $(S_t, A_t)$, after which we optimize the empirical version of the above optimization problem. Define $\F_1$ and $\G_1$ as two specific classes of functions over the state-action space, where we use $\F_1$ to model the shifted relative value function $\tilde Q^{\pi, \beta}$ and thus require $f (s^*, a^*) = 0$ for all $f \in \F$, and use $\G_1$ to model $\EE[\delta^{\pi, \beta}(Z_t; \eta, Q)|S_t, A_t]$.  In addition, let $J_1: \F_1 \goes \R^+$ and $J_2: \G_1 \goes \R^+$ be two penalty functions that measure the complexities of these two functional classes respectively. Distinct from $\widehat{M}_N(\beta, \pi)$ constructed from the EIF \eqref{estimating equation}, we use $\widehat{\eta}^{\pi,\beta}_N$ to denote the resulting estimator of $M(\beta, \pi)$ obtained from the Bellman equation \eqref{Bellman equation Value}. Therefore given two tuning parameters $\lambda_{1N}$ and $\mu_{1N}$, we can obtain the estimator $(\widehat{\eta}^\pi_N, \widehat Q^{\pi, \beta}_N)$ by minimizing the square of the projected Bellman equation error:
	\begin{align}
	(\widehat \eta_N^{\pi, \beta}, \widehat Q_N^{\pi, \beta}) = \argminb_{(\eta, Q) \in \R \times \F_1} \Pn \left[\frac{1}{T_0} \sum_{t=1}^{T_0} \gn(S_t, A_t; \eta, Q)^2\right] + \lambda_{1N} J_1^2(Q), \label{estimator}
	\end{align}
	where $\gn(\cdot, \cdot; \eta, Q)$ is the projected Bellman error with respect to $(\eta, Q)$, the policy $\pi$ and $\beta$. which is computed by
	\begin{align}
	%\g: (\eta, Q) \goes
	& \gn(\cdot, \cdot; \eta, Q) = \argminb_{g \in \G_1} \Pn\Big[ \frac{1}{T_0}\sum_{t=1}^{T_0} \big( \delta^{\pi, \beta}(Z_t; \eta, Q) - g(S_t, A_t)\big)^2  \Big] + \mu_{1N} J_2^2(g). \hspace{-3ex} \label{gn value}
%	& \text{where } L_n^\pi(g)   \notag
	\end{align}
%	\peng{Do we allow the tuning parameters to depend on the policy??? If not, we need to say that. For the theory we don't need. }
Such an estimator is called the coupled estimator in \cite{liao2020batch}.
	Finally, we can estimate $U^{\pi, \beta}$ by $\widehat U_N^{\pi, \beta} (s, a, s') = \sum_{a'} \pi(a'|s') \widehat Q^{\pi, \beta}_N(s', a') -  \widehat Q_N^{\pi, \beta}(s, a)$ for any $(s, a, s')$.

	\paragraph{Estimation of the ratio function.} Next we use another coupled estimator proposed by \citep{liao2019off} to estimate the ratio function $\omega^\pi$. This can be achieved by first estimating $e^\pi$, a scaled version of the ratio function defined as
	%It has been shown that a scaled version of ratio function can be identified as a minimizer of some objective function, based on which we construct the estimator. 
	%Specifically, we first estimate $e^\pi$, a scaled version of ratio function defined as follows.
	\begin{align}
%	e^\pi(s, a) =  \frac{\omega^\pi(s, a)}{1+\sigma_\pi^2}, \label{epi defn}
	e^\pi(s, a) =  \frac{\omega^\pi(s, a)}{\int \omega^\pi(s, a) d^\pi(s) \pi(a|s)dsda}\label{epi defn}.
	\end{align}
By treating $e^\pi$ as a new reward function, we can see that the long-term average reward is 1 under the induced Markov chain. Based on this, define a ``new" relative value function $H^{\pi} (s, a) = 	\operatorname*{lim}_{t^* \goes \infty} \frac{1}{t^*}\sum_{t=1}^{t^*} \EE_\pi\left[\sum_{k=1}^{t} \left\{1 -e^\pi(S_k, A_k)\right\} \Given S_1 = s, A_1 = a\right],$ which we assume is well defined, and a ``new" temporal difference as
	 $\Delta^\pi(Z_t; H) = 1 - H(S_t, A_t) + \sum_{a'} \pi(a'|S_{t+1})  H(S_{t+1}, a')$, where $H$ is an arbitrary function over $\S \times \A$. 
	\begin{comment}
	Analogous to \eqref{Bellman equation Value} in the estimation of the relative value difference function, we can show the following bellman equation holds.
	\begin{align}
	e^\pi (s, a) = \EE[\Delta^\pi(Z_t; H^\pi)|S_t=s, A_t=a]. \label{epi qpi}
	\end{align}
	Note that we cannot use the same strategy as that in the first nuisance estimation because we do not know $H^\pi$ and $e^\pi$. Relying on the invariant property of the stationary distribution, it can be shown that for any state-action function $H$,
	\begin{align*}
	\EE\left[\frac{1}{T_0}\sum_{t = 1}^{T_0}\left(\sum_{a'} \EE\left[\pi(a'|S_{t+1})  H(S_{t+1}, a')\, | \, S_t, A_t\right] - H(S_t, A_t)\right)e^\pi (S_t, A_t)   \right] = 0.
	\end{align*}
	%which implies that
	%\begin{align*}
	%\EE\left[\frac{1}{T_0}\sum_{t = 1}^{T_0}\left(\sum_{a'} \EE\left[\pi(a'|S_{t+1})  H(S_{t+1}, a')\, | \, S_t, A_t\right] - H(S_t, A_t)\right)\EE\left[\Delta^\pi(Z_t; H^\pi)|S_t, A_t\right]   \right] = 0,
	%\end{align*}
	\end{comment}
    It can be seen that $e^\pi (s, a) = \EE[\Delta^\pi(Z_t; H^\pi)|S_t=s, A_t=a]$. Relying on the invariant property of the stationary distribution $d^\pi$, one can show that $H^\pi$  satisfies: 
	\begin{align}
	H^\pi \in \argmin_{H }\EE\Big[\frac{1}{T_0}\sum_{t=1}^{T_0}\big( \EE[\Delta^\pi(Z_t; H)|S_t, A_t]\big)^2\Big], \label{qpi min}
	\end{align}
based on which we can develop a coupled estimator for $e^\pi$ in the same manner of estimating $U^{\pi, \beta}$. Specifically, define a function class $\F_2$ over $\S \times \A$ satisfying that $f (s^*, a^*) = 0$ for all $f \in \F_2$ (We can only identify $H^\pi$ up to a constant, so we target on a specific one denoted by $\tilde H^\pi$), and a specific class of functions $\G_2$ over $S \times \A$. 
	 %Note that the operator $\Delta^\pi(Z_t; H)$, similar to the operator $\delta^{\pi, \beta}(Z_t;\eta, Q)$, is invariant up to a constant shift of $H$. Since one $H^\pi$ suffices, in the estimator decribed below we only target on one specific $H^\pi$ with respect to some reference state-action pair, i.e., the shifted version $\tilde H^\pi(\cdot, \cdot) = H^\pi(\cdot, \cdot) - H^\pi(s^*, a^*)$, for a particular $(s^*, a^*)$ and we  assume a function class $\F_2$ over $\S \times \A$ satisfying that $f (s^*, a^*) = 0$ for all $f \in \F_2$. 
	 %, that is, $\Delta^\pi(\cdot; H) = \epsilon^\pi(\cdot; H+C)$ for any constant $C$. Thus by minimizing the above objective function we can only identify up to a constant of $H^\pi$.  This is similar to the case of the relative value function where the Bellman equation \eqref{Bellman equation Value} only identifies the true relative value function up to a constant. Fortunately, one specific $H^\pi$ suffices since the parameter of interest is $e^\pi$.   Consequently, in the estimator below we only target one specific version with respect to some reference state-action pair, i.e., the shifted version $\tilde H^\pi(\cdot, \cdot) = H^\pi(\cdot, \cdot) - H^\pi(s^*, a^*)$, for a particular $(s^*, a^*)$ and we  assume the function class $\F_2$ over $\S \times \A$ satisfying that $f (s^*, a^*) = 0$ for all $f \in \F_2$. 
	%The estimator of $H^\pi$ is constructed based projection of the function $(s, a) \mapsto \EE[\Delta^\pi(Z_t; H)|S_t=s, A_t=a]$ onto the functional space $\G_2$. 
	Then given tuning parameters $\lambda_{2N}$ and $\mu_{2N}$, the estimator  $\widehat H_N^\pi$ can be obtained by minimizing the square of the projected value with respect to $H\in \F_2$: 
	\begin{align}
	\widehat H_N^\pi = \argminb_{H \in \F_2} \Pn\Big[\frac{1}{T_0}\sum_{t=1}^{T_0}  \hn(S_t, A_t; H)^2\Big] + \lambda_{2N} J_1^2(H) \label{qhat}
	\end{align}
	where $\hn(\cdot, \cdot; H)$ is given by
	\begin{align}
	\widehat h_N(\cdot, \cdot; H) = \argminb_{h \in \G_2} \Pn\Big[\frac{1}{T_0} \sum_{t=1}^{T_0} \big(\Delta^\pi(Z_t; H)- h(S_t, A_t)\big)^2\Big]+  \mu_{2N} J_2^2(h).
	\label{minimizer of ratio}
	\end{align}
	%\peng{We use different function class but uses the same regularization functional...}
For the ease of presentation, we use the same penalty functions as that in estimating the relative value difference function. Given the estimator $\widehat H_N^\pi$, we obtain the estimator of $e^\pi$ as $\widehat e^\pi_N = \widehat h_N(\cdot, \cdot; \widehat H_N^\pi)$. 
%	\paragraph{Step 3}  
	By the definition of $\omega^\pi$, we have $\EE[(1/T_0) \sum_{t=1}^{T_0} \omega^\pi(S_t, A_t)] = 1$, which makes us to estimate $\omega^\pi$ by
	%Since $e^\pi$ is a scaled version of $\omega^\pi$ up to a constant, we construct the estimator of ratio $\omega^\pi$ by scaling $\widehat e^\pi_N$ such that its empirical mean is one, that is, 
	\begin{align}
	\widehat \omega^\pi_N(s, a) = \widehat e_N^\pi(s, a) /\Pn[(1/T_0)\sum_{t=1}^{T_0} \widehat e_N^\pi(S_t,A_t)], ~ \forall (s, a) \in \S \times \A.  \label{ratio estimator}
	\end{align}
\subsection{A Statistically Efficient Learning Method}\label{sec: algorithm}
For any given $\pi$ and $\beta$, after obtaining estimators for nuisance functions, we plug them in \eqref{pol.eval} for obtaining our estimator $\widehat{M}_N(\beta, \pi)$ of $M(\beta, \pi)$. The second step is to maximize $\widehat{M}_N(\beta, \pi)$ with respect to $\pi$ and $\beta$ for an estimated policy $\widehat \pi_{N}^c$ of $\pi^\ast_c$, i.e., 
\begin{align}
& \max_{\pi \in \Pi, \beta \in \mathbb{R}} && \frac{\Pn \left\{(1/T_0)\sum_{t=1}^{T_0} \widehat \omega^{\pi}_N(S_{t}, A_{t}) \left[\beta - \frac{1}{1-c}\left(\beta - R_{t}\right)_+ + \widehat U^{\pi, \beta}_N(S_t, A_t, S_{t+1})\right]\right\} }{\Pn \left\{(1/T_0)\sum_{t=1}^{T_0} \widehat \omega^{\pi}_N(S_{t}, A_{t})\right\} },\label{multi-level optimization} 
\end{align}
where $\widehat U^{\pi, \beta}_N$ and $ \widehat \omega^{\pi}_N$ are obtained via \eqref{estimator}-\eqref{gn value} and \eqref{qhat}-\eqref{ratio estimator} respectively.
%Overall, our batch policy learning algorithm consists of two steps. The first step is to obtain an efficient estimation of the objective function $M(\beta, \pi)$ in \eqref{equivalent maximization problem} using $\D_n$ for any policy $\pi$ and $\beta$, over which we search for the best in-class policy in the second step. 
These two steps form a bilevel (multi-level) optimization problem as we need to update two nuisance functions along with the update of $\pi$ in the objective of \eqref{multi-level optimization}. We defer the discussion of computation to Section 6 and the details of our full algorithm can be found in Section 2 of the Supplementary Material.
	\section{Theoretical Results}
	In this section, we provide theoretical justifications for our efficient learning method in estimating $\pi^\ast_c$. In particular, in Section \ref{sec: technical assumption}, we list all related technical assumptions.  In Section \ref{sec: nuisance function bound}, we derive uniform finite sample error bounds of our estimators $\widehat U^{\pi, \beta}_N$ and $\widehat \omega^\pi_N$ for $U^{\pi, \beta}$ and $\omega^\pi$ respectively over $\Pi$ and $[-R_{\max}, R_{\max}]$.  We then show our estimator $\widehat M_N(\beta, \pi)$ has doubly robust property and achieves the statistical efficiency bound in Section \ref{sec: semiparametric}. Finally, we establish a rate-optimal up to a logarithm factor finite sample upper bound on the regret of $\widehat{\pi}_N^c$, which is discussed in Section \ref{sec: regret bound}. All of these asymptotic (or finite-sample) results are derived in terms of the number of trajectories $n$ and the number of decision points $T$, which are novel.
	
	\paragraph{Notations.} Consider a state-action function $f(s, a)$.  Denote the conditional expectation operator by $\P^\pi f: (s, a) \mapsto \EE_\pi[f(S_{t+1}, A_{t+1})|S_t=s, A_t =a].$ Let the expectation under the stationary distribution induced by $\pi$ be $d^\pi(f) = \int f(s, a) d^\pi(s)\pi(a| s)dads$.  For a function $g(s, a, s')$ (or $g(s, a)$),  define $\norm{g}^2 = \EE\left\{(1/T_0) \sum_{t=1}^{T_0}g^2(S_t, A_t, S_{t+1})\right\}$ or ($\EE\left\{(1/T_0) \sum_{t=1}^{T_0}g^2(S_t, A_t)\right\}$). For a set $\X$ and $M > 0$, let $\B(X, M)$ be the class of bounded functions on $\X$ such that $\norm{f}_{\infty} \leq M$.  Denote  by $\N(\epsilon, \F, \text{dist}({\bullet}))$ the $\epsilon$-covering number of a set of functions $\F$, with respect to a certain metric, $\text{dist}({\bullet})$ (The definition of covering number can be found in Section 3 of Supplementary material). In addition, we use $\xrightarrow{d}$ to denote the weak convergence as $N \rightarrow \infty$. Before presenting our theoretical results, we need several technical assumptions stated below.
	
	 \subsection{Technical Assumptions}\label{sec: technical assumption}
	\begin{assumption}\label{ass: stationarity}
		The stochastic process $\{S_t, A_t\}_{t\geq 1}$ induced by the behavior policy $\pi^b$ is a stationary, exponentially $\boldsymbol{\beta}$-mixing stochastic process. The $\boldsymbol{\beta}$-mixing coefficient at time lag $k$ satisfies that $\beta_k \leq \beta_0 \exp(-\beta_1 k)$ for $\beta_0 \geq 0$ and $\beta_1 > 0$. 
		%The behavior policy is stationary depending on the current state, denoted by $\pi^b$.  
		In addition, there exists a positive constant $p_{\min}$ such that the behavior policy induced stationary density $d^{\pi_b}(s, a) \geq p_{\min}$ for every $(s, a) \in \S \times \A$.
		\end{assumption}
Assumption \ref{ass: stationarity}  characterizes the dependency among observations over time. The  $\boldsymbol{\beta}$-mixing coefficient at time lag $k$ basically means that the dependency between $\{S_t, A_t\}_{t \leq j }$	and $\{S_t, A_t\}_{t \geq (j+k) }$ decays to 0 at the exponential rate with respect to $k$. See \cite{bradley2005basic} for the exact definition of the exponentially $\boldsymbol{\beta}$-mixing. Note that these conditions are only imposed on the Markov chain induced by the behavior policy (the observed data) and thus independent of target ones (i.e., $\pi \in \Pi$). Therefore the mixing coefficients are fixed. Indeed, if the induced Markov chain is geometric ergodic and stationary (e.g., finite state, irreducible and aperiodic chains), then  $\{S_t, A_t\}_{t\geq 1}$ is at least exponentially $\boldsymbol{\beta}$-mixing. Furthermore, if we assume the induced Markov chain satisfies uniformly geometric ergodicity, then the process is $\phi$-mixing, which is stronger than $\boldsymbol{\beta}$-mixing. For detailed discussion, we refer to \cite{bradley2005basic}.  %\peng{I still don't see how this statement can hold. Also does the reader need to know about this here?} 
The stationary assumption on $\{S_t, A_t\}_{t\geq 1}$ is commonly assumed in the literature such as \citep{kallus2019efficiently}. In addition, this assumption may be further relaxed to so called asymptotically stationary stochastic processes \citep{agarwal2012generalization}. The generalization bounds related to this have been recently developed by \citep{kuznetsov2017generalization}. Since it is beyond the scope of this paper, we decide to leave it as a future work.  The lower bound requirement on $d^{\pi_b}$ is to make sure the ratio function is well defined and avoid the non-parametric identifiability issue for estimating $M(\beta, \pi)$. This is similar to the strict positivity assumption in causal inference. 
\begin{assumption}\label{ass: policy class}
	The policy class $\Pi$, with some distance metric $d_\Pi(\bullet, \bullet)$, satisfies:
	%\peng{I think we can merge the a-d. there exits constants such that for all (s,a, s'), pi1,pi2, beta1 beta2 bl abla holds }
	\begin{enumerate}[label=(\alph*)]
	\item There exists a positive constant $C_1$ such that for every $(s, a) \in \S \times \A$, and $\pi_1, \pi_2 \in \Pi$, and $\bar {\beta}_1, \bar {\beta}_2 \in [-R_{\max}, R_{\max}]$,
		\begin{align}
		| \pi_1(a | s) -\pi_2(a | s) | & \leq C_1d_\Pi(\pi_1, \pi_2), \label{lip pi}\\
			|\omega^{\pi_1}(s, a) - \omega^{\pi_2}(s, a)| &\leq C_1 d_{\Pi}(\pi_1, \pi_2), \label{lip ratio}\\
				|Q^{\pi_1, \bar \beta_1}(s, a, s') - Q^{\pi_2, \bar \beta_2}(s, a, s')|& \leq C_1\left(d_{\Pi}(\pi_1, \pi_2) + |\bar {\beta}_1 -\bar {\beta}_2|\right),\label{lips for U}\\
					|M(\bar {\beta}_1, \pi_1) - M(\bar \beta_2, \pi_2)|&  \leq C_1\left(d_{\Pi}(\pi_1, \pi_2) + \bar {\beta}_1 -\bar{ {\beta}}_2|\right).\label{lips for M}
		\end{align}
		\begin{comment}
		\item There exist some positive constant $C_1$ such that for every $(s, a) \in \S \times \A$, 
		$
		| \pi_1(a | s) -\pi_2(a | s) | \leq C_1d_\Pi(\pi_1, \pi_2).
		$
		\item There exist some positive constant $C_2$ such that for every $(s, a) \in (\S \times \A)$, and $\pi_1, \pi_2 \in \Pi$,
		\begin{align}
			|\omega^{\pi_1}(s, a) - \omega^{\pi_2}(s, a)| &\leq C_2 d_{\Pi}(\pi_1, \pi_2) \label{lip ratio} 
		\end{align}
		\item There exist some positive constants $C_3$ and $C_4$ such that for every $(s, a, s') \in (\S \times \A \times \S)$, $\pi_1, \pi_2 \in \Pi$, and ${\beta}_1, {\beta}_2 \in [-R_{\max}, R_{\max}]$,
		\begin{align}\label{lips for U}
			|Q^{\pi_1, \beta_1}(s, a, s') - Q^{\pi_2, \beta_2}(s, a, s')| \leq C_3d_{\Pi}(\pi_1, \pi_2) + C_4 |{\beta}_1 -{\beta}_2|
		\end{align}
		\item There exist some positive constants $C_5$ and $C_6$ such that for every $\pi_1, \pi_2 \in \Pi$, and $\bar{\beta}_1, \bar{\beta}_2 \in [-R_{\max}, R_{\max}]$,
		\begin{align}\label{lips for M}
		|M(\bar{\beta}_1, \pi_1) - M(\bar{\beta}_2, \pi_2)| \leq C_5d_{\Pi}(\pi_1, \pi_2) + C_6|\bar{\beta}_1 -\bar{\beta}_2|
		\end{align}
		\end{comment}
		\item There exists a positive  constant $C_2$ such that
		\begin{align}\label{vc entropy}
		\log \N(\epsilon, \Pi, d_{\Pi}) \leq C_2 \VC(\Pi) \log(\frac{1}{\epsilon}),
		\end{align}
		where $\VC(\Pi)$ is some positive index measuring the complexity of $\Pi$.
		\item 
		There exists some positive sequence  $\{\bar \alpha_t \}_{t \geq 1}$, where $\lim_{t \goes \infty}\bar{\alpha}_t = 0$, and a positive constant $C_3$, such that for every $\pi \in \Pi$ and $f$ over $\S \times \A$, the following holds for all $t \geq 1$:
		\begin{align}
		& \norm{(\P^\pi)^t f - d^\pi(f)}  \leq C_3 \norm{f} \bar  \alpha_t. \label{ass}
		\end{align}
		\item $\sup_{\pi \in \Pi} \norm{\omega^\pi}_\infty < \infty$.
	\end{enumerate}
\end{assumption}
Assumption \ref{ass: policy class} imposes structural assumptions on the policy class $\Pi$. In order to quantify the complexity of nuisance functions with respect to $\pi \in \Pi$ and $\beta \in [-R_{\max}, R_{\max}]$, we need to impose Lipschitz properties in Assumption \ref{ass: policy class} (a). The distance metric $d_\Pi$ is associated with the policy class. For example, if we consider a parametrized policy class indexed by $\theta$ (i.e., $\Pi = \{\pi_{\theta}, \theta \in \Theta\}$), then we can let $d_\Pi(\pi_{\theta_1}, \pi_{\theta_2}) = \norm{\theta_1 - \theta_2}_\infty$. If $\pi_\theta \in \Pi$ is Lipschitz continuous with respect to $\theta$, then \eqref{lip pi} is automatically satisfied for bounded state space. 
Moreover, for every $\pi \in \Pi$, if the induced Markov chain is uniformly geometric ergodic, then relying on the sensitivity bound such as \cite[Collary 3.1]{mitrophanov2005sensitivity}, \eqref{lip ratio}-\eqref{lips for M} will hold. Similar results and related proofs can be found in \citep{liao2020batch}.  Assumption \ref{ass: policy class}~(b) imposes an entropy condition on $\Pi$, which is commonly assumed in the finite-horizon settings such as \cite{athey2017efficient}. When we consider $\Pi$ parametrized by $\theta$, this condition can be replaced by restricting $\theta$ in a compact set. Assumption \ref{ass: policy class}~(c), which is mild, is related to the mixing-time of the induced Markov chain $\mathbb{P}^\pi$. Such an assumption is used to quantify the estimation errors of two nuisance functions by $\norm{\bullet}$ metric (i.e., $L^2$ norm with respect to the data generating process), bridging the gap between the target policy $\pi$ and the behavior one. We remark that Assumption \ref{ass: policy class}~(c) is much weaker than those
in \citep{van1998learning,liao2019off}. The last condition of Assumption \ref{ass: policy class} ensures the uniform upper bound for the true ratio function. This requires that all the target policies in $\Pi$ have some uniform overlap with the behavior policy. We believe this can be relaxed to some finite moment conditions by using some concentration inequalities for the suprema of unbounded empirical processes in the dependent data setting. 

We also need several technical assumptions on $(\F_j, \G_j)$ for $j = 1, 2$, which are the function classes used in the estimating the nuisance functions $U^{\pi, \beta}$ and $\omega^\pi$ respectively.

\begin{assumption} 
	\label{assumption: function classes}
	The following conditions are satisfied for $(\F, \G) = (\F_j, \G_j)$ with $j = 1,2$:
	
	\begin{enumerate}[label={(\alph*)}]
		\item  \label{b1} $\F \subset \B(\S \times \A, F_{\max})$ and $\G \subset \B(\S \times \A, G_{\max})$
		
		\item \label{b2} $f(s^*, a^*) = 0, f \in \F$. 
		
		\item  The regularization functionals, $J_1$ and $J_2$, are pseudo norms and induced by the inner products $J_1(\bullet, \bullet)$ and $J_2(\bullet, \bullet)$, respectively. 
		
		\item \label{b6}
		Let $\F_M = \{f \in  \F: J_1(f) \leq M\}$ and $\G_M = \{g \in \G: J_2(g) \leq M\}$. There exists some positive constant $C_4$  and $\alpha \in (0, 1)$ such that for any $\epsilon, M > 0$,
		\begin{align*}
		& \max \bigdkh{\log {\N} (\epsilon, \G_M, \norm{\cdot}_\infty), \log {\N}(\epsilon, \F_M, \norm{\cdot}_\infty)} \leq C_4 \left(\frac{M}{\epsilon}\right)^{2\alpha}.
		\end{align*}
		
	\end{enumerate}

\end{assumption}
We assume that functions in $\F_j$ and $\G_j$ are uniformly bounded to avoid some technical difficulty, while this can be relaxed by some truncation techniques. The requirement of $f(s^*, a^*) = 0$ for all $f \in \F$ is used for  identifying $\tilde{Q}^{\pi, \beta}$ and $\tilde{H}^\pi$. Such requirement does not create difficulty in computing our nuisance function estimators. See Section 2 of Supplementary Material for details. The last two technical conditions measure the complexity of functional classes. Similar assumptions have been used in the literature  such as \citep{farahmand2012regularized} and \citep{steinwart2008support}.
	
	\subsection{Finite Sample Error Bounds for Nuisance Functions}\label{sec: nuisance function bound}
	We first develop the uniform error finite sample bound for the relative value difference function. Define the projected Bellman error 	operator as
	\begin{align}\label{surrogate BelErr Op}
		g^\ast_{\pi, \beta}(\cdot, \cdot; \eta, Q) := \argminb_{g \in \G_1} \EE\left[\frac{1}{T_0} \sum_{t=1}^{T_0} \left\{\delta^{\pi, \beta}(Z_t; \eta, Q)- g(S_t, A_t)\right\}^2 \right].
	\end{align}
	We need the following additional assumptions to obtain the error bound.
		\begin{assumption}
		\label{assumption: value function}
		In the estimation of the relative value difference function, the following conditions are satisfied.
		\begin{enumerate}[label={(\alph*)}]
			
			\item \label{b3_2}
			$\tilde Q^{\pi, \beta} \in \F_1$ for $\pi \in \Pi$ and  $\sup_{\pi \in \Pi, |\beta| \leq R_{\max}} J_1(\tilde Q^{\pi, \beta})  < \infty$. 
			
			\item \label{b3_3} 
			$0 \in \G_1$. 
			
			\item \label{b4_2}
			There exists $\kappa > 0$, such that $\inf\{ \norm{g^\ast_{\pi, \beta}(\cdot, \cdot; \eta, Q)}_{} :   \norm{\EE[\delta^{\pi, \beta}(Z_t; \eta, Q)|S_t=\bullet, A_t=\bullet]}_{ } = 1, \abs{\eta} \leq R_{\max}, \abs{\beta} \leq R_{\max}, Q \in \F_1, \pi \in \Pi \} \geq \kappa. $

			\item \label{b5_2} There exists some positive constant $C_5$ such that $J_2\left\{g^\ast_{\pi, \beta}
			(\cdot, \cdot; \eta, Q)\right\} \leq C_5( 1 + J_1(Q))$ holds for all $\beta, \eta \in \R$, $Q \in \F_1$ and $\pi \in \Pi$. 
			
		\end{enumerate}	
	\end{assumption}
Then we have the following theorem that gives the finite sample error bound of our estimator for the relative value difference function.
		\begin{theorem}
		\label{thm: value}
		Suppose the tuning parameters $\mu_{1N} \simeq \lambda_{1N} \simeq (1+\VC(\Pi))(\log N)^{\frac{2+\alpha}{1+\alpha}}N^{-\frac{1}{1+\alpha}}$ and Assumptions   \ref{ass: stationarity}-\ref{assumption: value function} hold.  
		Then there exists some positive constant $C_{6}$ such that, for sufficiently large $N$, the following holds with  probability at least {$1-\frac{1}{N}$:}
		\begin{align*}
		\sup_{\pi \in \Pi, |\beta| \leq R_{\max}} \norm{\widehat U^{\pi, \beta}_N- \Upi}^2 \leq C_{6} (1+\VC(\Pi))^{\frac{1}{1+\alpha}}\log(N)^{\frac{2+\alpha}{1+\alpha}}N^{-\frac{1}{1+\alpha}},
		\end{align*}
		where the constant $C_6$ depends on $\beta_0, \beta_1, p_{\min}, \sup_{\pi \in \Pi} \norm{\omega^\pi}_\infty,\sup_{\pi \in \Pi, |\beta| \leq R_{\max}} J_1(\tilde Q^{\pi, \beta}), R_{max},  \kappa, F_{\max}, G_{\max}$ and constants $C_1$ to $C_5$.
	\end{theorem}
\begin{remark}
	Assumption \ref{assumption: value function}(a) assumes $\F_1$ contains true $\tilde{Q}^{\pi,\beta}$ and the penalty term is uniformly bounded. Assumption \ref{assumption: value function}(b)-(c) basically assume that the projected Bellman error is able to identify the true $M(\beta, \pi)$ and $Q^{\pi,\beta}$. Theorem \ref{thm: value} generalizes the results in \citep{liao2020batch} by deriving the finite-sample error bound in terms of both the sample size $n$ and the number of decision points $T_0$ in each trajectory. This error bound indicates that the estimator of the relative value difference function is consistent as long as either $n$ or $T_0$ goes to infinity. More importantly, our error bound can achieve the optimal rate  $N^{-\frac{1}{1+\alpha}}$ in the classical setting of nonparametric regression  up to a logarithm factor \citep{stone1982optimal}. The additional term $(1+\VC(\Pi))$ appears because this error bound is established uniformly over $\beta \in [-R_{\max}, R_{\max}]$ and $\pi \in \Pi$. Our proof uses the independent block techniques from \citep{yu1994rates} and is inspired by proof techniques in \citep{gyorfi2006distribution,farahmand2012regularized,liao2019off,liao2020batch}. 
\end{remark}
	
	Next, we discuss the uniform finite sample error bound for the ratio function $\omega^\pi$.   For $\pi \in \Pi$ and $H \in \F_2$, define the projected error as
	 $$ \h(\cdot, \cdot; H) = \argminb_{h \in \G_2} \EE\left[\frac{1}{T_0} \sum_{t=1}^{T_0} \left\{\Delta^\pi(Z_t; H) - h(S_t, A_t)\right\}^2 \right].$$  
	To derive the error bound, we need the following conditions similar as Assumption \ref{assumption: value function}.
	\begin{assumption}
		\label{assumption: ratio function}
		We assume that
		
		\begin{enumerate}
			[label={(\alph*)}]
			
			\item  \label{b3}
			For $\pi \in \Pi$, $\tilde H^{\pi}(\cdot, \cdot) \in \F_2$, and $\sup_{\pi \in \Pi} J_1(\tilde H^{\pi}) < \infty$. 
			
			\item \label{b3-2} $e^\pi \in \G_2$, for every $\pi \in \Pi$. 
			
			\item \label{b4}
			
			There exits $\kappa' > 0$, such that $\inf \{\norm{\g(\cdot, \cdot; H) - \g(\cdot, \cdot; \tilde H^{\pi}) }_{} :  \norm{(\I - \P^\pi)(H - \tilde H^{\pi})}_{ } = 1, H \in \F, \pi \in \Pi\}\geq \kappa'$.
			
			\item \label{b5}
			
			There exists some constant $C_{7}$ such that $J_2\left\{\h
			(\cdot, \cdot; H)\right\} \leq C_{7}(1 +  J_1(H))$ holds for $H \in \F$ and $\pi \in \Pi$.

		\end{enumerate}

	\end{assumption}
	\begin{theorem}
		\label{thm: ratio}
		Suppose Assumptions  \ref{ass: stationarity}-\ref{assumption: function classes}, and \ref{assumption: ratio function} hold.  Let $\widehat \omega^\pi_N$ be the estimated ratio function with tuning parameters $\mu_{2N} \simeq \lambda_{2N} \simeq (1+\VC(\Pi))(\log N)^{\frac{2+\alpha}{1+\alpha}}N^{-\frac{1}{1+\alpha}}$ defined in \eqref{ratio estimator}. 
		For any $m \geq 1$, there exists some positive constant $C_{8}$ such that with sufficiently large $N$, the following  holds with probability at least {$1-\frac{m}{N}$}
		\begin{align*}
		\sup_{\pi \in \Pi} \norm{\widehat \omega^\pi_N- \omega^\pi}^2 \leq C_{8} \left(1+\VC(\Pi)\right)\log(N)^{\frac{2+\alpha}{1+\alpha}}N^{-r_m},
		\end{align*}
		where  $r_{m} = \frac{1}{1+\alpha} - \frac{(1-\alpha)2^{-(m-1)}}{1+\alpha}$ and $C_8$ depends on $\beta_0, \beta_1, p_{\min}, \sup_{\pi \in \Pi} \norm{\omega^\pi}_\infty,  R_{max}, \kappa', F_{\max}, G_{\max}, \sup_{\pi \in \Pi} J_1(\tilde H^{\pi}), m$ and constants $C_1, C_2, C_4$ and $C_7$.
	\end{theorem}
	\begin{remark}
		Theorem \ref{thm: ratio} implies that our ratio estimator can achieve the near-optimal nonparametric convergence rate in the dependent data setting when $m$ is large enough, up to some logarithm factor. Again we have an additional term with respect to $\VC(\Pi)$ because our error bound is uniform over $\Pi$. While the derived rate may not be optimal compared with the classical non-parametric regression, as long as we can guarantee $r_m > \frac{1}{2}$ (e.g., $m\geq 3$), we are able to demonstrate the statistical efficiency of our estimator of $M(\beta, \pi)$ and establish the rate-optimal regret bound up to some logarithm factor. See the following two subsections. 
		%Our results generalize the results from \citep{liao2020batch}, which considers i.i.d setting and can only prove the error bound in terms of sample size $n$. We strengthen their results by showing the error bound in terms of the number of trajectories and decision points.
	\end{remark}

 %If assuming $\D$ is $\phi$-mixing, for example $\D$ satisfies uniform geometric ergodicity, we may also be able to obtain the similar bound. For sake of space, we do not discuss it further.
	\subsection{Statistical Efficiency}\label{sec: semiparametric}
	%\peng{Here we are talking about the efficiency of $M_N(\beta, \pi)$. But the previous section is about the regret. So it does not flow well. Put this section before regret?}
	In this section, we demonstrate the efficiency of our proposed estimators. In the i.i.d case, the variance of any asymptotic unbiased estimator is greater than or equal to the Cramer-Rao lower bounds. In the classic semi-parametric setting, the efficient bound is defined as the supremum of Cramer-Rao lower bounds over all parametric submodel. See \cite{van2000asymptotic} for details. Since our observations on each trajectory are dependent, we discuss the statistical efficiency of our proposed estimator $\widehat{M}_N(\beta, \pi)$ and $\max_{\beta \in \mathbb{R}} \widehat{M}_N(\beta, \pi)$ for $M(\beta, \pi)$ and $\max_{\beta \in \mathbb{R}} M(\beta, \pi)$ respectively under the notion of \citep{komunjer2010semiparametric} and \citep{kallus2019efficiently}. 
	
	Recall that the stochastic process $\{S_t, A_t\}_{t\geq1}$ is stationary by Assumption \ref{ass: stationarity}. Denote by $L(\{D_i\}_{i =1}^n; \varpi)$ as the likelihood function of a parametric sub-model indexed by a parameter $\varpi$ and
	\[
	L(\D_n; \varpi) = \Pi_{i=1}^n\{ d^{\pi_b}_\varpi(S^i_{1})\Pi_{t = 1}^{T_0} \pi^b_{\varpi} (A^i_{t}|S^i_{t}) P_\varpi (S^i_{t+1}|S^i_{t}, A^i_{t})\}. 
	\]
	The score function at the parameter $\varpi$ is then given by
	\begin{align*}
	\triangledown L_\varpi(\D_n) & = \frac{d \log L(\{D_i\}_{i=1}^n; \varpi)}{d \varpi}.
	\end{align*}
	We first discuss the efficiency of our estimator $\widehat{M}(\beta, \pi)$. Clearly, for a fixed $\beta$ and $\pi$, $M(\beta, \pi)
$ is a function of $\varpi$ and we denote its gradient with respect to $\varpi$ as
$$
\triangledown M(\varpi) = \frac{d \, M(\beta, \pi)}{d \, \varpi}.
$$
Denote the true parameter as $\varpi_0$. The semi-parametric efficiency bound for $M(\beta, \pi)$ can be defined as
	\begin{align}\label{EB bound}
		\EB(N) = N\sup \left\{\triangledown^T M(\varpi_0) \left\{\EE\left[ \triangledown L_{\varpi_0}(\D_n)\triangledown^T L_{\varpi_0}(\D_n)\right]\right\}^{-1} \triangledown M(\varpi_0)\right\},
	\end{align}
	where the supremum is taken over all parametric submodels that contain the true parameter. 
	Then we have the following theorem that shows the statistical efficiency of our estimator.
	\begin{theorem}\label{thm: efficiency bound}
		Under Assumptions \ref{ass: stationarity}-\ref{assumption: ratio function} and some regularity conditions, we have for any $\pi \in \Pi$ and $|\beta| \leq R_{\max}$,
		\begin{align}\label{EB}
			\frac{\sqrt{N}\left(\widehat M_N(\beta, \pi) - M(\beta, \pi)\right)}{\sqrt{\EB(N)}} \xrightarrow{d} \N(0, 1).
		\end{align}
		In particular,  we can show that the efficient bound
		$
		\EB(N) =  \EE\left[\psi^2(Z;U^{\pi,\beta}, \omega^\pi)\right],
		$,
		where $\psi(Z;U^{\pi, \beta}, \omega^\pi) ={\omega}^\pi(S, A) \left[\beta - \frac{1}{1-c}\left(\beta - \RR(S) \right)_+ +  U^{\pi, \beta}(S, A, S')-  M(\beta, \pi) \right]$ and $Z$ has the same distribution of $(S_t, A_t, S_{t+1})$.
	\end{theorem}
\begin{remark}
	The derivation of statistical efficient bound does not require the process $D$ to be stationary. However, in order to show our estimator is efficient, i.e., achieve this bound, we need to impose the stationarity assumption on the trajectory in order to show the in-sample bias decays faster than $N^{-\frac{1}{2}}$ in probability. This relies on the uniform finite sample error bounds for two nuisance functions and the doubly robust structure of our estimator. Finally, the martingale central limit theorem is applied to show its asymptotic normality.
\end{remark}
In the following, we show that $\max_{\beta \in \mathbb{R}} \widehat M_N(\beta, \pi)$ is also an efficient estimator for $\max_{\beta \in \mathbb{R}} M(\beta, \pi)$, which corresponds to the CVaR objective of the reward function under the stationary measure induced by $\pi$.
\begin{theorem}\label{thm: efficiency bound for CVaR}
	Suppose Assumptions \ref{ass: stationarity}-\ref{assumption: ratio function} hold and for every $\pi \in \Pi$, the optimal solution of $\max_{\beta \in \mathbb{R}} M(\beta, \pi)$ is unique, denoted by $\beta^\ast(\pi)$. Then under some regularity conditions, we have
	\begin{align}\label{EB for CVaR}
	\frac{\sqrt{N}\left(\max_{\beta \in \mathbb{R}}\widehat M_N(\beta, \pi) - \max_{\beta \in \mathbb{R}}M(\beta, \pi)\right)}{\sigma} \xrightarrow{d} \N(0, 1),
	\end{align}
	where $
	\sigma^2=  \EE\left[\psi^2(Z;U^{\pi,\beta^\ast(\pi)}, \omega^\pi)\right]. 
	$ Furthermore $\sigma^2$ is the statistical efficiency bound for estimating $\max_{\beta \in \mathbb{R}} M(\beta, \pi)$.
\end{theorem}
Note that $\max_{\beta \in \mathbb{R}}M(\beta, \pi) = \max_{|\beta| \leq R_{\max}} M(\beta, \pi)$. The proof of this theorem, which can be found in the Supplementary Material, relies on Danskin Theorem (e.g., \citep{danskin2012theory}) and the functional delta theorem (e.g., Theorem 5.7 of \citep{shapiro2021lectures}).  Finally we demonstrate the doubly robust property of $\widehat{M}(\beta, \pi)$, i.e., as long as one of the nuisance functions is estimated consistently, the proposed estimator is consistent. This is given by the following corollary.
	\begin{corollary}\label{lemma: doubly-rob} 
		Suppose the estimator $\UpiN$ and $\widehat{\omega}_N^\pi$ satisfy that $\norm{\UpiN-\bar {U}}$ and $\norm{\widehat{\omega}_N^\pi-\bar \omega}$ converge to $0$ in probability for some $\bar U$ and $\bar \omega$.  If either $\bar U= U^{\pi,\beta}$ or $\bar{\omega}= \omega^\pi$, then $\widehat M_N(\beta, \pi)$ converges to $M(\beta, \pi)$ and $\max_{\beta \in \mathbb{R}}\widehat M_N(\beta, \pi)$ converges to $\max_{\beta \in \mathbb{R}}M(\beta, \pi)$ in probability as $N\rightarrow \infty$. 
	\end{corollary}
The proof is similar to those in \citep{liao2020batch} and \citep{kallus2019efficiently}, so we omit here.	
	\subsection{Regret Guarantee}\label{sec: regret bound}
	Based on the uniform finite-sample error bounds for the two nuisance function estimations, we can derive the finite sample bound for the regret of $\widehat \pi_{N}^c$ defined in terms of $M(\beta, \pi)$:
\begin{align}\label{regret}
	\text{Regret}(\widehat \pi_{N}^c) &=\max_{\pi \in \Pi} \min_{u \in \Lambda^\pi_c} \, \, \EE_u\left[\RR(S)\right]- \min_{u \in \Lambda^{\widehat \pi_{N}^c}_c} \, \, \EE_u\left[\RR(S)\right].\notag \\
	 & =\max_{\pi \in \Pi, |\beta| \leq R_{\max}} M(\beta, \pi) - \max_{|\beta| \leq R_{\max}} M(\beta, \widehat \pi_{N}^c),
\end{align}
where the second equality is given by Theorem \ref{duality}. 
This regret bound can be interpreted as the difference between the smallest reward among the probability uncertainty set under the in-class optimal policy $\pi^\ast$ and that under the estimated policy $\widehat \pi_{N}^c$. 
	\begin{theorem}
		
		\label{thm: policy learning} 	
		Suppose the condition in Theorem \ref{duality} and Assumptions \ref{ass: stationarity} to \ref{assumption: ratio function} hold.  Let $\widehat \pi_{N}^c$ be the estimated policy obtained from \eqref{multi-level optimization} in which the nuisance functions are estimated with tuning parameters $\mu_{1N} \simeq \lambda_{1N} \simeq \mu_{2N} \simeq \lambda_{2N} \simeq (1+\VC(\Pi))^{\frac{1}{1+\alpha}}(\log N)^{\frac{2+\alpha}{1+\alpha}}N^{-\frac{1}{1+\alpha}}$. Then there exists a positive constant $C_{9}$ such that  for sufficiently large $N$, with probability at least $1- 1/N$, we have
		\begin{align*}
		\Regret(\widehat \pi_{N}^c) \leq C_{9} \log(N)\sqrt{\frac{\left(\VC(\Pi)+1\right)\underset{\pi \in \Pi,|\beta| \leq R_{\max}}{\sup} \EE\left[\psi^2(Z;U^{\pi,\beta}, \omega^\pi) \right]}{N}},
		\end{align*}
		where $C_{9}$ depends on $\beta_0, \beta_1, p_{\min}, \sup_{\pi \in \Pi} \norm{\omega^\pi}_\infty,  R_{max}, F_{max}, c$ and constants $C_1$ and  $C_2$.
		%where $\Sigma = \sup_{\pi \in \Pi} \EE\left[\phi^2(Z;U^\pi, \omega^\pi) \right]$.
	\end{theorem}
\begin{remark}
	Theorem \ref{thm: policy learning} gives, up to a logarithm factor, the rate-optimal regret bound of our learning method, compared with the rate-optimal regret bound developed in terms of the sample size $n$ in the infinite-horizon setting such as \citep{liao2020batch} and that in the finite-horizon setting such as \citep{athey2017efficient}. The logarithm factor is due to the dependence among observations.  One key reason why we are able to get the strong regret guarantee is because our estimator has the doubly robust property and achieves the statistical efficiency bound. To the best of our knowledge, this is the first regret bound in terms of total decision points in the batch RL. Such results imply that as long as the sample size or the horizon $T_0$ goes to infinity, the regret converges to 0, thus efficiently breaking the curse of horizon. When $c = 1$, we obtain the regret results for the estimated policy with respect to the long-term average reward MDP, which may be of independent interest. In addition, our theoretical results can be extended to discounted sum of rewards setting.
\end{remark}
\section{Numerical Study}
	
	%\subsection{Simulation Study} 
	In this section, we evaluate the performance of our proposed method via a simulation study. Our goal is to demonstrate the robustness of the learned policy under the proposed criterion in improving the generalizability, compared with two existing algorithms.  In order to promote a fast computation, we consider $\F_j$ and $\G_j$ for $1 \leq j \leq 2$ as some reproducing kernel Hilbert spaces (RKHSs). Then by the representer property, the optimization problem in \eqref{multi-level optimization} can be simplified as
	\begin{align}
	\max_{\pi \in \Pi, \beta \in \mathbb{R}} \quad \frac{\left(\widehat \mynu (\pi)\right)^T L \left(R^\beta_N - \tilde K (\pi)\widehat \alpha(\pi, \beta)\right)}{\widehat \mynu (\pi)^T L 1_N}\label{final optimization},
	\end{align}
	where $L$ is a kernel matrix generated by $\{S_t^{i}, A_t^{i}\}_{1\leq i \leq n, 1 \leq t \leq T_0}$, $\tilde{K}(\pi)$ is some shifted kernel matrix generated by $D_n$ and $\pi$, $\widehat \alpha(\pi, \beta)$ and $\widehat \mynu(\pi)$ are the corresponding coefficients for the estimators of $U^{\pi, \beta}$ and $\omega^\pi$ respectively, and $R_N^\beta = (\beta - \frac{1}{1-c}(\beta - \mathcal{R}(S_t^i)))_{1 \leq i \leq n, 1 \leq t \leq T_0}$, a vector with length $N$. A block coordinate ascent algorithm is then proposed to solve the optimization problem \eqref{final optimization}. The details can be found in Section 2 of the Supplementary Material.

	We consider the following simulation setting, which is similar as that in \cite{luckett2019estimating} (while their goal is  to learn an in-class optimal policy that maximizes the cumulative sum of discounted rewards). Specifically, we initialize two dimensional state vector $S_1 = (S_{1, 1}, S_{1, 2})$ by a standard multivariate Gaussian distribution. Consider a two-arm setting, where $\A = \{0, 1\}$. Given the current action $A_t \in \A$ and state $S_t$, the next state is generated by:
	\begin{align*}
	& S_{t+1, 1} = \frac{3}{4} (2A_t - 1) S_{t, 1} + \frac{1}{4}S_{t, 1}S_{t, 2} + \varepsilon_{t, 1}, \\
	&S_{t+1, 2} = \frac{3}{4} ( 1 - 2A_t) S_{t, 2} - \frac{1}{4}S_{t, 1}S_{t, 2} + \varepsilon_{t, 2},
	\end{align*}
	where each $\varepsilon_{t, j}$ follows independently $N(0, 1/2)$ for $j = 1, 2$. The reward function $R_{t}$ is given as
	$R_{t} = 2S_{t, 1} + S_{t, 2},$
	for $t = 1, \cdots, T_0$.
	We consider the behavior policy to be uniformly random, i.e., choosing each action with the same probability.
	
	Based on this generative model, we generate multiple trajectories with $n= 25$ and $T_0 = 24$ as our training data. Then we apply our method with $c$ ranging from $0.1, 0.3$ and $0.5$ to learn three different robust policies. Specifically, for each $c$, we consider RKHS with Gaussian kernels due to its universal consistency. The bandwidth is selected based on the median heuristic, e.g., median of pairwise distance \citep{fukumizu2009kernel}. Other tuning parameters are selected based on a min-max cross-validation procedure, which can be found in Algorithm 2 of  the Supplementary Material. To model policies, we consider the following stochastic parametrized policy class $\Pi$ indexed by $\theta$:
	\[
	\Pi = \left\{\pi \; \vert \; \pi(1 \,|\, s, \theta) = \frac{\exp(s^T\theta)}{1 + \exp(s^T\theta)}, \; \; \|\theta\|_{\infty} \leq c_0, \; \; \theta \in \mathbb{R}^p
	\right\},
	\]
	for some pre-specified constant $c_0>0$. Then our goal becomes optimizing $\theta$ to get the robust optimal in-class policy. Here the box constraint in $\Pi$ is used to regularize the policy class from overfitting. In our numerical study, we choose $c_0 = 10$. Furthermore, the constraint on $\theta$ can be naturally incorporated into the proposed block coordinate ascent algorithm. See Section 2 of Supplementary Material. Due to the non-convexity, it is not trivial to analyze the computational complexity of our algorithm. As this is beyond the scope of this paper, we leave it  as future work.

	 For comparison, we also implement two policy learning methods. One uses the long-term average reward for policy optimization proposed by \citep{liao2020batch} and the other is called  V-learning by \citep{luckett2019estimating} for optimizing the discounted sum of rewards. To test the performance of different methods and demonstrate the robustness of our method for varying length of horizon, we compute the  average rewards of each learned policy over different time-horizons using independent test dataset based on the above generative models. More specifically, we generate a test dataset with $1000$ trajectories, where actions are generated by each learned policy and compute the average rewards over $T= 50, \cdots, 100$. The results are shown in Figure 1 (left). As we can see, the performances of all methods are similar while our method performs slightly better. The slightly better performance comes from the uncertain set we consider to robustify the performance of our learned policy in improving the average rewards over different lengths of time horizons.

	To further test the performance of our method subject to the distributional change in the initial distribution, we change the initial state distribution from the standard bivariate normal distribution to a $t$-distribution with the degree of freedom $2$. We make this choice because $t$-distribution is a heavy-tailed distribution, significantly different from the normal distribution. Therefore we may be able to see how the performance of each method differs under this change. Again for each policy, we use the same method before to calculate its corresponding average rewards over different horizons. The corresponding results are provided in Figure 1 (right). It can be observed that our method outperforms the other two methods in terms of the average rewards over time horizons ranging from $50$ to $100$, demonstrating the robust performance of our method. The superior performance of our algorithm in this case mainly comes from the uncertainty set $\Lambda_c^\pi$, which considers all the initial distribution. By improving the worst case performance, the resulting policy can protect against the potential distributional change in the initial distribution. Lastly, we remark that the average rewards reported here is much larger than the previous scenario because of the heavy-tailed initial state distribution.
	\begin{figure}[H]
		\centering
		\includegraphics[scale=0.7]{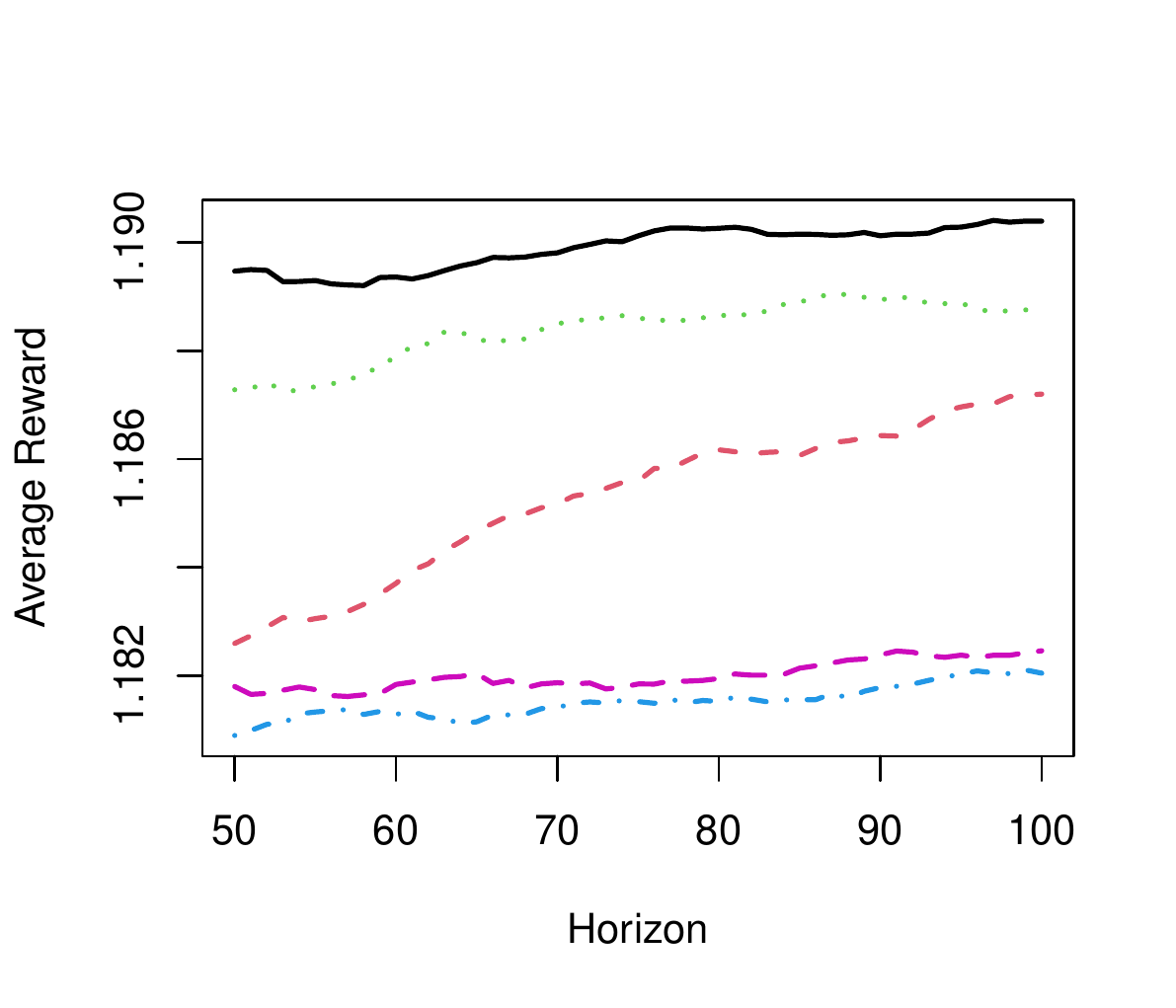}%
		\includegraphics[scale=0.7]{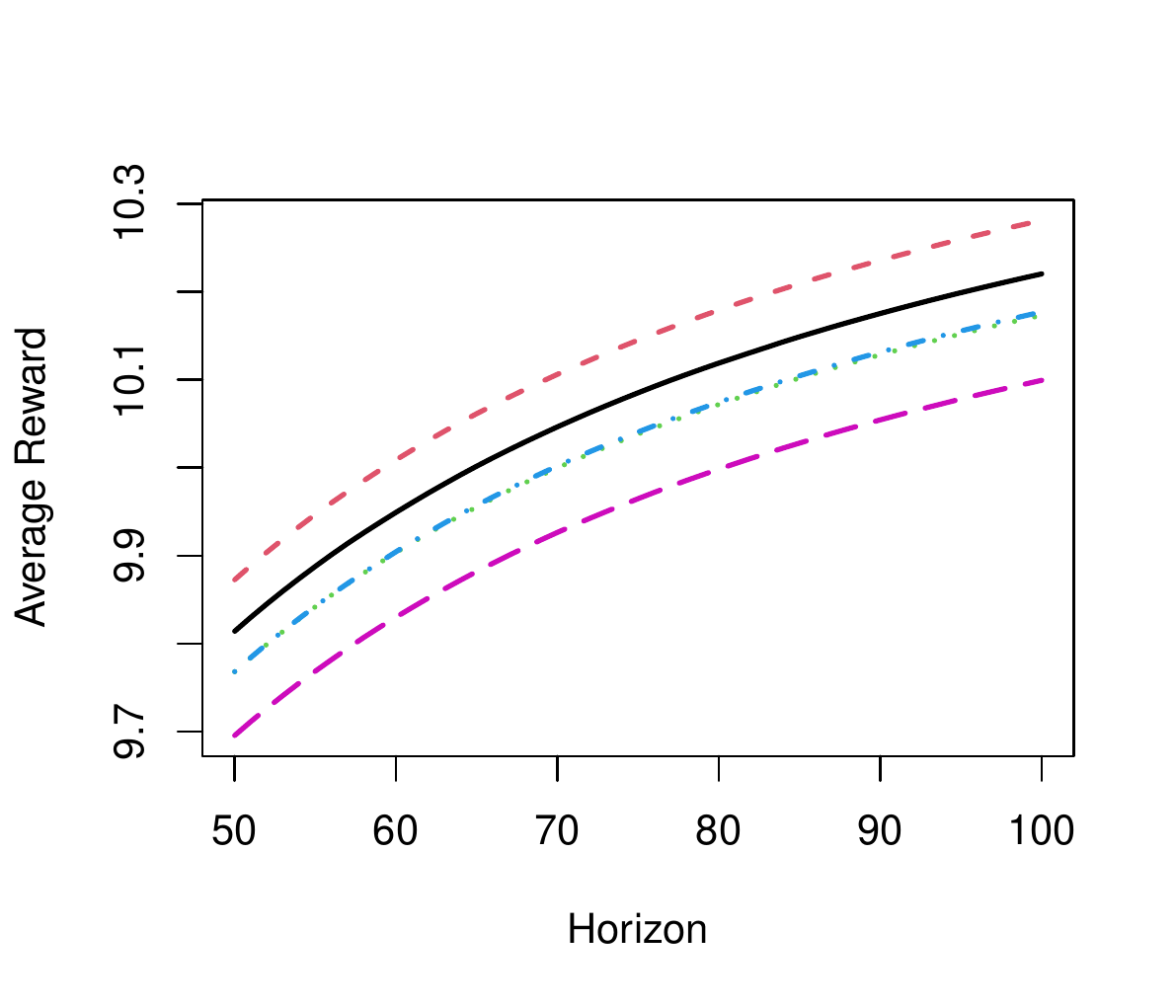}
		\label{fig:norm}
		\caption{The average rewards of five policies over horizons ranging from $50$ to $100$ when the initial state distribution is bivariate standard normal distribution (\textbf{left}) and the initial state distribution is $t$-distribution with degree of freedom $2$ (\textbf{right}). The black solid curve corresponds to the proposed robust policy using $c=0.5$, red short-dashed curve using $c=0.3$ and green dotted curve using $c=0.1$. The blue dashed curve with dots corresponds to the policy using the long-term average reward and purple long-dashed curve corresponds to the policy using the average cumulative  discounted rewards with discount rate $\gamma = 0.9$.
		}
	\end{figure}

	\section{Discussion}

	In this work, we propose a robust criterion to evaluate a policy by average rewards with respect to a set of distributions centered at the policy induced stationary distribution. It can be shown to contain average rewards across varying planning horizons with different reference distributions. Based on this criterion, we developed a data-efficient learning method to estimate the corresponding optimal policy that can maximize the worst case performance of some uncertainty set, improving the generalizability of the learned policy. A rate-optimal regret bound, up to a logarithm factor, was established in terms of the number of trajectories and decision points in each trajectory. A numerical study demonstrates the decent performance of our proposed method.
	
	In the following, we discuss the setting where the reward $R_t$ also depends on the current action. Define the expected reawrd by $r(s, a) = \EE[R_t|S_t=s, A_t=a]$. If we consider $\Pi$ as a class of deterministic policies, then we can correspondingly define $\U^\pi_c$ as
	$
	 \left\{\EE^\pi_u\left[r(S, A)\right] \, | \, u \in \Lambda^\pi_c   \right\}.
$
	To obtain $\pi^\ast_c$ using the modified $\U^\pi_c$, under the assumption that the essential minimums of $r(s, a)$ under $d^\pi$ are the same for every $\pi \in \Pi$,  one can show that it is  equivalent to solving
	$
	\max_{\pi \in \Pi, \beta \in \mathbb{R}} \left\{\beta -\frac{1}{\left(1-c\right)}\EE^\pi_{d^\pi}\left[\left(-r(S, A)+\beta \right)_+\right]\right\}.
	$
    If we consider a stochastic policy class, then we need to solve
	$
	\max_{\pi \in \Pi, \beta \in \mathbb{R}} \left\{\beta -\frac{1}{\left(1-c\right)}\EE_{d^\pi}\left[\left(-\sum_{a\in \A}\pi(a \given S)r(S, a)+\beta \right)_+\right]\right\}.
	$
	To obtain estimators for the above two objective functions, we need to implement an additional step by estimating the conditional reward function $r(s,a)$. This can be done by using some standard supervised learning techniques. In some applications, it may be more natural to define the reward as a function of the next state. In this case, one can include the reward into the state and still use our reward formulation. 
	
	Lastly, we discuss some future research directions. From the theoretical perspective, it will be interesting to derive the finite sample regret bound for the batch policy learning in the infinite-horizon MDP without stationarity and positivity assumptions. From the optimization perspective, our current algorithm requires a moderate computation and large memory due to the nonparametric estimation and the policy-dependent structure of nuisance functions. It is thus desirable to develop a more computationally efficient algorithm. One possible remedy is to consider  zero-order optimization method. In the proposed algorithm, we consider tuning parameters independent of the policy. It will be interesting to investigate a more general setting and study how to perform model selection in  batch RL, which seems far less studied in the literature. Another possible line of the research is to extend our proposed efficient policy learning method from the batch setting to the online setting. One challenging question is how to design an online algorithm to balance the evaluation of a policy and the search for a new policy given that all nuisance functions are policy dependent. Studying two-timescale stochastic algorithms such as \cite{konda2004convergence} may be a good starting point.

\bibliographystyle{abbrvnat}
\bibliography{reference}

\newpage
\begin{appendices}
\section{Introduction}
%In Section 2 of this supplementary material, we elaborate more about our nuisance function estimation.  We provide details about our computation for $\widehat{\pi}^\ast_c$, and how to select tuning parameters and specify the constant $c$ in determining the size of $\Lambda_c^\pi$. In Section 3, we provide an analytic example to justify the use of our proposed robust criterion. In Section 4, we give all our technical proofs.
In Section 2 of this Supplementary Material, we provide details about computing $\widehat{\pi}^c_N$, and show how to select all tuning parameters in our learning method and also the constant $c$ in determining the size of $\Lambda_c^\pi$. In Section 3, we give all our technical proofs of theoretical results in the main text.
\section{Optimization and related computation}
We start with our overall optimization problem.

\noindent \underbar{Upper level optimization task}:	\begin{align}
& \max_{\pi \in \Pi, \beta \in \mathbb{R}} && \frac{\Pn \left\{(1/T_0)\sum_{t=1}^{T_0} \widehat \omega^{\pi}_N(S_{t}, A_{t}) \left[\beta - \frac{1}{1-c}\left(\beta - R_{t}\right)_+ + \widehat U^{\pi, \beta}_N(S_t, A_t, S_{t+1})\right]\right\} }{\Pn \left\{(1/T_0)\sum_{t=1}^{T_0} \widehat \omega^{\pi}_N(S_{t}, A_{t})\right\} }\label{multi-level optimization app} 
\end{align}

\noindent \underbar{Lower level optimization task 1}:\begin{align}
& && (\widehat \eta_N^{\pi, \beta}, \widehat Q_N^{\pi, \beta}) = \argminb_{(\eta, Q) \in \R \times \F_1} \Pn \left[\frac{1}{T_0} \sum_{t=1}^{T_0} \left[\widehat g_N^{\pi, \beta}(S_t, A_t; \eta, Q)\right]^2\right] + \lambda_{1N} J_1^2(Q)\label{lower level1 app}\\
& &&\text{such that} \; \; \widehat g_N^{\pi, \beta}(\cdot, \cdot; \eta, Q) = \argmin_{h \in \G_1} \Pn\Big[ \frac{1}{T_0}\sum_{t=1}^{T_0} \big( \delta^{\pi, \beta}(Z_t; \eta, Q) - g(S_t, A_t)\big)^2  \Big] + \mu_{1N} J_2^2(g) \label{lower level2 app}
\end{align}

\noindent \underbar{Lower level optimization task 2}:\begin{align}
& &&\widehat H_N^\pi(\cdot,\cdot) = \argminb_{H \in \F_2} \Pn\Big[\frac{1}{T_0}\sum_{t=1}^{T_0}  \widehat h^2_N(S_t, A_t; H)\Big] + \lambda_{2N} J_1^2(H)\label{lower level3 app}\\[0.1in]
& && \text{such that} \, \, \widehat h_N(\cdot, \cdot; H) = \argminb_{h \in \G_2} \Pn\Big[\frac{1}{T_0} \sum_{t=1}^{T_0} \big(\Delta^\pi(Z_t; H)- h(S_t, A_t)\big)^2\Big]+  \mu_{2N} J_2^2(h) \label{lower level4 app},
\end{align}
where we recall that $\delta^{\pi, \beta}(Z_t; \eta, Q) =\beta - \frac{1}{1-c}\left(\beta - R_{t}\right)_+ + \sum_{a'} \pi(a'|S_{t+1}) Q(S_{t+1}, a) - Q(S_t, A_t) - \eta$, $\widehat{U}^{\pi, \beta}_N(s, a, s') = \sum_{a \in \cal A} \pi(a|s') \widehat{Q}^{\pi, \beta}_N(s', a) - \widehat{Q}^{\pi, \beta}_N(s, a)$, $\Delta^{\pi}(Z_t; H) = 1 - H(S_t, A_t) + \sum_{a'} \pi(a'|S_{t+1})  H(S_{t+1}, a')$, and $\widehat \omega^\pi_N(s, a) = \widehat h_N(s, a; \widehat H_N^\pi) /\Pn[(1/T)\sum_{t=1}^{T} \widehat h_N(s, a; \widehat H_N^\pi)]$.
\subsection{Optimization Algorithm}
As discussed at the end of Section 4 of the main text, the overall optimization problem is bi-level, where the upper level serves for searching an optimal robust policy and the lower level represents feasible sets, i.e., the estimation of our nuisance functions. In order to compute \eqref{multi-level optimization app}, we first need to specify spaces i.e., $\F_1, \F_2, \G_1$ and $\G_2$. For simplicity, we assume $\F_1 = \F_2$ and $\G_1 = \G_2$ and consider all these spaces as reproducing kernel Hilbert spaces (RKHSs) with radial basis function. This kernel has the universal property that can approximate any continuous functions under some mild conditions. In addition, considering RKHSs promotes efficient computations due to the representer theorem. Note that two parallel lower level problems \eqref{lower level1 app}-\eqref{lower level2 app} and \eqref{lower level3 app}-\eqref{lower level4 app} can be regarded as two nested kernel ridge regressions. By using the representer theorem, we can compute closed-form solutions for all our nuisance functions. Next, we specify the policy class, where we consider a class of stochastic paramterized policies indexed by $\theta$. For example, if we consider the binary-action space, i.e., $\A = \{0, 1\}$, then we can model $\Pi$ as
\[
\Pi = \left\{\pi \; \middle\vert \; \pi(1 \,|\, s, \theta) = \frac{\exp(s^T\theta)}{1 + \exp(s^T\theta)}, \; \; \|\theta\|_{\infty} \leq \bar{c}, \; \; \theta \in \mathbb{R}^d, \; \; s \in \S
\right\},
\]
where $\| \bullet\|_\infty$ is the infinity norm, $\bar{c}$ is some positive constant for keeping stochasticity of the learned policy. We remark that multiple action cases and other models for the policy class can be defined similarly.

For the remaining of this section, we describe our optimization algorithm  to obtained our estimated policy $\widehat{\pi}_N^c$ and the estimated auxiliary parameter $\widehat{\beta}$. We propose to use the block update algorithm. For each iteration, we first fixed $\pi$ (or equivalently $\theta$), and maximize $\widehat{M}_N(\beta, \pi)$ over $\beta$. Note that this is an one-dimensional optimization problem, which thus can be solved efficiently with the guarantee of finding a minimum. We remark that $\widehat{M}_N(\beta, \pi)$ is a piecewise linear function with respect to $\beta$ and thus an optimal solution must be one element in the vector $\{R_t^i\}_{1\leq i \leq n, 1 \leq t \leq T_0}$. Next, we fixed $\beta$ and maximize $\widehat{M}_N(\beta, \pi)$ over $\pi$. We  use a limited-memory Broyden-Fletcher-Goldfarb-Shanno  algorithm with box constraints (L-BFGS-B) to compute the solution $\widehat{\theta}$ \citep{liu1989limited}. To avoid bad solutions, in this step, we randomly select multiple initial points and search for the best solution. The full procedure can be found in Algorithm \ref{alg:optimization}.
\begin{algorithm}
	\caption{Maximize $\widehat{M}_N(\beta, \pi)$}
	\label{alg:optimization}
	\textbf{Input:} Data $\D_n$, initial $\theta_0$ and $\beta_0$, a constant $\bar{c}$, and tolerance $\epsilon_{tol} > 0$.
	
	Repeat for $t = 0, \cdots,$ do till $\left\|\theta_{t+1} - \theta_{t} \right\| \leq \min\left\{\left\|\theta_{t} \right\|_2, 1 \right\} \epsilon_{tol}$
	
	\Indp Compute $\beta_{t+1}$ by maximizing $\widehat{M}_N(\beta, \pi_{\theta_t})$ over $\beta$ with an initial $\beta_t$;
	
	Compute $\theta_{t+1}$ by maximizing $\widehat{M}_N(\beta_{t+1}, \pi)$ over $\Pi$ with an initial $\theta_t$. 
	
	\Indm	\textbf{Output:} $\widehat{\theta}$ and $\widehat{\beta}$.
\end{algorithm}

\textbf{Discussion of Step 4 in Algorithm \ref{alg:optimization}:} It is noted that Step 4 in Algorithm \ref{alg:optimization} only involves optimization over $\theta$ while keeping $\beta$ fixed. We rewrite the training data $\D_n$ into tuples $Z_h = \{S_h, A_h, R_h, S_{h+1}\}$ for $h = 1, \dots, N$, where $h$ indexes the tuple of transition sample in the training set $\D_n$,   $S_h$ and $S_{h+1}$ are the current and next states and $R_h$ is the associated reward. Let $W_h = (S_h, A_h)$  be one state-action pair, and $W_h' = (S_h, A_h, S_{h+1})$ be one state-action-next-state pair. Denote the kernel function for the state as $k_0(s_1, s_2)$, where $s_1, s_2 \in \S$. Then the state-action kernel function can be define as $k((s_1, a_1), (s_2, a_2)) = \indicator{a_1 = a_2} k_0(s_1, s_2)$. Recall that we have to restrict the function space $\F_1$ such that $Q(s^*, a^*) = 0$ for all $Q \in \F_1$ and $\F_2$ such that $H(s^*, a^*) = 0$ for all $H \in \F_2$ respectively so as to avoid the identification issue. For ease of presentation, in the following, we omit the subscript for $\F_j$ and $\G_j$ when there is no confusion.   Thus for any given kernel function $k$ defined on $\S \times \A$, we make the following transformation by defining $k(W_1, W_2) = k(W_1, W_2) - k(( s^*, a^*), W_2) - k(W_1, (s^*, a^*)) + k((s^*, a^*), (s^*, a^*))$ with some abuse of notations. One can check that the induced RKHS by this $k(\cdot, \cdot)$ satisfies the constraint in $\F$ automatically.    

We denote kernel functions for $\F$ and $\G$ by $k(\cdot, \cdot), l(\cdot, \cdot)$ respectively. The corresponding inner products are defined as $\innerprod{\cdot}{\cdot}_\F$ and $\innerprod{\cdot}{\cdot}_\G$. In terms of the inner minimization problem \eqref{lower level1 app}-\eqref{lower level2 app}, the closed form solution can be obtained by representer theorem. For example, $\widehat g^{\pi, \beta}_N(\cdot, \cdot; \eta, Q) =  \sum_{h=1}^{N} l(W_h, \cdot) \widehat{\gamma}(\eta, Q)$, where $\widehat{\gamma}(\eta, Q) = (L + \mu_1 I_N)^{-1} \delta_N^{{\pi}}(\eta, Q)$, $L$ is the kernel matrix of $l$, $\mu_1 = \mu_{1N} N$, and $\delta^\pi_N(\eta, Q) = (\delta^\pi(W'_h; \eta, Q))_{h=1}^N$ is a vector of TD error.  Moreover, each temporal difference error can be further written as $\delta^{\pi, \beta}(W'_h; \eta, Q) = \beta - \frac{1}{1-c}\left(\beta - R\right)_+ - \eta- \innerprod{Q}{\tilde k_{W'}}_\G$, where 
\[
\tilde k_{W'}(\cdot) = k(W, \cdot) - \sum_{a'} \pi(a'|S') k((S', a'), \cdot) \in \F
\]
One can demonstrate that $\widehat Q^{\pi, \beta}_N$ in \eqref{lower level1 app} can be expressed by the linear span: $\dkh{\sum_{h=1}^N \alpha_h \tilde k_{W_h'}(\cdot):  \alpha_h \in \R, h = 1, \dots, N}$ according to the representer property.
%Note that 
%\begin{align*}
%L(\beta, Q) = (1/N) \sum_{h=1}^{N} \left[\theta(\eta, Q)^\transpose l(W_h, U)\right]^2 \\
%= (1/N) \theta(\eta, Q)^\transpose \left[\sum_{h=1}^{N}  l(U, W_h) l(W_h, U)\right] \theta(\eta, Q) \\
%= (1/N) \theta(\eta, Q)^\transpose \left[\sum_{h=1}^{N}  L\left[, h\right] L\left[h, \right]\right] \theta(\eta, Q) \\
%= (1/N) \delta_N(\eta, Q)^\transpose (L + \mu)^{-1}  L^2 (L + \mu)^{-1} \delta_N(\eta, Q) \\
%= (1/N) \delta_N(\eta, Q)^\transpose M \delta_N(\eta, Q)
%\end{align*}
%where $M = (L + \mu)^{-1}  L^2 (L + \mu)^{-1}$. And
%\\left[
%\delta_N(\eta, Q) = R_N - F \beta  - \tilde K \alpha
%\\right]
Then the optimization problem \eqref{lower level1 app}-\eqref{lower level2 app} is equivalent to solving
\begin{align}
\label{value estimation}
(\widehat \eta_N^{\pi, \beta}, \widehat \alpha(\pi, \beta)) = \argminb_{\eta \in \R, \alpha \in \R^N} (R^\beta_N - \eta 1_N  - \tilde K(\pi) \alpha)^\transpose M (R^\beta_N - \eta 1_N  - \tilde K(\pi) \alpha) + \lambda_1 \alpha^\transpose \tilde K(\pi) \alpha
\end{align} 
where $R^\beta_N = (\beta - \frac{1}{1-c}\left(\beta - R_{h}\right)_+)_{h = 1}^N$, $\tilde K(\pi) = ( \innerprod{\tilde k_{W_h'}}{\tilde k_{W_{j}'}}_\F)_{j, h = 1}^N$, $M = (L + \mu_1 I_N)^{-1}  L^2 (L + \mu_1 I_N)^{-1}$, $1_N$ is a length-$N$ vector of all ones, $\lambda_1 = \lambda_{1N} N$ and $\alpha = (\alpha_h)_{h = 1}^N$ is a vector of length $N$. Note that the $(h, k)$-th element of the matrix $\tilde K(\pi)\left[h, k\right]$ can be further calculated as
\begin{align*}
\innerprod{\tilde k_{W_h'}}{\tilde k_{W_{j}'}}_\F  & = k(W_h, W_{j}) - \sum_{a'} \pi(a'|S_h') k((S_h', a'), W_j) - \sum_{a'} \pi(a'|S_j') k((S_j', a'), W_h)  \\
& \qquad + \sum_{a_h'} \sum_{a_{j}'} \pi(a_h'|S_h') \pi(a_j'|S_j') k((S_h', a_h'), (S_j', a_j')).
\end{align*}
We make $\tilde K(\pi)$ and $ \widehat \alpha(\pi, \beta)$ as functions of $\pi$ and $\beta$ to explicitly indicate their  dependency on the policy $\pi$ and the auxiliary parameter $\beta$. The first-order optimality implies that $(\widehat \eta_N^{\pi, \beta}, \widehat \alpha(\pi, \beta))$ satisfies
\begin{align*}
1_N^\transpose M 1_N \widehat \eta_N^{\pi, \beta}  = 1_N^\transpose M(R^\beta_N - \tilde K(\pi) \widehat \alpha(\pi, \beta)) \\
%(\tilde K M \tilde K  + \lambda \tilde K) \alpha = \tilde K M (R_n - F \beta) \iff 
( M \tilde K(\pi)  + \lambda_1 I_N) \widehat \alpha(\pi, \beta) =  M (R^\beta_N - M 1_N \widehat \eta_N^{\pi, \beta}), 
\end{align*}
which gives
\begin{align}
&(M \tilde{K}(\pi)+\lambda_1 I_N - MF(1_N^TM1_N)^{-1}1_N^TM\tilde{K}(\pi))\widehat \alpha (\pi, \beta)\label{for gradient 0}\\
=& (I_N - 1_N(1_N^TM1_N)^{-1}1_N^T)MR^\beta_N
\label{for gradient 1}
\end{align}
and thus the corresponding $\{\widehat{U}_N^\pi(W'_h)\}_{h = 1}^N =-\tilde{K}(\pi)  \widehat \alpha(\pi, \beta)$. In order to apply L-BFGS-B, we need to compute the Jacobian matrix of the vector $\{\widehat{U}_N^\pi(W'_h)\}_{h = 1}^N$ with respect to $\theta$. Based on the above equations, we know
\begin{align*}
\frac{\partial \{\widehat{U}_N^\pi(W'_h)\}_{h = 1}^N}{\partial \theta} = -\frac{\partial \tilde{K}(\pi) }{\partial \theta} \otimes \widehat{\alpha}(\pi, \beta) - \tilde{K}(\pi) \frac{\partial \widehat{\alpha}(\pi, \beta) }{\partial \theta},
\end{align*}
where $\otimes$ is denoted as a tensor product. Here $\frac{\partial \tilde{K}(\pi) }{\partial \theta}$ is a $\R^N \otimes \R^N \otimes \R^p$ tensor, where the $(i, j, k)$-th  element is the partial derivative $\frac{\partial \left[ \tilde{K}(\pi)\right]_{i, j}}{\partial \theta_k}$. In addition, $\frac{\partial \widehat \alpha(\pi, \beta) }{\partial \theta}$ can be calculated via implicit theorem based on the equation \eqref{for gradient 0}-\eqref{for gradient 1}, i.e.,
\begin{align}
&\left(M \otimes \frac{\partial \tilde{K}(\pi) }{\partial \theta}+\lambda I_N - MF(1_N^TM1_N)^{-1}1_N^TM\otimes \frac{\partial \tilde{K}(\pi) }{\partial \theta}\right)\widehat \alpha (\pi, \beta)\\
=& -(M \tilde{K}(\pi)+\lambda I_N - MF(1_N^TM1_N)^{-1}1_N^TM\tilde{K}(\pi))\frac{\partial \widehat \alpha (\pi, \beta)}{\partial \theta},
\label{for gradient 2}
\end{align}
which gives the expression of $\frac{\partial \widehat \alpha (\pi)}{\partial \theta}$, a $N$ by $p$ matrix.

We can use the same approach to get the closed-form solution for the problem \eqref{lower level3 app}-\eqref{lower level4 app} and compute its corresponding gradient with respect to $\theta$.  By some calculation, we can get $\{ \hn(W_j, \Hn )\}_{j=1}^N = L\widehat \mynu(\pi)$, where $\hn(W_j, \Hn ) = \sum_{j=1}^N \widehat\mynu_j(\pi) l(W_j, \cdot) $ and $\nu = (\widehat\mynu_j(\pi))_{j=1}^N $ satisfying the following two equations:
\begin{align}
&(M \tilde K(\pi)  + \lambda_2 I_N) \widehat \varphi(\pi) &&=  M 1_N \label{for gradient 3}\\
&(L + \mu_2 I_N) \widehat \mynu (\pi)&& =  1_N - \tilde K(\pi) \widehat \varphi(\pi) \label{for gradient 4},
\end{align}
again by the representer theorem,
where $\widehat{\varphi}(\pi)$ is estimated coefficient associated with $\tilde{K}(\pi)$, $\lambda_2 = \lambda_{2N}N$ and $\mu_2 = \mu_{2N}N$.
The Jacobian matrix of $\{ \hn(W_j, \Hn )\}_{j=1}^N$ can be computed by again using the implicit theorem on equations \eqref{for gradient 3} and \eqref{for gradient 4}. More specifically, we need to solve $\frac{\partial \widehat \mynu(\pi)}{\partial \theta}$ based on the following two equations.
\begin{align}
&\left(M \otimes \frac{\partial \tilde{K}(\pi) }{\partial \theta}\right)  \widehat \varphi(\pi) + \left(M \tilde K (\pi) + \lambda_2 I_N\right) \frac{\partial \widehat \varphi(\pi)}{\partial \theta}&&= 0. \label{for gradient 5}\\
&(L + \mu_2 I_N) \frac{\partial \widehat \mynu(\pi)}{\partial \theta} +  \tilde K (\pi) \frac{\partial \widehat \varphi(\pi)}{\partial \theta} + \frac{\partial \tilde{K}(\pi) }{\partial \theta} \otimes \frac{\partial \widehat \varphi(\pi)}{\partial \theta} && =  0. \label{for gradient 6}
\end{align}
Then we have
\begin{align*}
\frac{\partial \{ \gn(W_h, \Hn )\}_{j=1}^N}{\partial \theta} = L \frac{\partial \widehat \mynu(\pi)}{\partial \theta}.
\end{align*}

Summarizing together by plugging all the intermediate results into the objective function of our upper optimization problem \eqref{multi-level optimization app}, we can simplify step 4 in Algorithm \ref{alg:optimization} as
\begin{align}
& \max_{\pi \in \Pi_{\Theta}} && \frac{\left(\widehat \mynu (\pi)\right)^T L \left(R^\beta_N - \tilde K (\pi)\widehat \alpha(\pi, \beta)\right)}{\widehat \mynu (\pi)^T L 1_N}\label{final optimization app}.
\end{align}
The corresponding gradient with respect to $\theta$ can be computed directly as
\begin{align*}
&\frac{\left(\frac{\partial \widehat \mynu(\pi)}{\partial \theta}\right)^T L \left(R^\beta_N - \tilde K (\pi)\widehat \alpha(\pi, \beta)\right) - \left(\widehat \mynu (\pi)\right)^T L \left(\frac{\partial \tilde{K}(\pi) }{\partial \theta} \otimes \widehat \alpha(\pi, \beta) + \tilde{K}(\pi) \frac{\partial \widehat \alpha(\pi, \beta) }{\partial \theta}\right) \left(\widehat \mynu (\pi)^T L 1_N\right)}{\left(\widehat \mynu (\pi)^T L 1_N\right)^2} \\[0.1in]
-&\frac{\left(\frac{\partial \widehat \mynu(\pi)}{\partial \theta}\right)^T L \left(R^\beta_N - \tilde K (\pi)\widehat \alpha(\pi, \beta)\right)\left(\widehat \mynu (\pi)^T L 1_N\right)}{\left(\widehat \mynu (\pi)^T L 1_N\right)^2}.
\end{align*}
%\textbf{Discussion of Step 3 in Algorithm \ref{alg:optimization}.} Note that $\widehat{M}_N(\beta, \pi)$ is a piecewise linear function with respect to $\beta$ as we can see in \eqref{final optimization app}, $\widehat{\alpha}(\pi, \beta)$ is a linear function of $R_N^\beta$ and $R_N^\beta$ is a piecewise linear mapping of $\beta$. Since $|\beta| \leq R_{\max}$, then $\widehat \beta$ must be one of $\{R_h\}_{h = 1}^{N}$. Many one-dimensional optimization algorithms can be used such as golden-section search.

\subsection{Selection of Tuning Parameters}
In this subsection, we discuss the choice of tuning parameters in our method. 
The bandwidths in the Gaussian kernels are selected using median heuristic, e.g., median of pairwise distance \citep{fukumizu2009kernel}. The tuning parameters $(\lambda_{1N}, \mu_{1N})$ and $(\lambda_{2N}, \mu_{2N})$ are selected based on 3-fold cross-validation. We assume that all these tuning parameters are independent of the policy $\pi$ and $\beta$ so that we can select them based the estimation of ratio and relative value functions using some randomly generated policies and $\beta$. We adopt ideas from \citep{farahmand2011model} and \citep{liao2020batch}. Specifically, for the tuning parameters $(\lambda_{1N}, \mu_{1N})$ in the estimation of relative value function, we focus on \eqref{lower level1 app}-\eqref{lower level2 app} using cross-validation. For the tuning parameters $(\lambda_{2N}, \mu_{2N})$ in the estimation of ratio function, we focus on \eqref{lower level3 app}-\eqref{lower level4 app}. Both of the cross-validation procedures are based on choosing the tuning parameters that have the smallest estimated projected bellman errors on the validation set among a pre-specified tuning set. The details of selecting these tuning parameters can be found in Algorithm \ref{alg:cross-validation}.

\begin{algorithm}[!ht]
	\caption{Tuning parameters selection via cross-validation}
	\label{alg:cross-validation}
	\textbf{Input:} Data $ \{Z_h \}_{h=1}^{N} $, a set of $M$ policies $\left\{\pi_1, \cdots, \pi_M \right\} \subset \Pi$, a set of $\left\{\beta_1, \cdots, \beta_L \right\}$,  a set of $J$ candidate tuning parameters $\{(\mu^{(j)}_{1N}, \lambda^{(j)}_{1N})\}_{j=1}^J$ in the relative value function estimation, and a set of $J$ candidate tuning parameters $\{(\mu^{(j)}_{2N}, \lambda^{(j)}_{2N})\}_{j=1}^J$ in the ratio function estimation.
	
	Randomly split Data into $K$ subsets: $\{Z_h \}_{h=1}^{N} = \left\{D_k\right\}_{k = 1}^K$
	
	Denote $e^{(1)}(m, l, j)$ and $e^{(2)}(m, l,j)$ as the total validation error for $m$-th policy, $l$-th $\beta$ and $j$-th pair of tuning parameters in value and ratio function estimation respectively, for $m = 1, \cdots M$, $l= 1, \cdots, L$ and $j = 1, \cdots, J$. Set their initial values as $0$.
	
	Repeat for $ m = 1, \cdots, M$, 
	
	\Indp Repeat for $ l = 1,\cdots, L$, 
	
	\Indp Repeat for $ k = 1,\cdots, K$, 
	
	\Indp Repeat for $ j = 1,\cdots, J$
	
	\Indp Use $\{Z_h \}_{h=1}^{N} \backslash D_k $ to compute $(\widehat \eta_N^{\pi_{m}, \beta_l}, \widehat \alpha(\pi_m, \beta))$ and $\widehat{\mynu}(\pi_{m},\beta_l)$ by \eqref{lower level1 app}-\eqref{lower level2 app} and \eqref{lower level3 app}-\eqref{lower level4 app} using tuning parameters $(\mu^{(j)}_{1N}, \lambda^{(j)}_{1N})$ and $(\mu^{(j)}_{2N}, \lambda^{(j)}_{2N})$ respectively;
	
	Compute $\delta^{\pi_m, \beta_l}(\cdot; \widehat \eta^{\pi_m, \beta_l}, \widehat{Q}_N^{\pi_m, \beta_l})$ and $\Delta^{\pi_m}(\cdot; \widehat{H}_N^{\pi_m})$ and their corresponding squared Bellman errors $mse^{(1)}$ and $mse^{(2)}$ on the dataset	 $D_k$ by Gaussian kernel regression;
	
	Assign $e^{(1)}(m,l, j) = e^{(1)}(m,l, j) + mse^{(1)}$ and $e^{(2)}(m, j) = e^{(2)}(m, j) + mse^{(2)}$;
	
	\Indm \Indm \Indm
	
	Compute $j^{(1)\ast} \in \argmin_{j} \max_{m, l} e^{(1)}(m, l, j)$ and $j^{(2)\ast} \in \argmin_{j} \max_{m} e^{(2)}(m, j)$
	
	\textbf{Output:} $(\mu_{1N}^{j^{(1)\ast}}, \lambda_{1N}^{j^{(1)\ast}})$ and $(\mu_{2N}^{j^{(2)\ast}}, \lambda_{2N}^{j^{(2)\ast}})$.
\end{algorithm}
\subsection{Selection of Constant $c$ in $\Lambda^\pi_c$}
It is important to choose a proper constant $c$ in $\Lambda^\pi_c$ in order to protect against uncertainty in terms of the duration of use of the policy in future and different reference distributions. If we choose $c = 1$, \eqref{robust measure} in the main text becomes $R_{\min}$ for all $\pi \in \Pi$ and thus we are unable to distinct different policies because we are over conservative. If $c=0$, \eqref{robust measure} in the main text becomes the long-term average reward, where we  basically ignore any rewards happened in any finite period of time. If we know at least how long the learned policy will be implemented in the future, say $T_0$, and how fast the policy-induced Markov chain converges to the stationary distribution, we can choose $c$ properly. For example, we have the following uniform ergodic theorem given in Theorem 7.3.10 of \cite{hernandez2012further}.
\begin{theorem}[Uniform Geometric Ergodicity]\label{convergence theorem}
	If for any $\pi \in \Pi$, the induced Markov chain $P^\pi$ is $\psi$-irreducible and aperiodic and satisfies the geometric drift condition described in Theorem 7.3.1 of \citep{hernandez2012further}, then there exist constants $0 < \alpha(\pi) < 1$ and $C_0(\pi) >0$ such that,
	$$
	\max_{s \in \S} \left\|P^{\pi}_t(\bullet \, | \, S_1 = s) - d^\pi(\bullet) \right\|_{\text{TV}} \leq \min\left(1, C_0(\pi) \alpha^t(\pi)\right).
	$$	
\end{theorem}
Indeed, this theorem can further imply that for any $\tilde{\GG} \in \Lambda$
$$
\left\|\bar d^\pi_{T;\tilde{\GG}}(\bullet) - d^\pi(\bullet) \right\|_{\text{TV}} \leq \min\left(1,\frac{C_0(\pi)\alpha(\pi)}{1-\alpha(\pi)}\frac{1}{T}\right).
$$
As we can see, the average visiting distribution of the induced Markov chain converges sublinearly to the unique stationary distribution in terms of time $T$. If we  assume there exists some positive constants $\tilde \alpha < 1$ and $\tilde C_0$ independent of $\pi$ that the above inequality holds uniformly over $\Pi$, then we can choose $c$ based on these constants. For example, if we know $\tilde \alpha \leq 0.9$, $\tilde C_0\leq1$, and $T_0 =100$, then we can choose $c = \frac{0.9}{100(1-0.9)}= 0.09$, so that $\U^\pi_{T_0} \subseteq \U_c^\pi$, which satisfies our need. In practice, one also needs to consider the estimation error in terms of $c$. As we can see from \eqref{dual feasible set} in the main text, if we choose $c$ large, the set $\Lambda_c^\pi$ is large, thus requiring more data to estimate $\pi^\ast_c$ than that of smaller $c$ in order to achieve the same level of the accuracy. In contrast, a larger $c$ can guarantee a more robust policy than a smaller $c$ because we consider more uncertain scenarios. The remaining question is how to estimate $\tilde \alpha$ and $\tilde{C}_0$ using $D_n$, which we leave it as a future work. 
\section{Technical Proofs}
In this section, we provide all the technical proofs to the theoretical results in the main text. The notation $K(N) \leqconst L(N)$ (resp. $K(N) \geqconst L(N)$) means that there exist a sufficiently large constant (resp. small) constant $c_1>0$ (resp. $c_2>0$) such that $K(N) \geq c_1 L(N)$ (resp. $K(N) \leq c_2 L(N)$). Moreover, $K(N) \simeq L(N)$ means $K(N) \leqconst L(N)$ and $K(N) \geqconst L(N)$. All these constants do not depend on data, i.e., deterministic. For notational simplicity, we omit $\beta$, the auxiliary variable in the relative value function $Q^{\pi, \beta}$, its difference $U^{\pi, \beta}$, temporal difference $\delta^{\pi, \beta}$ and their related estimators when there is no confusion. Finally, we also denote $T=T_0$, $\mu_{jN} = \mu_N$, $\lambda_{jN} = \lambda_N$, $\F_j = \F$ and $\G_j = \G$ for the ease of presentation when there is no confusion. 

\begin{definition}[Covering Number]
	Let $\epsilon > 0$, $\F$ be a set of real-valued functions defined over some space $\X$. For a finite collection of $N_\epsilon = \{f_1, \cdots, f_{N_\epsilon}\}$ defined on $\X$ such that for every $f \in \F$, there exits a function $\bar f \in N_\epsilon$ that $\text{dist}(f - \bar f) \leq \epsilon$ is called an $\epsilon$-cover of $\F$ with respect to the metric $\text{dist}({\bullet})$. Let $\N(\epsilon, \F, \text{dist}({\bullet}))$ be the size of the smallest $\epsilon$-cover of $\F$ with respect to $\text{dist}({\bullet})$.
\end{definition}

\subsection{Proof of Theorem \ref{duality}}
It can be seen that the defined function $\phi(x)$ is convex. By the results in \citep{shapiro2017distributionally}~[section 3.2], we can show that
$$
\max_{u \in \Lambda^\pi_c} \, \, -\EE_{d^\pi}\left[\RR(S)\right] = \min_{\lambda \geq 0, \beta} \lambda c + \beta + \EE_{d^\pi}\left[\left(\lambda\phi\right)^\ast \left(-\RR(S) - \beta\right)\right],
$$
where the function $\left(\lambda\phi\right)^\ast (\bullet)$ refers to the conjugate of $\lambda \phi(\bullet)$. Note that we modify the left hand side above into a maximization problem to be consistent with results in \citep{shapiro2017distributionally}. Then by the definition of $\phi(x)$, we have that
\[ \left(\lambda\phi\right)^\ast (x)= \begin{cases} 
-\frac{\lambda}{2} + \left(x + \frac{\lambda}{2}\right)_+ & x \leq \frac{\lambda}{2} \\[0.1in]
+\infty & x > \frac{\lambda}{2}
\end{cases},
\]
where $\left(\bullet \right)_+ = \max(0, \bullet)$. Then we have the following equivalent formulation.
\begin{align*}
& \min_{\lambda \geq 0, \beta} \lambda c + \beta + \EE_{d^\pi}\left[\left(\lambda\phi\right)^\ast \left(-\RR(S) - \beta \right)\right] \\[0.1in]
=& \min_{\substack{\beta, \lambda \geq 0, \\ \lambda \geq -(2R_{\min} + 2\beta)}} \lambda c + \beta - \frac{\lambda}{2} +\EE_{d^\pi}\left[\left(-\RR(S) +\frac{\lambda}{2} - \beta \right)_+\right] \\[0.1in]
=& \min_{\substack{\beta, \lambda \geq 0, \\ \lambda \geq -(R_{\min} + \beta)}} \lambda c + \beta +\EE_{d^\pi}\left[\left(-\RR(S) - \beta \right)_+\right] \\[0.1in]
=& \min_{\lambda \geq 0} \lambda c - R_{\min} - \lambda +\EE_{d^\pi}\left[\left(-\RR(S)+R_{\min} + \lambda \right)_+\right] \\[0.1in]
=& - cR_{\min}  + \min_{\beta \geq R_{\min}} -\left(1-c\right)\beta +\EE_{d^\pi}\left[\left(-\RR(S)+\beta \right)_+\right] \\[0.1in]
=& - cR_{\min}  - \left(1-c\right)\max_{\beta \in \mathcal{R}} \left\{\beta -\frac{1}{\left(1-c\right)}\EE_{d^\pi}\left[\left(-\RR(S)+\beta \right)_+\right]\right\},
\end{align*}
where the first equality uses the definition of $\left(\lambda\phi\right)^\ast (x)$ and the assumption in this theorem, the second equality changes the variable $\beta \leftarrow (\beta - \frac{\lambda}{2})$, the third equality uses the monotonicity with respect with $\beta$, the fourth equality changes the variable $\beta \leftarrow (\lambda + R_{\min})$ and the last inequality is because the optimal solution is within the feasible set. Therefore, we have the first statement and the second statement follow immediately as below.
\begin{align*}
&\argmax_{\pi \in \Pi} \min_{u \in \Lambda^\pi_c}\EE_{d^\pi}\left[\RR(S)\right]\\[0.1in]
= &\argmax_{\pi \in \Pi} \max_{\beta \in \mathbb{R}} \left\{\beta -\frac{1}{\left(1-c\right)}\EE_{d^\pi}\left[\left(-\RR(S)+\beta \right)_+\right]\right\}.
\end{align*}
\subsection{Finite Sample Error Bound for the Relative Value  Difference Function}
\textbf{Proof of Theorem \ref{thm: value}} 
Denote $\bar B = [-R_{\max}, R_{\max}]$.
By Lemma \ref{lem: U bound} given by Assumption \ref{ass: stationarity}, we have
\begin{align*}
\sup_{\pi \in \Pi, \beta \in \bar B}\norm{\widehat U^{\pi,\beta}_N- U^{\pi,\beta}} \leq \left(1 +  \frac{1}{p_{\min}}\right)\sup_{\pi \in \Pi, \beta \in \bar B}\norm{\widehat Q_N^{\pi,\beta} - Q^{\pi,\beta}}.
\end{align*}
By the definition of $\widehat U^{\pi,\beta}$, we can assume that the expectation of $\widehat Q_N^{\pi,\beta}$ under the stationary distribution is 0, otherwise we can shift $\widehat Q_N^{\pi,\beta}$ by a constant to obtain it. Then we apply Lemma \ref{lemma:contraction} and Lemma \ref{lemma: upper bound of eta, Q} below to get
\begin{align*}
\sup_{\pi \in \Pi, \beta \in \bar B}\norm{\widehat U^{\pi,\beta}_N- U^{\pi,\beta}} &\lesssim 2\left(1 +  \frac{1}{p_{\min}}\right)\sup_{\pi \in \Pi, \beta \in \bar B}\norm{(\I - \P^\pi)\left(\widehat Q_N^{\pi,\beta} - Q^{\pi,\beta} \right)}\\
& \lesssim 2\left(1 +  \frac{1}{p_{\min}}\right)(1 + \sqrt{1 + \sup_{\pi \in \Pi}\sigma_\pi^2}) \sup_{\pi \in \Pi, \beta \in \bar B}\norm{\E_\pi(\widehat{\eta}_N^{\pi,\beta}, \widehat Q_N^{\pi,\beta})},
\end{align*}
where $\sigma^2_\pi$ is the variance of $\omega^\pi$and $\E_{\pi,\beta}$ is the bellman error, i.e.,
$$
\E_{\pi,\beta}(s, a; \eta, Q) \triangleq  \EE\left[\beta - \frac{1}{1-c}\left(\beta - R_{t}\right)_+ + \sum_{a'} \pi(a'|S_{t+1}) Q(S_{t+1}, a') - \eta - Q(s, a) \given S_t = s, A_t = a\right].
$$
Since $\norm{w^\pi}^2 = 1 + \sigma_\pi^2$ and by Assumption \ref{ass: policy class}~(d), $\sup_{\pi \in \Pi}\sigma_\pi^2 < \infty$.

Next, we derive the uniform error bound for $\sup_{\pi \in \Pi, \beta \in \bar B}\norm{\E_{\pi,\beta} (\widehat{\eta}_N^{\pi,\beta}, \widehat Q_N^{\pi,\beta})}$. Let $$\T_\pi(s, a; Q) = \EE\left[\beta - \frac{1}{1-c}\left(\beta - R_{t}\right)_+ + \sum_{a'} \pi(a'|S_{t+1}) Q(S_{t+1}, a') \given S_t = s, A_t = a\right].
$$
By the definition of $\kappa$ in Assumption \ref{assumption: value function}(c),
\begin{align}
& \sup_{\pi \in \Pi, \beta \in \bar B}\norm{\E_{\pi,\beta} (\widehat{\eta}_N^{\pi,\beta}, \widehat Q_N^{\pi,\beta})}^2 =\sup_{\pi \in \Pi, \beta \in \bar B}\norm{\T_\pi(\cdot, \cdot; \widehat Q_N^{\pi,\beta}) - \widehat \eta^{\pi,\beta}_N - \widehat Q_N^{\pi,\beta}(\cdot, \cdot)}_{}^2  \notag \\
& \leq \frac{1}{\kappa^2} \sup_{\pi \in \Pi, \beta \in \bar B}\norm{\g(\widehat \eta_N^{\pi,\beta}, \widehat Q_N^{\pi,\beta})}_{}^2  
\leq \frac{2}{\kappa^2} \left(\sup_{\pi \in \Pi, \beta \in \bar B}\norm{\g(\widehat \eta_N^{\pi,\beta}, \widehat Q_N^{\pi,\beta}) - \widehat g^{\pi,\beta}_{N}(\widehat \eta_N^{\pi,\beta}, \widehat Q_N^{\pi,\beta})}_{}^2 + \sup_{\pi \in \Pi, \beta \in \bar B}\norm{\widehat g^{\pi,\beta}_{N}(\widehat \eta_N^{\pi,\beta}, \widehat Q_N^{\pi,\beta})}_{}^2\right). \label{decom}
\end{align}
We first consider the second term in the RHS of the above inequality. Using Lemma \ref{lemma: second term value} and letting $\tau \leq \frac{1}{3}$, with $N$ sufficiently large, the following holds with probability at least $1-\delta$:
\begin{align*}
& \sup_{\pi \in \Pi, \beta \in \bar B}\norm{\widehat g_N^\pi(\widehat \eta_N^{\pi,\beta}, \widehat Q_N^{\pi,\beta})}_{}^2  \\
& \leqconst  \mu_N+  \mu_N \sup_{\pi \in \Pi, \beta \in \bar B}J_2^2 \left\{g^\ast_{\pi,\beta}(\eta^{\pi,\beta}, \tilde{Q}^{\pi,\beta})\right\}+ (\mu_N+\lambda_N) \sup_{\pi \in \Pi, \beta \in \bar B}J_1^2(\tilde{Q}^{\pi,\beta}) +
\frac{(VC(\Pi)+1)\left[\log\left(\max(1/\delta, N)\right)\right]^{\frac{1}{\tau}}}{N} \\
&+ N^{-\frac{1-(2+\alpha)\tau}{1+\alpha-\tau(2+\alpha)}} + \frac{1+VC(\Pi)}{N \mu_N^{\alpha/(1-\tau(2+\alpha))}} + \frac{(VC(\Pi)+1)\log^{\frac{\alpha/\tau}{1+\alpha - \tau(2+\alpha)}}(\max(N,1/\delta))}{N\mu_N^{\alpha/(1+\alpha- (2+\alpha)\tau)}}  + \frac{1}{N \lambda_N^{\frac{\alpha}{1 - \tau(2+ \alpha)}}}\\
& \leqconst  \mu_N\left(1 +  \sup_{\pi \in \Pi, \beta \in \bar B} J_1^2(\tilde{Q}^{\pi,\beta})\right) +
\frac{(VC(\Pi)+1)\left[\log\left(\max(1/\delta, N)\right)\right]^{\frac{1}{\tau}}}{N} + N^{-\frac{1-(2+\alpha)\tau}{1+\alpha-\tau(2+\alpha)}}\\
& + \frac{1+VC(\Pi)}{N \mu_N^{\alpha/(1-\tau(2+\alpha))}} + \frac{(VC(\Pi)+1)\log^{\frac{\alpha/\tau}{1+\alpha - \tau(2+\alpha)}}(\max(N,1/\delta))}{N\mu_N^{\alpha/(1+\alpha- (2+\alpha)\tau)}},
\end{align*}
where we use the condition that $\lambda_N \simeq \mu_N$ and Assumption \ref{assumption: value function}(d).

We now turn to the first term. By Lemma \ref{lemma: uniform l2 bound value} with the same $\tau$ used above and $N$ sufficiently large, we have at least probability $1-\delta$, 
\begin{align*}
& \sup_{\pi \in \Pi, \beta \in \bar B}\norm{\widehat g^\pi_{N} (\widehat \eta_N^{\pi,\beta}, \widehat Q_N^{\pi,\beta}) - \g (\widehat \eta_N^\pi, \widehat Q_N^\pi)}_{}^2 && \leqconst \mu_N + \mu_N \sup_{\pi \in \Pi, \beta \in \bar B}J_2^2 \left\{\g(\widehat \eta_N, \widehat Q_N^\pi)\right\}+ \mu_N \sup_{\pi \in \Pi, \beta \in \bar B}J_1^2(\widehat Q_N^\pi) \\
& &&+
\frac{(VC(\Pi)+1)\left[\log\left(\max(1/\delta, N)\right)\right]^{\frac{1}{\tau}}}{N} 
+ 	\frac{1+ VC(\Pi)}{N \mu_N^{\frac{\alpha}{1 - \tau(2+ \alpha)}}}.
%\\& \leq c_1 (1 + J_1^2(\widehat Q_N^\pi) +  J_2^2(\g (\widehat \eta_n^\pi, \widehat Q_N^\pi)) + \log(1/\delta))\mu_N
\end{align*}
Using Assumption \ref{assumption: value function}(d) again, this can be further bounded by
\begin{align}
\sup_{\pi \in \Pi, \beta \in \bar B}\norm{\widehat g^{\pi,\beta}_{N} (\widehat \eta_N^{\pi,\beta}, \widehat Q_N^{\pi,\beta}) - g^\ast_{\pi,\beta} (\widehat \eta_N^{\pi,\beta}, \widehat Q_N^{\pi,\beta})}_{}^2 & \leqconst \mu_N (1+ \sup_{\pi \in \Pi, \beta \in \bar B}J_2^2(\widehat Q^{\pi,\beta}_N)) + 
\frac{(VC(\Pi)+1)\left[\log\left(\max(1/\delta, N)\right)\right]^{\frac{1}{\tau}}}{N} \nonumber \\[0.1in]
& + 	\frac{1+ VC(\Pi)}{N \mu_N^{\frac{\alpha}{1 - \tau(2+ \alpha)}}}  \label{first term}
%\\ c_1 (1 + 2C_1^2 + (1 + 2C_2^2)J_1^2(\widehat Q_N^\pi) + \log(1/\delta))\mu_N 
\end{align} 
To bound $\sup_{\pi \in \Pi, \beta \in \bar B}J_1^2(\widehat Q_N^{\pi,\beta})$, the optimizing property of the estimators $(\widehat \eta_N^{\pi,\beta},\widehat Q_N^{\pi,\beta})$ in \eqref{lower level1 app} implies that 
\begin{align*}
\lambda_N J_1^2(\widehat Q_N^{\pi,\beta}) 
& \leq \Pn \left[\frac{1}{T} \sum_{t=1}^{T} \widehat g^{\pi,\beta}_{N}(S_t, A_t; \widehat \eta_N^{\pi,\beta}, \widehat Q_N^{\pi,\beta})^2\right] + \lambda_N J_1^2(\widehat Q_N^{\pi,\beta}) \\
& \leq \Pn \left[\frac{1}{T} \sum_{t=1}^{T} \widehat g^{\pi,\beta}_{N}(S_t, A_t; \eta^{\pi,\beta}, \tilde Q^{\pi,\beta})^2\right] + \lambda_N J_1^2(\tilde Q^{\pi,\beta})\\
& = \Pn \left[\frac{1}{T} \sum_{t=1}^{T} (\widehat g^{\pi,\beta}_{N}(S_t, A_t; \eta^{\pi,\beta}, \tilde Q^{\pi,\beta}) - g^\ast_{\pi,\beta}(S_t, A_t; \eta^{\pi,\beta}, \tilde Q^{\pi,\beta}))^2\right] + \lambda_N J_1^2(\tilde Q^{\pi,\beta})\\
& \leqconst \mu_N (1+J_1^2(\tilde Q^{\pi,\beta}))+
\frac{(VC(\Pi)+1)\left[\log\left(\max(1/\delta, N)\right)\right]^{\frac{1}{\tau}}}{N} 
+ 	\frac{1+ VC(\Pi)}{N \mu_N^{\frac{\alpha}{1 - \tau(2+ \alpha)}}}+ \lambda_N J_1^2(\tilde Q^\pi),
\end{align*}
where we use $g_{\pi,\beta}^\ast(\eta^{\pi,\beta}, \tilde Q^{\pi,\beta}) = 0$ in the third line and the last inequality follows by Lemma \ref{lemma: uniform l2 bound value} and the fact that $J_2(g^\ast_{\pi,\beta}(\eta^{\pi,\beta}, \tilde Q^{\pi,\beta})) = 0$. As a result, we have
\begin{align*}
\sup_{\pi \in \Pi, \beta \in \bar B} J_1^2(\widehat Q_N^{\pi,\beta})  &\leqconst \sup_{\pi \in \Pi, \beta \in \bar B}J_1^2(\tilde Q^{\pi,\beta}) + \frac{\mu_N}{\lambda_N} (1+\sup_{\pi \in \Pi, \beta \in \bar B}J_1^2(\tilde Q^{\pi,\beta})) \\[0.1in]
& +
\frac{(VC(\Pi)+1)\left[\log\left(\max(1/\delta, N)\right)\right]^{\frac{1}{\tau}}}{N\lambda_N} 
+ 	\frac{1+ VC(\Pi)}{\lambda_N N \mu_N^{\frac{\alpha}{1 - \tau(2+ \alpha)}}}.
\end{align*}
Combining with (\ref{first term}) and recalling that $\lambda_N \simeq \mu_N$ give
\begin{align*}
\sup_{\pi \in \Pi, \beta \in \bar B} \norm{\widehat g^{\pi,\beta}_{N} (\widehat \eta_N^{\pi,\beta}, \widehat Q_N^{\pi,\beta}) - g^\ast_{\pi,\beta} (\widehat \eta_N^{\pi,\beta}, \widehat Q_N^{\pi,\beta})}_{}^2  & \leqconst \mu_N (1+\sup_{\pi \in \Pi, \beta \in \bar B}J_1^2(\tilde Q^{\pi,\beta})) \\[0.1in]
& +
\frac{(VC(\Pi)+1)\left[\log\left(\max(1/\delta, N)\right)\right]^{\frac{1}{\tau}}}{N} 
+ 	\frac{1+ VC(\Pi)}{N \mu_N^{\frac{\alpha}{1 - \tau(2+ \alpha)}}} 
%\\ & + \mu_N \left(  1+J_1^2(\tilde Q^\pi) +
%\frac{\left[\log\left(\max(1/\delta, N)\right)\right]^{\frac{1}{\tau}}}{N\mu_N} 
%+ 	\frac{1+ VC(\Pi)}{\mu_N N \mu_N^{\frac{\alpha}{1 - \tau(2+ \alpha)}}} \right).
%\\ c_1 (1 + 2C_1^2 + (1 + 2C_2^2)J_1^2(\widehat Q_N^\pi) + \log(1/\delta))\mu_N 
\end{align*}
Summarizing together, we can show that for sufficiently large $N$ and if $\tau \leq \frac{1}{3}$, then with probability at least $1-2\delta$, we have
\begin{align*}
&\sup_{\pi \in \Pi, \beta \in \bar B}\norm{\T_{\pi,\beta}(\cdot, \cdot; \widehat Q_N^{\pi,\beta}) - \widehat \eta^{\pi,\beta}_N - \widehat Q_N^{\pi,\beta}(\cdot, \cdot)}_{}^2\\ & \leqconst  \mu_N+  \mu_N \sup_{\pi \in \Pi, \beta \in \bar B}J_1^2  (\tilde{Q}^{\pi,\beta})+
\frac{(VC(\Pi)+1)\left[\log\left(\max(1/\delta, N)\right)\right]^{\frac{1}{\tau}}}{N} \\
&+ N^{-\frac{1-(2+\alpha)\tau}{1+\alpha-\tau(2+\alpha)}} + \frac{1+VC(\Pi)}{N \mu_N^{\alpha/(1-\tau(2+\alpha))}} + \frac{(VC(\Pi)+1)\log^{\frac{\alpha/\tau}{1+\alpha - \tau(2+\alpha)}}(\max(N,1/\delta))}{N\mu_N^{\alpha/(1+\alpha- (2+\alpha)\tau)}}
%& \leqconst  \mu_N\left(1 +  \sup_{\pi \in \Pi, \beta \in \bar B} J_1^2(\tilde{Q}^\pi)\right) +
%\frac{(VC(\Pi)+1)\left[\log\left(\max(1/\delta, N)\right)\right]^{\frac{1}{\tau}}}{N} + N^{-\frac{1-(2+\alpha)\tau}{1+\alpha-\tau(2+\alpha)}}\\
%& + \frac{1+VC(\Pi)}{N \mu_N^{\alpha/(1-\tau(2+\alpha))}} + \frac{(VC(\Pi)+1)\log^{\frac{\alpha/\tau}{1+\alpha - \tau(2+\alpha)}}(\max(N,1/\delta))}{N\mu_N^{\alpha/(1+\alpha- (2+\alpha)\tau)}}
\end{align*}
To conclude our proof, we discuss how to choose $\mu_N$ and $\tau$ to obtain a reasonable upper bound. Observe the RHS of the above bound, we can see that when $\mu_N$ converges to $0$, the last term will decay faster than the last but the second term. Then we fix $\tau$ and let 
\begin{align*}
\mu_N \left(1 + \sup_{\pi \in \Pi, \beta \in \bar B} J_1^2(\tilde{Q}^{\pi,\beta})\right) = \frac{1+VC(\Pi)}{N \mu_N^{\alpha/(1-\tau(2+\alpha))}},
\end{align*}
which gives us that
$$
\mu_N = \left[\frac{1+VC(\Pi)}{N(1+ \sup_{\pi \in \Pi, \beta \in \bar B}J_1^2(\tilde{Q}^{\pi,\beta}))}\right]^{\frac{1-\tau(2+\alpha)}{1+\alpha - \tau(2+ \alpha)}}.
$$
Plugging into the bound, we can have
\begin{align*}
&\mu_N\left(1 +  \sup_{\pi \in \Pi, \beta \in \bar B} J_1^2(\tilde{Q}^{\pi,\beta})\right) +
\frac{(VC(\Pi)+1)\left[\log\left(\max(1/\delta, N)\right)\right]^{\frac{1}{\tau}}}{N} + N^{-\frac{1-(2+\alpha)\tau}{1+\alpha-\tau(2+\alpha)}}\\
&  +\frac{1+VC(\Pi)}{N \mu_N^{\alpha/(1-\tau(2+\alpha))}} + \frac{(VC(\Pi)+1)\log^{\frac{\alpha/\tau}{1+\alpha - \tau(2+\alpha)}}(\max(N,1/\delta))}{N\mu_N^{\alpha/(1+\alpha- (2+\alpha)\tau)}}\\
&\leqconst \frac{\left[(1+ \sup_{\pi \in \Pi, \beta \in \bar B}J_1^2(\tilde{Q}^{\pi,\beta}))\right]^{\frac{\alpha}{1+\alpha - \tau(2+ \alpha)}}\left(1+VC(\Pi)\right)^{\frac{1-\tau(2+\alpha)}{1+\alpha - \tau(2+ \alpha)}}}{N^{\frac{1 - \tau(2+ \alpha)}{1 + \alpha - \tau(2+ \alpha)}}}\\
&+ \frac{(VC(\Pi)+1)\left[\log\left(\max(1/\delta, N)\right)\right]^{\frac{1}{\tau}}}{N} + N^{- \frac{1 - \tau(2+ \alpha)}{1 + \alpha - \tau(2+ \alpha)}} \\
& + \frac{\left(1+ \sup_{\pi \in \Pi, \beta \in \bar B}J_1^2(\tilde{Q}^{\pi,\beta})\right)^{\frac{\alpha(1-\tau(2+\alpha))}{(1+\alpha - \tau(2 + \alpha))^2}}\left[\log(\max(N, 1/\delta))\right]^\frac{\alpha}{\tau\left(1+\alpha - \tau(2+\alpha)\right)}}{\left(1+ VC(\Pi)\right)^{\frac{\alpha(1-\tau(2+\alpha))}{(1+\alpha - \tau(2 + \alpha))^2}}N^{1 - \frac{\alpha(1-\tau(2+\alpha))}{(1+\alpha - \tau(2 + \alpha))^2}}}(VC(\Pi)+1)\\
& \leqconst \frac{\left[(1+ \sup_{\pi \in \Pi, \beta \in \bar B}J_1^2(\tilde{Q}^{\pi,\beta}))\right]^{\frac{\alpha}{1+\alpha - \tau(2+ \alpha)}}\left(1+VC(\Pi)\right)^{\frac{1-\tau(2+\alpha)}{1+\alpha - \tau(2+ \alpha)}}}{N^{\frac{1 - \tau(2+ \alpha)}{1 + \alpha - \tau(2+ \alpha)}}}\\
&+ \frac{(VC(\Pi)+1)\left[\log\left(\max(1/\delta, N)\right)\right]^{\frac{1}{\tau}}}{N} \\
&+ (VC(\Pi)+1)\frac{\left(1+ \sup_{\pi \in \Pi, \beta \in \bar B}J_1^2(\tilde{Q}^{\pi,\beta})\right)^{\frac{\alpha(1-\tau(2+\alpha))}{(1+\alpha - \tau(2 + \alpha))^2}}\left[\log(\max(N, 1/\delta))\right]^\frac{\alpha}{\tau\left(1+\alpha - \tau(2+\alpha)\right)}}{\left(1+ VC(\Pi)\right)^{\frac{\alpha(1-\tau(2+\alpha))}{(1+\alpha - \tau(2 + \alpha))^2}}N^{1 - \frac{\alpha(1-\tau(2+\alpha))}{(1+\alpha - \tau(2 + \alpha))^2}}}\\
& \leqconst \frac{\left[(1+ \sup_{\pi \in \Pi, \beta \in \bar B}J_1^2(\tilde{Q}^{\pi,\beta}))\right]^{\frac{\alpha}{1+\alpha - \tau(2+ \alpha)}}\left(1+VC(\Pi)\right)}{N^{\frac{1 - \tau(2+ \alpha)}{1 + \alpha - \tau(2+ \alpha)}}}\\
&+ \frac{\left(1+VC(\Pi)\right)\left[\log\left(\max(1/\delta, N)\right)\right]^{\frac{1}{\tau}}}{N} \\
&+ \frac{\left(1+VC(\Pi)\right)\left(1+ \sup_{\pi \in \Pi, \beta \in \bar B}J_1^2(\tilde{Q}^{\pi,\beta})\right)^{\frac{\alpha(1-\tau(2+\alpha))}{(1+\alpha - \tau(2 + \alpha))^2}}\left[\log(\max(N, 1/\delta))\right]^\frac{\alpha}{\tau\left(1+\alpha - \tau(2+\alpha)\right)}}{N^{1 - \frac{\alpha(1-\tau(2+\alpha))}{(1+\alpha - \tau(2 + \alpha))^2}}}.\\
\end{align*}
To minimize the RHS of the above bound, we first consider
\begin{align*}
&\frac{\left[1+ \sup_{\pi \in \Pi, \beta \in \bar B}J_1^2(\tilde{Q}^{\pi,\beta})\right]^{\frac{\alpha}{1+\alpha - \tau(2+ \alpha)}}}{N^{\frac{1 - \tau(2+ \alpha)}{1 + \alpha - \tau(2+ \alpha)}}} \\
= & \frac{\left(1+ \sup_{\pi \in \Pi, \beta \in \bar B}J_1^2(\tilde{Q}^{\pi,\beta})\right)^{\frac{\alpha(1-\tau(2+\alpha))}{(1+\alpha - \tau(2 + \alpha))^2}}\left[\log(\max(N, 1/\delta))\right]^\frac{\alpha}{\tau\left(1+\alpha - \tau(2+\alpha)\right)}}{N^{1 - \frac{\alpha(1-\tau(2+\alpha))}{(1+\alpha - \tau(2 + \alpha))^2}}}.
\end{align*}
This is equivalent to letting
\begin{align}\label{rate equation}
\left[N\left(1+ \sup_{\pi \in \Pi, \beta \in \bar B}J_1^2(\tilde{Q}^{\pi,\beta})\right)\right]^{\frac{\alpha}{1+\alpha - \tau(2+ \alpha)}} = \left[\log(\max(N, 1/\delta))\right]^\frac{1}{\tau}.
\end{align}
Denote
\begin{align*}
A & = N\left(1+ \sup_{\pi \in \Pi, \beta \in \bar B}J_1^2(\tilde{Q}^{\pi,\beta})\right)\\
B & = \log(\max(N, 1/\delta)).
\end{align*}
Then we can obtain $\tau$ by solving
$$
\frac{\alpha}{1+\alpha - \tau(2+ \alpha)} \log(A) = \frac{1}{\tau} \log(B),
$$
which gives
$$
\tau = \frac{(1+\alpha)\log(B)}{\alpha \log(A) + (2+ \alpha)\log(B)}.
$$
Next, we consider
\begin{align*}
&\frac{\left[1+ \sup_{\pi \in \Pi, \beta \in \bar B}J_1^2(\tilde{Q}^{\pi,\beta})\right]^{\frac{\alpha}{1+\alpha - \tau(2+ \alpha)}}}{N^{\frac{1 - \tau(2+ \alpha)}{1 + \alpha - \tau(2+ \alpha)}}} = \frac{\left[\log(\max(N, 1/\delta))\right]^\frac{1}{\tau}}{N},
\end{align*}
which again gives us that
$$
\tau = \frac{(1+\alpha)\log(B)}{\alpha \log(A) + (2+ \alpha)\log(B)}.
$$
Based on these two observation, we will let 
$$
\tau = \frac{(1+\alpha)\log(B)}{\alpha \log(A) + (2+ \alpha)\log(B)}.
$$
Clearly, when $N$ is sufficiently large, $\log(A)$ dominates $\log(B)$ and then $\tau$ can be arbitrarily small, thus eventually satisfying $\tau \leq \frac{1}{3}$. In such case,  we can show that
\begin{align*}
&\norm{\T_{\pi,\beta}(\bullet, \bullet; \widehat Q_N^{\pi,\beta}) - \widehat \eta^{\pi,\beta}_N - \widehat Q_N^{\pi,\beta}(\bullet, \bullet)}_{}^2  \leqconst  \frac{\left[(1+ \sup_{\pi \in \Pi, \beta \in \bar B}J_1^2(\tilde{Q}^{\pi,\beta}))\right]^{\frac{\alpha}{1+\alpha - \tau(2+ \alpha)}}\left(1+VC(\Pi)\right)}{N^{\frac{1 - \tau(2+ \alpha)}{1 + \alpha - \tau(2+ \alpha)}}}\\
+&  \frac{(VC(\Pi)+1)\left[\log\left(\max(1/\delta, N)\right)\right]^{\frac{1}{\tau}}}{N} \\
& + (VC(\Pi)+1)\frac{\left(1+ \sup_{\pi \in \Pi, \beta \in \bar B}J_1^2(\tilde{Q}^{\pi,\beta})\right)^{\frac{\alpha(1-\tau(2+\alpha))}{(1+\alpha - \tau(2 + \alpha))^2}}\left[\log(\max(N, 1/\delta))\right]^\frac{\alpha}{\tau\left(1+\alpha - \tau(2+\alpha)\right)}}{N^{1 - \frac{\alpha(1-\tau(2+\alpha))}{(1+\alpha - \tau(2 + \alpha))^2}}}\\
\leqconst& \left(1+VC(\Pi)\right)(1+ \sup_{\pi \in \Pi, \beta \in \bar B}J_1^2(\tilde{Q}^{\pi,\beta}))^{\frac{\alpha}{1+\alpha}}\left[\log\left(\max(1/\delta, N)\right)\right]^{\frac{2+\alpha}{1+\alpha}} N^{-\frac{1}{1+\alpha}}.
\end{align*}
Correspondingly, we can choose 
$$
\mu_N \simeq \left(1+VC(\Pi)\right)(1+ \sup_{\pi \in \Pi, \beta \in \bar B}J_1^2(\tilde{Q}^{\pi,\beta}))^{-\frac{1}{1+\alpha}}\left[\log\left(\max(1/\delta, N)\right)\right]^{\frac{2+\alpha}{1+\alpha}} N^{-\frac{1}{1+\alpha}}.
$$
Putting all together, we can conclude that
$$
\sup_{\pi \in \Pi, \beta \in \bar B} \norm{\widehat U_N^{\pi,\beta}- U^{\pi,\beta}}^2 \leqconst \left(1+VC(\Pi)\right)(1+ \sup_{\pi \in \Pi, \beta \in \bar B}J_1^2(\tilde{Q}^{\pi,\beta}))^{\frac{\alpha}{1+\alpha}}\left[\log\left(\max(1/\delta, N)\right)\right]^{\frac{2+\alpha}{1+\alpha}} N^{-\frac{1}{1+\alpha}},
$$
with probability at least $1-3\delta$. Letting $\delta = \frac{1}{3N}$, we obtain the desired result.\\[0.1in]

Denote $U^{\pi}(Q) = Q(s, a) - \sum_{a' \in \cal A}\pi(a' | s') Q(s', a)$ and $U^{\pi,\beta}= U^\pi(Q^\pi)$.
\begin{lemma}\label{lem: U bound}
	Under Assumption \ref{ass: stationarity}, for any state-action function $Q$, we have 
	\begin{align*}
	\norm{U^{\pi}(Q)- U^{\pi,\beta}} \leq \left(1 +  \frac{1}{p_{\min}}\right)\norm{Q - Q^{\pi,\beta}}.
	\end{align*}
\end{lemma}
\textbf{Proof of Lemma \ref{lem: U bound}}
We omit $\beta$ for the ease of presentation in this proof.
\begin{align}
& \norm{U^\pi(Q) - U^\pi} = \sqrt{\EE[(1/T) \sum_{t=1}^{T} (U^\pi(S_t, A_t, S_{t+1};Q)- U^\pi(S_t, A_t, S_{t+1}))^2]} \notag \\
& \leq \sqrt{\EE[(1/T) \sum_{t=1}^{T} (Q(S_t, A_t)- Q^\pi(S_t, A_t))^2]} \notag \\& \qquad + \sqrt{\EE\Big[(1/T) \sum_{t=1}^{T} \big(\sum_{a} \pi(a|S_{t+1}) (Q(S_{t+1}, a)- Q^\pi(S_{t+1}, a))\big)^2\Big]}
\notag \\
& \leq \sqrt{\EE[(1/T) \sum_{t=1}^{T} (Q(S_t, A_t)- Q^\pi(S_t, A_t))^2]} \notag \\& \qquad + \sqrt{\EE\Big[(1/T) \sum_{t=1}^{T} \big(\sum_{a} \frac{\pi(a|S_{t+1})}{\pi_b(a|S_{t+1})}\pi_b(a|S_{t+1}) (Q(S_{t+1}, a)- Q^\pi(S_{t+1}, a))\big)^2\Big]}\notag\\ 
& \leq \sqrt{\EE[(1/T) \sum_{t=1}^{T} (Q(S_t, A_t)- Q^\pi(S_t, A_t))^2]} \notag \\& \qquad +\frac{1}{p_{\min}} \sqrt{\EE\Big[(1/T) \sum_{t=1}^{T} \big(\sum_{a} \pi_b(a|S_{t+1}) (Q(S_{t+1}, a)- Q^\pi(S_{t+1}, a))\big)^2\Big]}\notag\\
&= (1 + \frac{1}{p_{\min}})
\norm{Q - Q^\pi}
\label{key},
\end{align}
where the last inequality is based on $\pi_b(a | s) \geq p_{\min}$ for every $(s, a) \in \S \times \A$ and the last equality is based on the stationarity of the trajectory $\D$ given in Assumption \ref{ass: stationarity}.

\begin{lemma} \label{lemma:contraction}
	Suppose Assumption \ref{ass: policy class} (e) holds. Then for any state-action function $\tilde Q$ such that $d^\pi(\tilde Q) = 0$, we have $\norm{\tilde  Q - Q^{\pi,\beta}} \leq 2\left( 1+  C_4 \bar{\alpha}/(1-\bar{\alpha}) \right) \norm{(\I - \P^\pi) (\tilde Q - Q^{\pi,\beta})}$, for some constant $C_4$.
\end{lemma}
\textbf{Proof of Lemma \ref{lemma:contraction}}
We omit $\beta$ in $Q^{\pi,\beta}$ in this proof. Let $\P^\pi_t := (\P^\pi)^t$ be the $t$-step transition kernel. Choose a $t$ sufficiently large such that $C_4\bar{\alpha}_{t} \leq 1/2$ for every $\pi \in \Pi$, then we can obtain that
\begin{align*}
\norm{\tilde  Q - Q^\pi}_{{}} & \leq 	\norm{(\I - \P^\pi_t)(\tilde  Q - Q^\pi)}_{{}}  + \norm{\P^\pi_t(\tilde  Q - Q^\pi)}_{} \\ & \leq \norm{(\I - \P^\pi_t)(\tilde  Q - Q^\pi)}_{{}}  + C_3\bar{\alpha}^t \norm{\tilde  Q - Q^\pi}_{} \\
& \leq \norm{(\I - \P^\pi_t)(\tilde  Q - Q^\pi)}_{{}}  + (1/2) \norm{\tilde  Q - Q^\pi}_{} .
\end{align*}
This implies that
\begin{align*}
\norm{\tilde  Q - Q^\pi}_{{}} & \leq 2 \norm{(\I - \P^\pi_t)(\tilde  Q - Q^\pi)}_{{}}  \\
&= 2 \norm{(\I - \P_1^\pi + \P_1^\pi - \P_2^\pi + \cdots + \P^\pi_{t-1} - P^\pi_t)(\tilde  Q - Q^\pi)}_{{}} \\
& \leq 2 (\norm{(\I - \P_1^\pi )(\tilde  Q - Q^\pi)}_{{}}  + \norm{(\P_1^\pi - \P_2^\pi )(\tilde  Q - Q^\pi)}_{{}} + \cdots  \norm{(\P^\pi_{t-1} - P^\pi_t)(\tilde  Q - Q^\pi)}_{{}}).
\end{align*}
Denote $h = (\I - \P^\pi) (\tilde Q - Q^\pi)$. It can be seen that $d^\pi(h) = 0$. Now for each $k$, we can have  
\begin{align*}
\norm{(\P_{k-1}^\pi - \P^\pi_{k} )(\tilde Q - Q^\pi)}_{} = \norm{\P_{k-1}^\pi  (\I - \P^\pi)(\tilde Q - Q^\pi)}_{} = \norm{\P_{k-1}^\pi  h}_{} \leq C_0\norm{h}_{} \bar{\alpha}_{k-1},
\end{align*}
by again Assumption \ref{ass: policy class}~(e).
Hence 
\begin{align*}
\norm{\tilde  Q - Q^\pi}_{{}} 
& \leq 2 (\norm{h}_{} + C_0\norm{h}_{} \bar{\alpha} + C_0\norm{h}_{} \bar{\alpha}_{2} + \cdots + C_0\norm{h}_{} \bar{\alpha}_{t-1}) \\
& = 2C_0\norm{h} \sum_{k = 0}^{t-1} \bar{\alpha}_{k-1}\\
& \leq \tilde C_0\norm{h} ,
\end{align*}
for some constant $\tilde C_0$ since $t$ is uniformly bounded.
\begin{lemma}
	\label{lemma: upper bound of eta, Q}
	For all $(\eta, Q) \in \R \times \F$, $ \abs{\eta - M(\pi, \beta)} \leq \sqrt{1+\sigma_\pi^2} \norm{\E_{\pi,\beta}(\eta, Q)}$ and $\norm{(\I - \P^\pi)(Q -  Q^{\pi,\beta})} \leq (1+\sqrt{1+\sigma_\pi^2}) \norm{\E_{\pi,\beta}(\eta, Q)}$, where $\I$ is the identity operator. 
\end{lemma}
\textbf{Proof of Lemma \ref{lemma: upper bound of eta, Q}}
Denote $M(\beta, \pi) = \eta^\pi$ for the ease of presentation. The Bellman error can be written as
\begin{align*}
& \E_{\pi,\beta}(s, a; \eta, Q)  = \EE\left[\beta - \frac{1}{1-c}\left(\beta - R_{t}\right)_+ + \sum_{a'} \pi(a'|S_{t+1}) Q(S_{t+1}, a') - \eta - Q(s, a) \given S_t = s, A_t = a\right]\\
& =  (\eta^{\pi,\beta} - \eta) + ( Q^{\pi,\beta} - Q)(s, a) - \P^\pi ( Q^{\pi,\beta} - Q)(s, a)\\
& =  (\eta^{\pi,\beta} - \eta)e^\pi(s, a) + (\eta^{\pi,\beta} - \eta)(1-e^\pi(s, a))  + ( Q^{\pi,\beta} - Q)(s, a) - \P^\pi ( Q^{\pi,\beta} - Q)(s, a)\\
& =  (\eta^{\pi,\beta} - \eta)e^\pi(s, a) + (\eta^{\pi,\beta} - \eta)(H^\pi(s, a) - \P^\pi H^\pi(s, a))  + ( Q^{\pi,\beta} - Q)(s, a) - \P^\pi ( Q^{\pi,\beta} - Q)(s, a)\\
& =   (\eta^{\pi,\beta} - \eta)e^\pi(s, a) + w(s, a) -  \P^\pi w(s, a),
\end{align*}
where the fourth inequality is based on the bellman equation of the scaled ratio function and in the last equality we define $w =  Q^{\pi,\beta} - Q + (\eta^{\pi,\beta} - \eta) H^\pi$.
Using the orthogonality property of the stationary distribution, we have $\norm{\E_{\pi,\beta}(\eta, Q)}_{}^2 = (\eta - \eta^{\pi,\beta})^2 \norm{e^\pi}_{}^2 + \norm{(\I - \P^\pi)w}_{}^2$ and thus $\abs{\eta - \eta^\pi} \leq \norm{e^\pi}_{}^{-1} \norm{\E_\pi(\eta, Q)}$.  Furthermore,  we have
\begin{align*}
& \norm{(\I - \P^\pi)(Q -  Q^{\pi,\beta})} = \norm{\E_{\pi,\beta}(\eta, Q) + (\eta - \eta^{\pi,\beta})} \\
& \leq \norm{\E_{\pi,\beta}(\eta, Q)} + \abs{\eta - \eta^{\pi,\beta}} \leq (1+\norm{e^{\pi}}^{-1}) \norm{\E_{\pi,\beta}(\eta, Q)} .
\end{align*}
Note that by the definition of (scaled) ratio functions, $\norm{e^\pi} = \norm{\omega^\pi}/(1+\sigma^2_\pi) = (1+\sigma^2_\pi)^{-1/2}$ (since $\norm{\omega^\pi}^2 = 1 + \sigma_\pi^2$) and thus we have 
\begin{align*}
& \abs{\eta - \eta^{\pi,\beta}} \leq \sqrt{1+\sigma_\pi^2} \norm{\E_{\pi,\beta}(\eta, Q)}\\
& \norm{(\I - \P^\pi)(Q -  Q^\pi)} \leq (1+\sqrt{1+\sigma_\pi^2}) \norm{\E_{\pi,\beta}(\eta, Q)}.
\end{align*}

\begin{lemma}
	\label{lemma: uniform l2 bound value}
	Let $g^\ast_{\pi,\beta}(\eta, Q)$ be the projected Bellman error operator defined in (\ref{surrogate BelErr Op}) and $\widehat g^{\pi,\beta}{N}(\eta, Q)$ be the estimated  Bellman error defined in (\ref{lower level2 app}) with the tuning parameter $\mu_N$. Suppose Assumptions \ref{ass: stationarity}, \ref{ass: policy class}, Assumption \ref{assumption: function classes}, and \ref{assumption: value function} hold. For any $0 < \tau \leq \frac{1}{3}$ and sufficiently large $N$, with probability at least $1-\delta$, the following inequalities hold for all $\eta,\beta \in \bar B$, $Q \in  \F$ and $\pi \in \Pi$:
	\begin{align*}
	& \norm{\widehat g^{\pi,\beta}{N} (\eta, Q) - g^\ast_{\pi,\beta}(\eta, Q)}_{}^2 && \leqconst \mu_N + \mu_N J_2^2 \left\{\g(\eta, Q)\right\}+ \mu_N J_1^2(Q) \\
	& &&+
	\frac{(VC(\Pi)+1)\left[\log\left(\max(1/\delta, N)\right)\right]^{\frac{1}{\tau}}}{N} 
	+ 	\frac{1+ VC(\Pi)}{N \mu_N^{\frac{\alpha}{1 - \tau(2+ \alpha)}}},\\
	&  J_2^2(\widehat g^{\pi,\beta}_{N}(\eta, Q)) && \leqconst 1 + J_2^2 \left\{g^\ast_{\pi,\beta}(\eta, Q)\right\}+ J_1^2(Q) \\
	& &&+
	\frac{(VC(\Pi)+1)\left[\log\left(\max(1/\delta, N)\right)\right]^{\frac{1}{\tau}}}{N\mu_N} 
	+ 	\frac{1+VC(\Pi)}{N \mu_N^{\frac{1 - \tau(2+ \alpha) + \alpha}{1 - \tau(2+ \alpha)}}},\\
	& \norm{\widehat g^{\pi,\beta}_{N}(\eta, Q) - g^\ast_{\pi,\beta}(\eta, Q)}_{N}^2 && \leqconst \mu_N + \mu_N J_2^2 \left\{g^\ast_{\pi,\beta}(\eta, Q)\right\}+ \mu_N J_1^2(Q) \\
	& &&+
	\frac{(VC(\Pi)+1)\left[\log\left(\max(1/\delta, N)\right)\right]^{\frac{1}{\tau}}}{N} 
	+ 	\frac{1+ VC(\Pi)}{N \mu_N^{\frac{\alpha}{1 - \tau(2+ \alpha)}}}.
	\end{align*}
\end{lemma}

\textbf{Proof of Lemma \ref{lemma: uniform l2 bound value}}
We omit $\beta$ in $Q^{\pi,\beta}$, $U^{\pi,\beta}$ and their relative quantities for the ease of presentation. Notice that Assumption \ref{ass: stationarity} implies that $\{(S_{it}, A_{it})\}_{i\geq1, t\geq 1}$ is also exponentially $\boldsymbol{\beta}$-mixing. We start with decomposing the error as
\begin{align*}
& \norm{\widehat g_N^\pi (\eta, Q) - \g (\eta, Q) }_{}^2  
= \frac{1}{T}\sum_{t=1}^{T}\EE\left[\left\{\widehat g_N^\pi (S_t, A_t; \eta, Q) - \g(S_t, A_t;\eta, Q)\right\}^2\right] \\
& =  \frac{1}{T} \sum_{t=1}^{T}\EE\left[ \left\{\widehat g_N^\pi (S_t, A_t; \eta, Q) - \delta_t^\pi(\eta, Q) + \delta_t^\pi(\eta, Q) -  \g(S_t, A_t;\eta, Q)\right\}^2\right]\\
& = \frac{1}{T} \sum_{t=1}^{T} \EE\left[ \left\{\delta_t^\pi(\eta, Q) - \widehat g_N^\pi (S_t, A_t; \eta, Q)\right\}^2\right] + \frac{1}{T} \sum_{t=1}^{T} \EE\left[ \left\{\delta_t^\pi(\eta, Q) -  \g(S_t, A_t;\eta, Q)\right\}^2\right] \\
& \qquad \quad  +  \frac{2}{T} \sum_{t=1}^{T}\EE\left[ \left\{\widehat g_N^\pi (S_t, A_t; \eta, Q) - \delta_t^\pi(\eta, Q)\right\} \left\{\delta_t^\pi(\eta, Q) -  \g(S_t, A_t;\eta, Q)\right\}\right].
\end{align*}
Since $ \sum_{t=1}^{T}\EE\left[ \left\{\E_\pi(S_t, A_t; \eta, Q) -  \g(S_t, A_t;\eta, Q)\right\} g(S_t, A_t)\right] = 0$ for all $g \in \G$ due to the optimizing property of $\g$, the last term above can be simplified as
\begin{align*}
&  \frac{2}{T} \sum_{t=1}^{T} \EE\Big[ \bigdkh{\widehat g_N^\pi (S_t, A_t; \eta, Q) - \g(S_t, A_t;\eta, Q)  \\
	& \qquad \qquad +   \g(S_t, A_t;\eta, Q)- \delta_t^\pi(\eta, Q)} \bigdkh{\delta_t^\pi(\eta, Q) -  \g(S_t, A_t;\eta, Q)}\Big]\\
& = \frac{2}{T} \sum_{t=1}^{T} \EE\left[ \bigdkh{\g(S_t, A_t;\eta, Q)- \delta_t^\pi(\eta, Q)} \bigdkh{\delta_t^\pi(\eta, Q) -  \g(S_t, A_t;\eta, Q)}\right]\\
& = - \frac{2}{T} \sum_{t=1}^{T}\EE\left[ \left\{\delta_t^\pi(\eta, Q) -  \g(S_t, A_t;\eta, Q)\right\}^2\right].
\end{align*}
As a result, we can have
\begin{align*}
& \norm{\widehat g_N^\pi (\eta, Q) - \g (\eta, Q) }_{}^2 \\
& =   \EE\Big[ \frac{1}{T} \sum_{t=1}^{T}  \bigdkh{\delta_t^\pi(\eta, Q) -  \widehat g_N^\pi (S_t, A_t; \eta, Q)}^2  -  \bigdkh{\delta_t^\pi(\eta, Q) - \g(S_t, A_t;\eta, Q)}^2\Big].
\end{align*}
For $g_1, g_2 \in \G, \eta \in \R, Q \in  \Q, \pi \in \Pi, \beta \in \bar B$, we define the following two functions:
%and a functional,  $\mathbf{J}$, by 
\begin{align*}
& f_1^\pi(g_1, g_2, \eta, Q): (S, A, S') \mapsto  \left\{\delta^\pi(\eta, Q) - g_1(S, A)\right\}^2 -  \left\{\delta^\pi(\eta, Q) -  g_2(S, A)\right\}^2\\
& f_2^\pi(g_1, g_2, \eta, Q): (S, A, S') \mapsto  \left\{\delta^\pi(\eta, Q) - g_2(S, A)\right\}\left\{g_1(S, A) - g_2(S, A)\right\},
%	\\ & \mathbf{J}^2(g_1, g_2, Q) = J_2^2(g_1) + \frac{2}{3}J_2^2(g_2) + \frac{2}{3}J_1^2(Q)
\end{align*}
where the underlying distribution of $(S, A, S')$ is the same as $(S_t, A_t, S_{t+1})$. Recall that $\left\{S_t, A_t, S_{t+1}\right\}_{t = 1}^{T}$ is a stationary process by Assumption \ref{ass: stationarity}.

%$$f(g_1, g_2, \eta, Q): \D \mapsto \frac{1}{T} \sum_{t=1}^{T} (\delta_t^\pi(\eta, Q) - g_1(S_t, A_t))^2 -  (\delta_t^\pi(\eta, Q) -  g_2(S_t, A_t))^2$$ 
%and $\mathbf{J}^2(g_1, g_2, Q) = J_2^2(g_1) + (2/3)J_2^2(g_2) + (2/3)J_1^2(Q)$.  
With these notations, we know that
\begin{align*}
& \norm{\widehat g_N^\pi (\eta, Q) - \g (\eta, Q) }_{}^2 = \EE \left[ f_1^\pi\left\{\widehat g_N^\pi(\eta, Q), \g (\eta, Q), \eta, Q\right\}\right],  \\
& \norm{\widehat g_N^\pi (\eta, Q) - \g (\eta, Q) }_{N}^2 = \PN \left[f_1^\pi\left\{\widehat g_N^\pi(\eta, Q), \g (\eta, Q), \eta, Q\right\}  + 2 f_2^\pi\left\{\gn(\eta, Q), \g(\eta, Q), \eta, Q\right\}\right].
\end{align*}
In the following, we introduce the decomposition for each pair of $(\eta, Q)$: 
\begin{align*}
& \norm{\widehat g_N^\pi (\eta, Q) - \g (\eta, Q) }_{}^2 + \norm{\widehat g_N^\pi (\eta, Q) - \g (\eta, Q) }_{N}^2 + \mu_N J_2^2\left\{\widehat g_N^\pi(\eta, Q)\right\} \\
%& = P f_1^\pi\left\{\widehat g_N^\pi(\eta, Q), \g (\eta, Q), \eta, Q\right\} + \Pn  f_1^\pi\left\{\widehat g_N^\pi(\eta, Q), \g (\eta, Q), \eta, Q\right\} + \mu_N J_2^2\left\{\widehat g_N^\pi(\eta, Q)\right\} \\
& = I_1(\eta, Q) + I_2(\eta, Q),
\end{align*}
where 
\begin{align*}
& I_1(\eta, Q) = 3\PN f_1^\pi\left\{\widehat g_N^\pi(\eta, Q), \g (\eta, Q), \eta, Q\right\} +  \mu_N [3J_2^2\left\{\widehat g_N^\pi (\eta, Q)\right\} + 2J_2^2 \left\{\g(\eta, Q)\right\} + 2J_1^2(Q)]\\
& I_2(\eta, Q) =  (\PN + P) f_1^\pi\left\{\widehat g_N^\pi(\eta, Q), \g (\eta, Q), \eta, Q\right\} + \mu_N J_2^2\left\{\widehat g_N^\pi(\eta, Q)\right\} \\
& \qquad \qquad \qquad + 2 \PN f_2^\pi\left\{\gn(\eta, Q), \g(\eta, Q), \eta, Q\right\}- I_1(\eta, Q).
\end{align*}
For the first term, the optimizing property of $\widehat g_N^\pi(\eta, Q)$ implies that 
\begin{align*}
\frac{1}{3}I_1(\eta, Q)  
%	 & = \Pn f_1^\pi\left\{\widehat g_N^\pi(\eta, Q), \g (\eta, Q), \eta, Q\right\} + \mu_N \mathbf{J}^2(\widehat g_N^\pi (\eta, Q), \g(\eta, Q), Q )\\
& =  \PN \Big[  \left\{\delta_t^\pi(\eta, Q) - \widehat g_N^\pi (S, A; \eta, Q)\right\}^2 -  \left\{\delta^\pi(\eta, Q) -  \g(S, A;\eta, Q)\right\}^2\Big] \\
& \qquad + \mu_N J_2^2\left\{\widehat g_N^\pi (\eta, Q)\right\} + \frac{2}{3}\mu_N J_2^2 \left\{\g(\eta, Q)\right\} + \frac{2}{3}\mu_N J_1^2(Q)\\
& = \PN \Big[  \left\{\delta_t^\pi(\eta, Q) - \widehat g_N^\pi (S, A; \eta, Q)\right\}^2\Big]+ \mu_N J_2^2\left\{\widehat g_N^\pi (\eta, Q)\right\} \\
&  \qquad - \Pn \big[  \left\{\delta^\pi(\eta, Q) -  \g(S, A;\eta, Q)\right\}^2\big] + \frac{2}{3} \mu_N J_2^2 \left\{\g(\eta, Q)\right\}+ \frac{2}{3}  \mu_N J_1^2(Q)  \\
& \leq \frac{5}{3} \mu_N J_2^2 \left\{\g(\eta, Q)\right\}+ \frac{2}{3} \mu_N J_1^2(Q).
\end{align*}
Thus, $I_1(\eta, Q) \leq 5 \mu_N J_2^2 \left\{\g(\eta, Q)\right\}+ 2\mu_N J_1^2(Q)$ holds for all $(\eta, Q)$.

Next we derive the uniform bound of $I_2(\eta, Q)$ over all $(\eta, Q)$. We use the independent block techniques \citep{yu1994rates} and the peeling device with the exponential inequality for the relative deviation of the empirical process developed in \citep{farahmand2012regularized}. The key step is to develop an individualized independent block for each peeling component.

First of all, we  apply the peeling device. Note that $\EE [f_2^\pi\left\{\gn(\eta, Q), \g(\eta, Q), \eta, Q\right\}(S, A)] = 0$ and recall that the process $\{S_t, A_t\}_{t= 1}^T$ is stationary. We can then write $I_2(\eta, Q)$ as
\begin{align*}
& I_2 (\eta, Q)  
= (\PN + P) f_{1}^\pi\left\{\widehat g_N^\pi(\eta, Q), \g (\eta, Q), \eta, Q\right\} + \mu_N J_2^2\left\{\widehat g_N^\pi(\eta, Q)\right\}  \\
& \qquad\qquad\qquad + 2 \PN f_2^\pi\left\{\gn(\eta, Q), \g(\eta, Q), \eta, Q\right\} -  3 \PN f_1^\pi\left\{\widehat g_N^\pi(\eta, Q), \g (\eta, Q), \eta, Q\right\} \\
& \qquad\qquad\qquad  -  \mu_N [3J_2^2\left\{\widehat g_N^\pi (\eta, Q)\right\} + 2J_2^2 \left\{\g(\eta, Q)\right\} + 2J_1^2(Q)]\\
& = 2 (P - \PN) (f_1^\pi-f_2^\pi)\left\{\widehat g_N^\pi(\eta, Q), \g (\eta, Q), \eta, Q\right\} - P (f_1^\pi-f_2^\pi)\left\{\widehat g_N^\pi(\eta, Q), \g (\eta, Q), \eta, Q\right\}  \\
& \qquad - 2\mu_N [J_2^2\left\{\widehat g_N^\pi(\eta, Q)\right\} + J_2^2(\g(\eta, Q)) + J_1^2(Q)].
\end{align*} 
For simplicity, we denote $f^\pi = f_1^\pi - f_2^\pi$ by
\begin{align*}
f^\pi(g_1, g_2, \eta, Q): (S, A, S')\mapsto  (g_2-g_1)(S, A) \cdot \left(3\delta^\pi(\eta, Q) - 2g_2(S, A) - g_1(S, A)\right),
\end{align*}
and the functional $$\mathbf{J}^2(g_1, g_2, Q) = J_2^2(g_1) + J_2^2(g_2) + J_1^2(Q),$$ for any $g_1, g_2 \in \G$ and $Q \in \Q$. 
Fix some $t > 0$. 
\begin{align*}
& \Pr\left\{\exists (\beta, \pi, \eta, Q) \in \bar B \times \Pi \times \bar B  \times \Q,  I_2(\eta, Q) > t\right\}\\
& = \sum_{l=0}^{\infty} \Pr\Big(\exists (\beta, \pi, \eta, Q) \in \bar B \times \Pi \times \bar B  \times \Q,  ~ 2\mu_N \mathbf{J}^2\left\{\widehat g_N^\pi (\eta, Q), \g(\eta, Q), Q\right\} \in ~ [2^l t \indicator{l\neq 0}, 2^{l+1}t),   \\
& \hspace{14ex} 2(P-\PN) f^\pi\left\{\widehat g_N^\pi(\eta, Q), \g (\eta, Q), \eta, Q\right\} > P f^\pi\left\{\widehat g_N^\pi(\eta, Q), \g (\eta, Q), \eta, Q\right\} \\
& \hspace{20ex}+ 2\mu_N \mathbf{J}^2\left\{\widehat g_N^\pi (\eta, Q), \g(\eta, Q), Q\right\} + t\Big)\\
& \leq \sum_{l=0}^{\infty} \Pr\Big(\exists (\beta, \pi, \eta, Q) \in \bar B  \times \Pi \times \bar B  \times \Q,  ~ 2\mu_N \mathbf{J}^2\left\{\widehat g_N^\pi (\eta, Q), \g(\eta, Q), Q\right\} \leq 2^{l+1}t,   \\
& \hspace{14ex} 2(P-\PN) f^\pi\left\{\widehat g_N^\pi(\eta, Q), \g (\eta, Q), \eta, Q\right\} > P f^\pi\left\{\widehat g_N^\pi(\eta, Q), \g (\eta, Q), \eta, Q\right\} + 2^l t\Big)\\
& \leq \sum_{l=0}^\infty \Pr\left( \sup_{h \in \F_l}  \frac{(P-\PN) \left\{h(Z)\right\}}{P \left\{h(Z)\right\} + 2^l t} > \frac{1}{2} \right),
\end{align*}
where the function class $\F_l = \{ f^\pi\left\{g, \g (\eta, Q), \eta, Q\right\}: J_2^2(g) \leq \frac{2^{l} t}{\mu_N}, J_2^2(\g(\eta, Q)) \leq \frac{2^{l} t}{\mu_N},  J_1^2(Q) \leq \frac{2^{l} t}{\mu_N},  \eta \in \bar B , Q \in \Q, \pi \in \Pi, \beta \in \bar B  \}.$ In addition, it is also easy to see that  for any $h =  f^\pi\left\{g, \g (\eta, Q), \eta, Q\right\} \in \F_{l}$,
\begin{align}\label{C1 condition}
\norm{f^\pi\left\{g, \g (\eta, Q), \eta, Q\right\}}_\infty \leq  6G_{\max}(\frac{2}{1-c}R_{\max}+ 2Q_{\max} + 3G_{\max}) \triangleq K_1, 
\end{align}

Next, we bound each term of the above probabilities by using the independent block technique. We define a partition by dividing the index $\{1, \cdots, N\}$ into $2v_N$ blocks, where each block has an equal length $x_N$. The residual block is denoted by $R_N$, i.e., $\{(j-1)x_N+1, \cdots, (j-1)x_N + x_N \}_{j= 1}^{2v_N}$ and $R_N = \left\{2v_Nx_N +1, \cdots, N \right\}$. Then it can be seen that  $N- 2x_N < 2 v_N x_N \leq N$ and the cardinality $|R_N|< 2x_N$. 
%	We define $(x_N, v_N)$-independent block  of a stationary stochastic process $\{Z_t\}_{t=1}^{N}$ as $\{Z(H_j): H_j = \{(j-1)+1, \cdots, (j-1) + x_N\} \}_{j= 1}^{2v_N}$ and $R_N$. 

For each $l \geq 0$, we will use a different independent block sequence denoted by $(x_{N, l}, v_{N, l})$ with the residual $R_{l}$ and then optimize the probability bound by properly choosing $(x_{N, l}, v_{N, l})$ and $R_{l}$. More specifically, we choose
$$
x_{N, l} = \lfloor x'_{N, l} \rfloor \epc \mbox{and } \epc  v_{N, l} = \lfloor \frac{N}{2x_{N,l}} \rfloor,
$$
where $x'_{N, l} = (\frac{Nt}{VC(\Pi)+1})^\tau (2^l)^p$ and $v'_{N, l} = \frac{N}{2x'_{N, l}}$ with some positive constants $\tau$ and $p$  determined later. We require $\tau \leq p \leq \frac{1}{2 + \alpha} \leq \frac{1}{2}$. We also need $t \geq \frac{VC(\Pi)+1}{N}$ so that $x'_{N, l}\geq 1$. Suppose $N$ is sufficiently large such that
\begin{align}\label{sample size constraint for lemma of value function}
N \geq c_1 \triangleq 4 \times 8^2 \times K_1 \times (VC(\Pi)+1). 
%	\geq 4^{\frac{p}{1-p}}8^{\frac{1}{1-p}}K_1^{\frac{p}{1-p}}(VC(\Pi)+1) .
\end{align}

In the following, we consider two cases. The first case considers any $l$ such that $x'_{N,l} \geq \frac{N}{8(VC(\Pi)+1)}$. In this case, since $\tau \leq p$, we can show that $x'_{N, l} \leq (\frac{Nt2^l}{VC(\Pi)+1})^p$. Combining with the sample size requirement,  we can obtain that  $(\frac{Nt2^l}{VC(\Pi)+1}) \geq (\frac{\frac{N}{8}}{VC(\Pi)+1})^\frac{1}{p} \geq 4NK_1$. Then we can show that in this case,
$$
\frac{(P-\PN) \left\{h(Z)\right\}}{P \left\{h(Z)\right\} + 2^l t} \leq  \frac{2K_1}{2^lt} \leq  \frac{1}{2}.
$$
Therefore, when $t \geq \frac{(VC(\Pi)+1)}{N}$ and $x'_{N, l} \geq \frac{N}{8(VC(\Pi)+1)}$, 
$$
\Pr\left( \sup_{h \in \F_l}  \frac{(P-\PN) \left\{h(Z)\right\}}{P \left\{h(Z)\right\} + 2^l t} > \frac{1}{2} \right) = 0.
$$

The second case we consider is when $x'_{N, l} < \frac{N}{8(VC(\Pi)+1)}$. We apply the relative deviation concentration inequality for the exponential $\boldsymbol{\beta}$-mixing stationary process given in Theorem 4 of \cite{farahmand2012regularized}, which combined results in \cite{yu1994rates} and Theorem 19.3 in \cite{gyorfi2006distribution}. To use their results, it then suffices to verify conditions (C1)-(C5) in Theorem 4 of \cite{farahmand2012regularized} with $\F = \F_l$, $\epsilon = 1/2$ and $\eta = 2^{l} t$ to get an exponential inequality for each term in the summation. First of all, Condition (C1) has been verified in \eqref{C1 condition}.

For (C2), recall  $f^\pi  = f_1^\pi - f_2^\pi$ and thus
\begin{align*}
\EE[f^\pi\left\{g, \g (\eta, Q), \eta, Q\right\}^2] \leq 2\EE [f_1^\pi\left\{g, \g (\eta, Q), \eta, Q\right\}(S, A, S')^2] + 2\EE [f_2^\pi\left\{g, \g (\eta, Q), \eta, Q\right\}(S, A, S')^2].
\end{align*}
For the first term of RHS above:
\begin{align*}
& \EE [f_1^\pi\left\{g, \g (\eta, Q), \eta, Q\right\}(S, A, S')^2]  \\
& = \EE\left[\left\{ \left\{\delta_t^\pi(\eta, Q) - g(S, A)\right\}^2 -  \left\{\delta_t^\pi(\eta, Q) - \g(S, A; \eta, Q)\right\}^2\right\}^2\right] \\
& = \EE\left[ \left\{2\delta^\pi(\eta, Q) - g(S, A)  - \g(S, A; \eta, Q)\right\}^2 \left\{\g(S, A; \eta, Q) - g(S, A)\right\}^2\right] \\
& \leq \left\{2(\frac{2}{1-c}R_{\max} + 2 Q_{\max}) + 2G_{\max}\right\}^2   \EE[\left(\g(S, A; \eta, Q) - g(S, A)\right)^2] \\
& = 4\left(\frac{2}{1-c}R_{\max} + 2 Q_{\max} + G_{\max}\right)^2 \EE \left[f^\pi\left\{g, \g (\eta, Q), \eta, Q\right\}(S, A, S')\right],
\end{align*}
and the second term:
\begin{align*}
&  \EE [f_2^\pi\left\{g, \g (\eta, Q), \eta, Q\right\}(S, A, S')^2] \\
& = \EE\left[\Bigdkh{ \dkh{\delta_t^\pi(\eta, Q) - \g(S, A; \eta, Q)}\dkh{g(S, A) - \g(S, A; \eta, Q)}}^2\right] \\
& \leq  \EE\left[ \dkh{\delta^\pi(\eta, Q) - \g(S, A; \eta, Q)}^2\dkh{g(S, A) - \g(S, A; \eta, Q)}^2\right] \\
& \leq \left(\frac{2}{1-c}R_{\max} + 2 Q_{\max} + G_{\max}\right)^2 \EE[\left(\g(S, A; \eta, Q) - g(S, A)\right)^2] \\
&  = \left(\frac{2}{1-c}R_{\max} + 2 Q_{\max} + G_{\max}\right)^2 \EE [f^\pi\left\{g, \g (\eta, Q), \eta, Q\right\}(S, A, S')],
\end{align*}
where we use again the fact that $\EE [f_2^\pi\left\{\gn(\eta, Q), \g(\eta, Q), \eta, Q\right\}(S, A, S')] = 0$. Putting together, we have shown that
\begin{align*}
\EE[f^\pi\left\{g, \g (\eta, Q), \eta, Q\right\}(S, A, S')^2]
%	\\
%	& \leq (2\cdot 4 + 2)\left(2R_{\max} + 2 Q_{\max} + G_{\max}\right)^2 \EE [f_1^\pi\left\{g, \g (\eta, Q), \eta, Q\right\}(\D)] \\
& \leq K_2 \EE [f^\pi\left\{g, \g (\eta, Q), \eta, Q\right\}(S, A, S')],
\end{align*}
where $K_2 = \left(\frac{2}{1-c}R_{\max} + 2 Q_{\max} + G_{\max}\right)^2$. This shows that Condition (C2) is satisfied.

%$K_1 ＝ 6G_{\max}(2R_{\max}+ 2Q_{\max} + 3G_{\max})$ and  $K_2$
%6G_{\max}(2R_{\max}+ 2Q_{\max} + 3G_{\max})$. 
%There exists some $K_1, K_2$ (depending on $\pi_{\min}$, $R_{\max}, Q_{\max}$) such that $\norm{f}_\infty \leq K_1$ and $\EE\left[f(\D)^2\right] \leq K_2 \EE\left[f(\D)\right]$. 

To verify the condition  (C3), without loss of generality, we assume $K_1 \geq 1$. Otherwise, let $K_1 = \max(1, K_1)$. Then we know that $2K_1x_{N, l} \geq \sqrt{2K_1x_{N, l}}$ since $x_{N, l} \geq 1$. We need to verify $\sqrt{N} \epsilon\sqrt{1-\epsilon} \sqrt{\eta} \geq 1152K_1x_{N, l}$, or it suffices to have $\sqrt{N} \epsilon\sqrt{1-\epsilon} \sqrt{\eta} \geq 1152K_1x'_{N, l}$ since $x'_{N, l} \geq x_{N, l}$ by definition. Recall that $\epsilon = 1/2$ and $\eta = 2^{l} t$. To show this, it is enough to show that
$$
\sqrt{N}\frac{\sqrt{2}}{4}  \sqrt{2^lt} \geq 1152K_1(\frac{Nt2^l}{VC(\Pi)+1})^p,
$$
since $(\frac{Nt2^l}{VC(\Pi)+1})^p \geq x'_{N, l}$. Recall that $p \leq \frac{1}{2+\alpha}$, then it is sufficient to 
let $t \geq \frac{2304\sqrt{2}K_1}{N}\triangleq \frac{c_1'}{N}$ so that the above inequality holds for every $l \geq 0$.

Next we verify (C4) that $\frac{|R_l|}{N} \leq \frac{\epsilon\eta}{6K_1}$.  Recall that $|R_l| < 2 x_{N, l} \leq 2x'_{N, l} = 2(\frac{Nt}{VC(\Pi)+1})^\tau(2^l)^p$. So if $t \geq \frac{c_2}{N}$ for some positive constant $c_2$ depending on $K_1$, we can have
$$
\frac{\epsilon\eta}{6K_1} = \frac{2^lt}{12K_1} \geq \frac{2(\frac{Nt}{VC(\Pi)+1})^\tau (2^l)^p}{N} = \frac{2x'_{N,l}}{N}> \frac{|R_l|}{N}.
$$
In addition, $|R_l| \leq 2 x'_{N, l} < \frac{N}{2}$.

Lastly we verify condition (C5). Define $$\Q_M = \{c + U: \abs{c} \leq R_{\max}, U = Q(s, a)-\sum_{a' \in \cal A}\pi(a'|s')Q(s', a'), Q \in \Q, J_1(Q) \leq M\}$$
and $\G_M = \{g: g \in \G, J_2(g) \leq M\}$.
%and the standard result for the entropy number for bounded set in euclidean space
%$t > c_1 n^{-1}$ for some constant $c_1$ depending on $K_1$, i.e., $(Q_{\max}, R_{\max}, G_{\max})$. 
It is not hard to verify that with $M = \sqrt{\frac{2^l t}{\mu_N}}$,
\begin{align*}
&  \log\left(\mathcal{N}(\epsilon , \F_l, \norm{\cdot}_N) \right)\\
\leqconst &  \log\left(\cal \mathcal{N}(\epsilon, \Q_{M}, \norm{\cdot}_\infty) \mathcal{N} (\epsilon, \G_{M}, \norm{\cdot}_\infty) \right) + \log(\N(\epsilon, \Pi, d_{\Pi}(\bullet))\N(\epsilon, \bar B, \norm{\cdot}_\infty)),\\
\end{align*}
by Assumption \ref{ass: policy class}.
As a result of the entropy condition in Assumption \ref{assumption: function classes}~(d) and \ref{ass: policy class}~(d), let $t \geq \mu_N$, we have 
\begin{align*}
& \log  \mathcal{N}(\epsilon, \F_l, \norm{\cdot}_N) \\
& \leqconst \log \mathcal{N}(\epsilon, \Q_{M}, \norm{\cdot}_\infty) + 2 \log \mathcal{N}(\epsilon, \G_{M}, \norm{\cdot}_\infty) + \log \mathcal{N}\left\{\epsilon, \Pi, d_{\Pi}(\bullet)\right\} + \log \mathcal{N}\left\{\epsilon, \bar B, \norm{\cdot}_\infty\right\}\\
& \leqconst   \left(\frac{2^l t}{\mu_N}\right)^{\alpha} \epsilon^{-2\alpha} + (VC(\Pi)+1)\log\left(1/\epsilon\right)\\
& \leq  c_3(1+VC(\Pi))  \left(\frac{2^l t}{\mu_N}\right)^{\alpha} \epsilon^{-2\alpha},
\end{align*}
for some constant $c_3\geq 1$ and $VC(\Pi)$ is the VC-index of the policy class $\Pi$. 
Then Condition (C5) is satisfied if the following inequality holds for all $x \geq (2^ltx_{N, l})/8$, 
\begin{align*}
\frac{\sqrt{v_{N, l}} (1/2)^2 x}{96x_{N, l}\sqrt{2} \max(K_1, 2K_2)} 
&\geq \int_{0}^{\sqrt{x}} \sqrt{c_3(1+VC(\Pi)) } \left(\frac{2^l t}{\mu_N}\right)^{\alpha/2} \left(\frac{u}{2x_{N, l}}\right)^{-\alpha}  du \\
&= x_{N, l}^{\alpha} x^{\frac{1-\alpha}{2}} \sqrt{2^\alpha c_3(1+VC(\Pi)) } \left(\frac{2^l t}{\mu_N}\right)^{\alpha/2}.
\end{align*}
It is enough to guarantee that
$$
\frac{\sqrt{v_{N, l}} (1/2)^2 x}{96x_{N, l}\sqrt{2} \max(K_1, 2K_2)} \geq x_{N, l}^{\alpha} x^{\frac{1-\alpha}{2}} \sqrt{2^\alpha c_3(1+VC(\Pi))} \left(\frac{2^l t}{\mu_N}\right)^{\alpha/2}.
$$
After some algebra, we can check that the above inequality holds if for some constant $c_4$,
\begin{align*}
t \geq c_4 \frac{(x_{N, l})^{1+\alpha}}{v'_{N, l}2^l \mu_N^\alpha}(1+VC(\Pi)),
\end{align*}
or equivalently,
\begin{align*}
t \geq c_5 \frac{1 + VC(\Pi)}{N \mu_N^{\frac{\alpha}{1 - \tau(2+ \alpha)}}\left(2^l\right)^{\frac{1 - p(2+\alpha)}{1 - \tau(2 + \alpha)}}},
\end{align*}
by the definition that $x_{N, l} \leq x'_{N, l}$ and $v'_{N, l} \leq v_{N, l}$. To summarize, if for any $l\geq 0$,
$$
t \geq \mu_N + c_5 \frac{1+ VC(\Pi)}{N \mu_N^{\frac{\alpha}{1 - \tau(2+ \alpha)}}\left(2^l\right)^{\frac{1 - p(2+\alpha)}{1 - \tau(2 + \alpha)}}},
$$
then the entropy inequality in Condition (C5) above holds. Since $0 < \tau \leq p \leq \frac{1}{1 + 2\alpha}$, the right hand side is a non-increasing function of $l$. Then
as long as, 
$$
t \geq \mu_N + c_5 \frac{1+ VC(\Pi)}{N \mu_N^{\frac{\alpha}{1 - \tau(2+ \alpha)}}},
$$
Condition (C5) holds .

To summarize, the conditions (C1-C5) in Theorem 4 in \cite{farahmand2012regularized} with $\F = \F_l$, $\epsilon = 1/2$ and $\eta = 2^{l} t$ hold for every $l \geq 0$ when $t \geq c'_2n^{-1}(VC(\Pi)+1)$ for some constant $c'_2\geq1$ and $t \geq \mu_N + c_5 \frac{1+ VC(\Pi)}{N \mu_N^{\frac{\alpha}{1 - \tau(2+ \alpha)}}}$. Thus when $N\geq c_1$,
\begin{align*}
& \Pr\left\{\exists (\beta, \pi, \eta, Q) \in B \times \Pi \times B \times \Q,  I_2(\eta, Q) > t\right\} \\
& \leq \sum_{l=0}^\infty \Pr \left[\sup_{h \in \F_l}  \frac{(P-\PN) \left\{h(Z)\right\}}{P \left\{h(Z)\right\} + 2^l t} > \frac{1}{2} \right] \\
& \leq \sum_{l=0}^\infty 120 \exp\left\{-c_6\frac{v_{N,l }^{'2}t2^l}{N}\right\} + 2\beta_{x_{N,l}}v_{N, l} \\
& \leq \sum_{l=0}^\infty 120 \exp\left\{-c_6\frac{v_{N,l }^{'2}t2^l}{N}\right\} + 2\beta_0\exp\left(-\beta_1 x_{N, l} + \log v_{N, l}\right) ,
\end{align*}
where the last inequality is based on exponential decay given in Assumption \ref{ass: stationarity}.  When $t \geq \frac{(VC(\Pi)+1)\left(4/\beta_1\log(N)\right)^{1/\tau}}{N}$, we have $\log v_{N, l} \leq \frac{1}{2}\beta_1 x_{N, l}$ by using $x_{N, l}' \leq 2x_{N, l}$ and $v_{N, l} \leq N$. This will further imply that $2\beta_{x_{N,l}}v_{N, l} \leq 2\beta_0\exp\left(-\beta_1 x_{N, l}/2 \right)$. Then we will have 
\begin{align*}
& \Pr\left\{\exists (\beta, \pi, \eta, Q) \in B \times \Pi \times B \times \Q,  I_2(\eta, Q) > t\right\} \\
& \leq \sum_{l=0}^\infty 120 \exp\left\{-c_6\frac{v_{N,l }^{'2}t2^l}{N}\right\} + 2\beta_0\exp\left(-\beta_1 x_{N, l} + \log v_{N, l}\right)\\
& \leqconst \sum_{l=0}^\infty 120 \exp\left(-c_7(Nt)^{1-2\tau}(2l)^{1-2p}(VC(\Pi)+1)^{2\tau} \right)+ 2\beta_0\exp\left(-\beta_1 (\frac{Nt}{VC(\Pi)+1})^{\tau}(2^l)^p\right)\\
& \leq c_8 \exp\left(-c_9(Nt)^{1-2\tau}(VC(\Pi)+1)^{2\tau}\right) + c_{10}\exp\left(-c_{11} (\frac{Nt}{VC(\Pi)+1})^{\tau}\right).
\end{align*}
As long as $t$ satisfies all the above constraints, 
$$
I_2(\eta, Q) \leq \frac{1}{N}\left\{\left(\frac{\log(\frac{2c_9}{\delta})}{c_8}\right)^{\frac{1}{1-2\tau}}\right\}  +\frac{VC(\Pi)+1}{N} \left\{\left(\frac{\log(\frac{2c_{11}}{\delta})}{c_{10}}\right)^{\frac{1}{\tau}} \right\},
$$
with probability at least $1-\delta$.
Collecting all the conditions on $t$ and combining with the bound of $I_1(\eta, Q)$, we have shown that with probability at least $1-\delta$, the following holds for all $(\beta, \pi, \eta, Q) \in B \times \Pi \times B \times \Q$:
\begin{align*}
& \norm{\widehat g_N^\pi (\eta, Q) - \g (\eta, Q) }_{}^2 + \norm{\widehat g_N^\pi (\eta, Q) - \g (\eta, Q) }_{N}^2 + \mu_N J_2^2\left\{\widehat g_N^\pi(\eta, Q)\right\}  \\
& \leq \mu_N + 5 \mu_N J_2^2 \left\{\g(\eta, Q)\right\}+ 2\mu_N J_1^2(Q) + \frac{1}{N}\left\{\left(\frac{\log(\frac{2c_9}{\delta})}{c_8}\right)^{\frac{1}{1-2\tau}}\right\}  +\frac{VC(\Pi)+1}{N} \left\{\left(\frac{\log(\frac{2c_{11}}{\delta})}{c_{10}}\right)^{\frac{1}{\tau}} \right\} \\
& + c_5 \frac{1+ VC(\Pi)}{N \mu_N^{\frac{\alpha}{1 - \tau(2+ \alpha)}}} + \frac{(VC(\Pi)+1)\left(4/\beta_1\log(N)\right)^{1/\tau}}{N} + c'_2\frac{VC(\Pi)+1}{N}.
\end{align*}
Recall that we require $0 < \tau \leq p \leq \frac{1}{2+\alpha}$. Consider any $\tau \leq \frac{1}{3}$ and pick $p$ any value between $\tau$ and $\frac{1}{1+2\alpha}$. Then the bound above can be simplified as
\begin{align*}
& \norm{\widehat g_N^\pi (\eta, Q) - \g (\eta, Q) }_{}^2 + \norm{\widehat g_N^\pi (\eta, Q) - \g (\eta, Q) }_{N}^2 + \mu_N J_2^2\left\{\widehat g_N^\pi(\eta, Q)\right\}  \\
& \leqconst (1+\mu_N) J_2^2 \left\{\g(\eta, Q)\right\}+ \mu_N J_1^2(Q) +
\frac{(VC(\Pi)+1)\left[\log\left(\max(1/\delta, N)\right)\right]^{\frac{1}{\tau}}}{N} + 	\frac{1+ VC(\Pi)}{N \mu_N^{\frac{\alpha}{1 - \tau(2+ \alpha)}}}.
\end{align*}

\begin{lemma}\label{lemma: second term value}
	Suppose the conditions in Lemma \ref{lemma: uniform l2 bound value}  hold. Let $(\widehat \eta_N^{\pi, \beta}, \widehat Q_N^{\pi, \beta})$ be the estimator in \eqref{lower level1 app}-\eqref{lower level2 app} with tuning parameter $\lambda_N$, and $\widehat g_N^{\pi, \beta}(\eta, Q)$ be the estimated Bellman error operator with the tuning parameter $\mu_N$. Up to some constant that, for sufficiently large $N$, the following holds with probability at least $1-2\delta$: 
	\begin{align*}
	& \norm{\widehat g_N^{\pi, \beta}(\widehat \eta_N^{\pi, \beta}, \widehat Q_N^{\pi, \beta})}_{}^2 + \norm{\widehat g_N^{\pi, \beta}(\widehat \eta_N^{\pi, \beta}, \widehat Q_N^{\pi, \beta})}_{N}^2  \\
	& \leqconst \mu_N+  \mu_N J_2^2 \left\{g^\ast_{\pi, \beta}(\eta^{\pi, \beta}, \tilde{Q}^{\pi, \beta})\right\}+ (\mu_N + \lambda_N) J_1^2(\tilde{Q}^{\pi, \beta}) +
	\frac{(VC(\Pi)+1)\left[\log\left(\max(1/\delta, N)\right)\right]^{\frac{1}{\tau}}}{N} \\
	&+ N^{-\frac{1-(2+\alpha)\tau}{1+\alpha-\tau(2+\alpha)}} + \frac{1+VC(\Pi)}{N \mu_N^{\alpha/(1-\tau(2+\alpha))}} + \frac{(VC(\Pi)+1)\log^{\frac{\alpha/\tau}{1+\alpha - \tau(2+\alpha)}}(\max(N,1/\delta))}{N\mu_N^{\alpha/(1+\alpha- (2+\alpha)\tau)}}  + \frac{1}{N \lambda_N^{\frac{\alpha}{1 - \tau(2+ \alpha)}}}.
	\end{align*}
\end{lemma}

\textbf{Proof of Lemma \ref{lemma: second term value}}
We omit $\beta$ in $Q^{\pi,\beta}$, $U^{\pi,\beta}$ and their relative quantities for the ease of presentation. Fix some $\delta > 0$.  Define a functional $f: (S, A) \mapsto g^2(S, A)$ for notational convenience, we decompose the error by 
\begin{align*}
\norm{\widehat g_N^\pi(\widehat \eta_N^\pi, \widehat Q_N^\pi)}_{}^2 + \norm{\widehat g_N^\pi(\widehat \eta_N^\pi, \widehat Q_N^\pi)}_{N}^2  = (P + \PN)f \left\{\widehat g_N^\pi (\widehat \eta_N^\pi, \widehat Q_N^\pi)\right\} 
%\\ & = 2 (\Pn f\left\{\widehat g_N^\pi(\widehat \eta_n^\pi, \widehat Q_N^\pi)right\} + \lambda_n J_1^2(\widehat Q_N^\pi)) + P f(\widehat g_N^\pi (\widehat \eta_n^\pi, \widehat Q_N^\pi) - 2 (\Pn f\left\{\widehat g_N^\pi(\widehat \eta_n^\pi, \widehat Q_N^\pi) \right\} + \lambda_n J_1^2(\widehat Q_N^\pi))\\
=I_{1} + I_{2},
\end{align*}
where
\begin{align*}
&I_{1} = 3 \left[\PN f\left\{\widehat g_N^\pi(\widehat \eta_N^\pi, \widehat Q_N^\pi) \right\} + (2/3)\lambda_n J_1^2(\widehat Q_N^\pi)\right]\\
&  I_{2} = (\PN +P) f\left\{\widehat g_N^\pi (\widehat \eta_N^\pi, \widehat Q_N^\pi)\right\}  -I_{1}.
\end{align*}
Denote $\widehat \eta_N^\pi$ as an estimation of $M(\beta, \pi)$ using the Bellman equation of the relative value function. We assume the average reward estimates $\widehat \eta_N^\pi \in \bar B$,  otherwise we can first show the consistency and use the high probability bound to focus on the truncation of this estimator. For the first term $I_{1}$, assumptions in Lemma \ref{lemma: uniform l2 bound value}, the optimizing property \eqref{lower level2 app} and the in-sample error bound in Lemma \ref{lemma: uniform l2 bound value} imply that for some fixed $\tau_1 \leq \frac{1}{3}$ and for sufficiently large $N$, the following holds with probability at least $1-\delta$, 
\begin{align*}
I_{1}
%& = \Pn f\left\{\widehat g_N^\pi(\widehat \eta_n^\pi, \widehat Q_N^\pi) \right\} + \lambda_n J_1^2(\widehat Q_N^\pi)\\
& \leq 3\PN f(\widehat g_N^\pi(\eta^\pi, \tilde Q^{\pi})) + 3\lambda_N J_1^2(\tilde Q^\pi)\\
& = 3\PN \left\{\widehat g_N^\pi(S, A; \eta^\pi, \tilde Q^{\pi})^2\right\}+ 3\lambda_N J_1^2(\tilde Q^\pi)\\
& = 3\PN \left[ \left\{\widehat g_N^\pi(S, A; \eta^\pi, \tilde Q^\pi) - \g(S, A; \eta^\pi, \tilde Q^\pi)\right\}^2\right]+ 3\lambda_N J_1^2(\tilde Q^\pi)\\
%& \leqconst \left(\mu_n J_1^2(\tilde Q^\pi)  + \frac{1}{n \mu_n^{\alpha}} + \frac{1}{n} +  \frac{\log(3/\delta)}{n}\right) + \lambda_n J_1^2(\tilde Q^\pi)\\
& = 3 \norm{\widehat g_N^\pi(\eta^\pi, \tilde Q^\pi) - \g(\eta^\pi, \tilde Q^\pi)}_{N}^2 + 3\lambda_N J_1^2(\tilde Q^\pi)\\
& \leqconst  \mu_N + \mu_N J_2^2 \left\{\g(\eta^\pi, \tilde{Q}^\pi)\right\}+ \left(\mu_N+\lambda_N\right) J_1^2(\tilde{Q}^\pi) +
\frac{(VC(\Pi)+1)\left[\log\left(\max(1/\delta, N)\right)\right]^{\frac{1}{\tau_1}}}{N} + 	\frac{1+VC(\Pi)}{N \mu_N^{\frac{\alpha}{1 - \tau_1(2+ \alpha)}}} ,
\end{align*}
where in the second equality we use $\g(\eta^\pi, \tilde Q^\pi) = 0$ from Assumption \ref{assumption: value function}~(b).  

The second term $I_{2}$ can be written as
\begin{align*}
I_{2} & = (\PN + P) f\left\{\widehat g_N^\pi (\widehat \eta_N^\pi, \widehat Q_N^\pi)\right\}  - 3 (\PN f\left\{\widehat g_N^\pi(\widehat \eta_N^\pi, \widehat Q_N^\pi) \right\} +  (2/3)\lambda_N J_1^2(\widehat Q_N^\pi))\\
& = 2(P-\PN) f\left\{\widehat g_N^\pi (\widehat \eta_N^\pi, \widehat Q_N^\pi)\right\} - P f\left\{\widehat g_N^\pi (\widehat \eta_N^\pi, \widehat Q_N^\pi)\right\} - 2\lambda_N J_1^2(\widehat Q_N^\pi).
\end{align*}
Define the constant 
\begin{align*}
\zeta^2(N, \mu_N, \delta, \tau_1) &= 1+ \frac{(VC(\Pi)+1)\left[\log\left(\max(1/\delta, N)\right)\right]^{\frac{1}{\tau_1}}}{N\mu_N} + \frac{1+VC(\Pi)}{N \mu_N^{\frac{1 - \tau_1(2+ \alpha) + \alpha}{1 - \tau_1(2+ \alpha)}}}.
\end{align*}

Using the probability bound on the complexity (i.e., $J_2(\widehat g_N^\pi(\eta, Q))$) developed in Lemma \ref{lemma: uniform l2 bound value} and Assumption \ref{assumption: value function}~(d), we can show that with probability at least $1 - \delta$, 
\begin{align*}
J_2\left\{\widehat g_N^\pi(\widehat \eta_N^\pi, \widehat Q_N^\pi)\right\}  &  \leqconst\left[ J_1(\widehat Q_N^\pi) +  J_2\left\{\g (\widehat \eta_N^\pi, \widehat Q_N^\pi)\right\} + \zeta(N, \mu_N, \delta, \tau_1)\right] \\
& \leqconst\left\{ J_1(\widehat Q_N^\pi) + J_1(\widehat Q_N^\pi) + \zeta(N, \mu_N, \delta, \tau_1)\right\} \\
& \leqconst\left\{J_1(\widehat Q_N^\pi)  + \zeta(N, \mu_N, \delta, \tau_1)\right\} \\
& = c_1 \left\{J_1(\widehat Q_N^\pi) +  \zeta(N, \mu_N, \delta, \tau_1)\right\},
\end{align*} 
for some constant $c_1$. 
For simplicity, we denote this event by $E =\Bigdkh{J_2\left\{\widehat g_N^\pi(\widehat \eta_N^\pi, \widehat Q_N^\pi)\right\} \leq c_1 \left\{J_1(\widehat Q_N^\pi) +  \zeta(N, \mu_N, \delta, \tau_1)\right\}}$. 
%%\peng{$$J_2(\widehat g_N^\pi(\eta, Q)) \leq c_0(J_1(Q) + J_2(\g(\eta, Q)) + \beta(n, \mu_n, \delta))$$ for some constant $c_0$ when $\mu_n = O_P(n^{-1/(1+\alpha)})$. }
%%Combining with Assumption \ref*{assumption: bounding J(EQ)}, for some constant $c_1$ (depending $c_0$ and $C_1, C_2$ in Assumption \ref*{assumption: bounding J(EQ)}), we have $\Pr(E) > 1-\delta/4$, where the event $E$ is given by
%\begin{align*}
%
%\end{align*}
Then we can have $\Pr(I_{2} > t) \leq \Pr \left\{(I_{2} > t) \jiao E\right\} + \delta$ and all we need to  bound is the first term using peeling device on $2\lambda_N J_1^2(\widehat Q_N^\pi)$ in $I_{2}$.  More specifically,
\begin{align*}
& \Pr \left\{(I_{2} > t) \jiao E\right\} 
= \sum_{l=0}^\infty \Pr \left[\{ I_{2} > t, ~2\lambda_N J_1^2(\widehat Q_N^\pi) \in [2^{l}t \indicator{t\neq 0}, 2^{l+1} t)\} \jiao E\right]\\
& \leq \sum_{l=0}^\infty \Pr\Big[~ 2(P-\PN)f \left\{\widehat g_N^\pi (\widehat \eta_N^\pi, \widehat Q_N^\pi)\right\} > P f\left\{\widehat g_N^\pi (\widehat \eta_N^\pi, \widehat Q_N^\pi)\right\} + 2\lambda_N J_1^2(\widehat Q_N^\pi) + t, \\
& \hspace{12ex} 2\lambda_N J_1^2(\widehat Q_N^\pi) \in [2^{l}t \indicator{t\neq 0}, 2^{l+1} t), J_2\left\{\widehat g_N^\pi(\widehat \eta_N^\pi, \widehat Q_N^\pi)\right\} \leq c_1 \left\{J_1(\widehat Q_N^\pi) + \zeta(N, \mu_N, \delta, \tau_1)\right\}\Big]\\
& \leq  \sum_{l=0}^\infty \Pr\Big[~2(P-\PN) f\left\{\widehat g_N^\pi (\widehat \eta_N^\pi, \widehat Q_N^\pi)\right\} > P f\left\{\widehat g_N^\pi (\widehat \eta_N^\pi, \widehat Q_N^\pi)\right\} + 2^{l}t \indicator{t\neq 0} + t, \\
& \hspace{12ex} 2\lambda_N J_1^2(\widehat Q_N^\pi) \leq 2^{l+1} t, J_2\left\{\widehat g_N^\pi(\widehat \eta_N^\pi, \widehat Q_N^\pi)\right\} \leq c_1 \left\{ \sqrt{(2^l t)/\lambda_N} + \zeta(N, \mu_N, \delta, \tau_1)\right\}\Big]\\
& \leq  \sum_{l=0}^\infty \Pr\Big[~2(P-\PN) f\left\{\widehat g_N^\pi (\widehat \eta_N^\pi, \widehat Q_N^\pi)\right\} > P f\left\{\widehat g_N^\pi (\widehat \eta_N^\pi, \widehat Q_N^\pi)\right\} + 2^{l}t, \\
& \hspace{12ex} J_2\left\{\widehat g_N^\pi(\widehat \eta_N^\pi, \widehat Q_N^\pi)\right\} \leq c_1  \left\{\sqrt{ (2^l t)/\lambda_N} + \zeta(N, \mu_N, \delta, \tau_1)\right\}\Big]\\
& \leq \sum_{l=0}^\infty \Pr\left[ \sup_{h \in \F_l}  \frac{(P-\PN) \left\{h(S, A)\right\}}{P \left\{h(S, A)\right\} + 2^l t} > \frac{1}{2} \right],
\end{align*}
where $\F_l = \left\{  f(g) : J_2(g) \leq c_1 \left\{ \sqrt{(2^l t)/\lambda_N} + \zeta(N, \mu_N, \delta, \tau_1) \right\}, g \in \G \right\}$. It is easy to see that  $\abs{f(g)(S, A)} \leq G_{\max}^2 \triangleq K_1$.

Similar to Lemma \ref{lemma: uniform l2 bound value}, we bound each term of the above probabilities by using the independent block technique. 
For each $l \geq 0$, we will use an  independent block sequence $(x_{N, l}, v_{N, l})$ with the residual $R_{ l}$. By controlling the size of these blocks, we can  optimize the bound. We let
$$
x_{N, l} = \lfloor x'_{N, l} \rfloor \epc \mbox{and } \epc  v_{N, l} = \lfloor \frac{N}{2x_{N, l}} \rfloor,
$$
where $x'_{N, l} = (Nt)^\tau (2^l)^p$ and $v'_{N, l} = \frac{N}{2x'_{N, l}}$ with some positive constants $\tau$ and $p$. Let $\tau \leq p \leq \frac{1}{2 + \alpha} \leq \frac{1}{2}$ and $N$ satisfies the following constraint:
\begin{align}\label{sample size constraint for Q 1}
N \geq c_1 \triangleq 4 \times 8^2 \times K_1 \geq 4^{\frac{p}{1-p}}8^{\frac{1}{1-p}}.
\end{align} 
By the definition of $x'_{N,l}$ and assuming $t \geq \frac{1}{N}$, $x_{N,l}\geq 1$. Then we consider two cases. The first case is  any $l$ such that $x'_{N,l} \geq \frac{N}{8}$. In such case, based on the assumption over $\tau$ and $p$, we can show that $x'_{N, l} \leq (Nt2^l)^p$, which further implies that $(Nt2^l) \geq 4NK_1$ by the sample constraint and $p \leq \frac{1}{2+\alpha}$. Then we can show that for this case,
$$
\frac{(P-\PN) \left\{h(S, A)\right\}}{P \left\{h(S, A)\right\} + 2^l t} \leq  \frac{2K_1}{2^lt} \leq  \frac{1}{2},
$$
for sufficiently large $N$.
Thus such terms does not contribute to the probability bound.

The second case we consider is any $l$ such that $x'_{N, l} < \frac{N}{8}$. We again apply the relative deviation concentration inequality for the exponential $\boldsymbol{\beta}$-mixing stationary process given in Theorem 4 of \cite{farahmand2012regularized}, which combined results in \cite{yu1994rates} and Theorem 19.3 in \cite{gyorfi2006distribution}. It then suffices to verify conditions (C1)-(C5) in Theorem 4 of \cite{farahmand2012regularized} with $\F = \F_l$, $\epsilon = 1/2$ and $\eta = 2^{l} t$ to get an exponential inequality for each term in the summation. The conditions (C1) has been verified. For (C2), we have $P f^2(g) \leq G_{\max}^2 P f(g)$ and thus (A2) holds by choosing $K_2 = G_{\max}^2$

%$K_1 ＝ 6G_{\max}(2R_{\max}+ 2Q_{\max} + 3G_{\max})$ and  $K_2$
%6G_{\max}(2R_{\max}+ 2Q_{\max} + 3G_{\max})$. 
%There exists some $K_1, K_2$ (depending on $\pi_{\min}$, $R_{\max}, Q_{\max}$) such that $\norm{f}_\infty \leq K_1$ and $\EE\left[f(\D)^2\right] \leq K_2 \EE\left[f(\D)\right]$. 

To verify the condition  (C3), without loss of generality, we assume $K_1 \geq 1$. Otherwise, let $K_1 = \max(1, K_1)$. Then we know that $2K_1x_{N, l} \geq \sqrt{2K_1x_{N, l}}$ since $x_{N, l} \geq 1$. We need to have $\sqrt{N} \epsilon\sqrt{1-\epsilon} \sqrt{\eta} \geq 1152K_1x_{N, l}$, or suffice to have $\sqrt{N} \epsilon\sqrt{1-\epsilon} \sqrt{\eta} \geq 1152K_1x'_{N, l}$. Recall that $\epsilon = 1/2$ and $\eta = 2^{l} t$. So it is enough to show that
$$
\sqrt{N}\frac{\sqrt{2}}{4}  \sqrt{2^lt} \geq 1152K_1(Nt2^l)^p.
$$
We can check that if $t \geq \frac{2304\sqrt{2}K_1}{N}$, the above inequality holds for every $l \geq$ since $p \leq \frac{1}{2+\alpha}$.

Next we verify (C4) that $\frac{|R_l|}{N} \leq \frac{\epsilon\eta}{6K_1}$.  Recall that $|R_l| \leq 2 x_{N, l} \leq 2x'_{N, l} = (Nt)^\tau(2^l)^p$. So if $t \geq \frac{c_2}{n}$ for some positive constant $c_2$, we can have
$$
\frac{\epsilon\eta}{6K_1} = \frac{2^lt}{12K_1} \geq \frac{2(Nt)^\tau (2^l)^p}{N} = \frac{2x'_{N,l}}{N}\geq \frac{|R_l|}{N}.
$$
In addition, $|R_l| \leq 2 x'_{N, l} < \frac{N}{2}$.

We now verify the final condition (C5).  First, we obtain an upper bound $\mathcal{N}_{}(u, \F_l; \norm{\cdot}_\infty)$  for all possible realization of $(S, A)$.  For any $g_1, g_2 \in \G$, 
\[
\PP_N\left[f(g_1)(S, A) - f(g_2)(S, A)\right]^2 \leq 4 G_{\max}^2 \norm{g_1 - g_2}_{N}^2.
\]
Thus applying Assumption \ref{assumption: function classes}~(d) implies that for some constant $c_3$, the metric entropy for each $l$ is bounded  by
\begin{align*}
& \log {\mathcal{N}}(u, \F_l, \norm{\cdot}_\infty) \\
& \leq \log {\mathcal{N}}\left(\frac{u}{2G_{\max}}, \{g: J_2(g) \leq c_1 ( \sqrt{(2^l t)/\lambda_N }+ \zeta(N, \mu_N, \delta, \tau_1)), g \in \G  \}, \norm{\cdot}_\infty \right)\\
& \leqconst \left[\frac{c_1 \left\{ \sqrt{(2^l t)/\lambda_N} + \zeta(N, \mu_N, \delta, \tau_1)\right\}}{u/(2G_{\max})}\right]^{2\alpha} \leq c_3 \left\{ \left(\frac{2^l t}{\lambda_N}\right)^\alpha + \zeta(N, \mu_N, \delta, \tau_1)^{2\alpha}\right\} u^{-2\alpha},
\end{align*} 
for some positive constant $c_3$.

Now we see the condition (C5) is satisfied if the following inequality holds for all $x \geq (2^ltx_{N, l})/8$ such that
\begin{align*}
\frac{\sqrt{v_{N, l}} (1/2)^2 x}{96x_{N, l}\sqrt{2} \max(K_1, 2K_2)} 
&\geq \int_{0}^{\sqrt{x}} \sqrt{c_3} \left\{ \left(\frac{2^l t}{\lambda_N}\right)^\alpha + \zeta(N, \mu_N, \delta, \tau_1)^{2\alpha}\right\}^{1/2}  \left(\frac{u}{2x_{N, l}}\right)^{-\alpha} du \\
&= x_{N, l}^{\alpha} x^{\frac{1-\alpha}{2}} \sqrt{2^\alpha c_3} \left( \left(\frac{2^l t}{\lambda_N}\right)^\alpha + \zeta(N, \mu_N, \delta, \tau_1)^{2\alpha}\right)^{1/2}.
\end{align*}
It is sufficient to let the following inequality hold:
\begin{align*}
\frac{\sqrt{v_{N, l}}  }{384x_{N, l}\sqrt{2} \max(K_1, 2K_2)}  x^{\frac{1+\alpha}{2}} \geq  \sqrt{c_3'}x_{N, l}^{\alpha}  \left\{ \left(\frac{2^l t}{\lambda_n}\right)^\alpha + \zeta(N, \mu_N, \delta, \tau_1)^{2\alpha} \right\}^{1/2},
\end{align*}
for some constant $c_3'$.
Using the inequality that $(a+b)^{1/2} \leq \sqrt{a} + \sqrt{b}$ and the fact that LHS is increasing function of $x$, it's enough to ensure that the following two inequalities hold:
\begin{align*}
& \frac{\sqrt{v_{N, l}}  }{384x_{N, l}\sqrt{2} \max(K_1, 2K_2)}  (x_{N, l}2^l t/8)^{\frac{1+\alpha}{2}} \geq  \sqrt{c_3'} x_{N, l}^{\alpha} \left(\frac{2^l t}{\lambda_N}\right)^{\alpha/2}\\
& \frac{\sqrt{v_{N, l}}  }{384x_{N, l}\sqrt{2} \max(K_1, 2K_2)}  (x_{N, l}2^l t/8)^{\frac{1+\alpha}{2}} \geq \sqrt{c_3'}x_{N, l}^{\alpha}  \zeta(N, \mu_N, \delta, \tau_1)^\alpha.
\end{align*}
By the definition of $v_{N, l}$ and $x_{N, l}$, after some algebra, we can see that the first inequality holds if
$$
t \geq c_5 \frac{1}{N \lambda_N^{\frac{\alpha}{1 - \tau(2+ \alpha)}}}.
$$
The second inequality holds if
$t$ satisfies
\begin{align*}
t & \geq  c_6 N^{-\frac{1-(2+\alpha)\tau}{1+\alpha- (2+\alpha)\tau}}\zeta(N, \mu_N, \delta, \tau_1)^{\frac{2\alpha}{1+\alpha - (2+\alpha)\tau}}.
\end{align*}
Choosing $\tau = \tau_1 \leq 1/3$, we can obtain that
\begin{align*}
&N^{-\frac{1-(2+\alpha)\tau}{1+\alpha-\tau(2+\alpha)}}\zeta(N, \mu_N, \delta, \tau_1)^{\frac{2\alpha}{1+\alpha - (2+\alpha)\tau}}\\ 
& = N^{-\frac{1-(2+\alpha)\tau}{1+\alpha-\tau(2+\alpha)}} \left[1+ \frac{(VC(\Pi)+1)\left[\log\left(\max(1/\delta, N)\right)\right]^{\frac{1}{\tau_1}}}{N\mu_N} + \frac{1+VC(\Pi)}{N \mu_N^{\frac{1 - \tau_1(2+ \alpha) + \alpha}{1 - \tau_1(2+ \alpha)}}} \right]^{\frac{\alpha}{1+\alpha - (2 +\alpha)\tau}} \\
& \leqconst N^{-\frac{1-(2+\alpha)\tau}{1+\alpha-\tau(2+\alpha)}} + \frac{1+VC(\Pi)}{N \mu_N^{\alpha/(1-\tau(2+\alpha))}} + \frac{(VC(\Pi)+1)\log^{\frac{\alpha/\tau}{1+\alpha - \tau(2+\alpha)}}(\max(N,1/\delta))}{N\mu_N^{\alpha/(1+\alpha- (2+\alpha)\tau)}}. 
%& \leqconst N^{-\frac{1-(2+\alpha)\tau}{1+\alpha-\tau(2+\alpha)}} + \frac{\log^{\frac{\alpha/\tau}{1+\alpha - \tau(2+\alpha)} }(\max(N,1/\delta))+ VC(\Pi)}{N \mu_N^{\alpha/(1-\tau(2+\alpha))}}. \mbox{\QZL{Don't proceed to this step}}
\end{align*}

Putting all together, all conditions (C1) to (C5) would be satisfied for all $l \geq 0$ when
\begin{align*}
t \geq & \frac{c_2}{N}  +c'_5 \left\{N^{-\frac{1-(2+\alpha)\tau}{1+\alpha-(2+\alpha)\tau}} + \frac{1+VC(\Pi)}{N \mu_N^{\alpha/(1-\tau(2+\alpha))}} + \frac{(VC(\Pi)+1)\log^{\frac{\alpha/\tau}{1+\alpha - \tau(2+\alpha)}}(\max(N,1/\delta))}{N\mu_N^{\alpha/(1+\alpha- (2+\alpha)\tau)}}\right\} +c_5 \frac{1}{N \lambda_N^{\frac{\alpha}{1 - \tau(2+ \alpha)}}},
\end{align*}
for some constant $c_5'$.

Applying Theorem 4 in \cite{farahmand2012regularized} with $\F = \F_l$, $\epsilon = 1/2$ and $\eta = 2^{l} t$, for sufficiently large $N$, we can obtain that 
\begin{align*}
& \Pr\left\{\exists (\beta, \pi, \eta, Q) \in B \times \Pi \times B \times \Q,  I_2(\eta, Q) > t\right\} \\
& \leq \sum_{l=0}^\infty \Pr \left[\sup_{h \in \F_l}  \frac{(P-\PN) \left\{h(S, A)\right\}}{P \left\{h(S, A)\right\} + 2^l t} > \frac{1}{2} \right] \\
& \leq \sum_{l=0}^\infty 120 \exp\left\{-c_6\frac{v_{N,l }^{'2}t2^l}{N}\right\} + 2\beta_{x_{N,l}}v_{N, l} \\
& \leq \sum_{l=0}^\infty 120 \exp\left\{-c_6\frac{v_{N,l }^{'2}t2^l}{N}\right\} + 2\beta_0\exp\left(-\beta_1 x_{N, l} + \log v_{N, l}\right) ,
\end{align*}
where the last inequality is based on Assumption \ref{ass: stationarity}.  When $t \geq \frac{\left(4/\beta_1\log(N)\right)^{1/\tau}}{N}$, we have $\log v_{N, l} \leq \frac{1}{2}\beta_1 x_{N, l}$. This will further imply that $2\beta_{x_{N,l}}v_{N, l} \leq 2\beta_0\exp\left(-\beta_1 x_{N, l}/2 \right)$. Then we will have 
\begin{align*}
& \Pr\left\{\exists (\beta, \pi, \eta, Q) \in B \times \Pi \times B \times \Q,  I_2(\eta, Q) > t\right\} \\
& \leq \sum_{l=0}^\infty 120 \exp\left\{-c_6\frac{v_{N,l }^2t2^l}{N}\right\} + 2\beta_0\exp\left(-\beta_1 x_{N, l} + \log v_{N, l}\right)\\
& \leqconst \sum_{l=0}^\infty 120\exp\left(-c_7(Nt)^{1-2\tau}(2l)^{1-2p} \right)+ 2\beta_0\exp\left(-\beta_1 (Nt)^{\tau}(2^l)^p\right)\\
& \leq c_8 \exp\left(-c_9(Nt)^{1-2\tau}\right) + c_{10}\exp\left(-c_{11} (Nt)^{\tau}\right).
\end{align*}
As long as $t$ satisfies all the above constraints, then 
$$
I_2(\eta, Q) \leq \frac{1}{N}\left\{\left(\frac{\log(\frac{2c_8}{\delta})}{c_9}\right)^{\frac{1}{1-2\tau}} + \left(\frac{\log(\frac{2c_{10}}{\delta})}{c_{11}}\right)^{\frac{1}{\tau}} \right\},
$$
with probability at least $1-\delta$.
Collecting all the conditions on $t$ and combining with the bound of $I_1(\eta, Q)$, we have shown that with probability at least $1-2\delta$, the following holds for all $(\beta, \pi, \eta, Q) \in B \times \Pi \times B \times \Q$:
\begin{align*}
& \norm{\widehat g_N^\pi(\widehat \eta_N^\pi, \widehat Q_N^\pi)}_{}^2 + \norm{\widehat g_N^\pi(\widehat \eta_N^\pi, \widehat Q_N^\pi)}_{N}^2  \\
& \leqconst  \mu_N+ \mu_N J_2^2 \left\{\g(\eta^\pi, \tilde{Q}^\pi)\right\}+ (\mu_N + \lambda_N) J_1^2(\tilde{Q}^\pi) +
\frac{(VC(\Pi)+1)\left[\log\left(\max(1/\delta, N)\right)\right]^{\frac{1}{\tau}}}{N} + 	\frac{1+ VC(\Pi)}{N \mu_N^{\frac{\alpha}{1 - \tau(2+ \alpha)}}} \\
& + N^{-\frac{1-(2+\alpha)\tau}{1+\alpha-\tau(2+\alpha)}} +  \frac{(VC(\Pi)+1)\log^{\frac{\alpha/\tau}{1+\alpha - \tau(2+\alpha)}}(\max(N,1/\delta))}{N\mu_N^{\alpha/(1+\alpha- (2+\alpha)\tau)}} + \frac{1}{N \lambda_N^{\frac{\alpha}{1 - \tau(2+ \alpha)}}}, 
\end{align*}
which concludes our proof.

\subsection{Finite Sample Error Bounds of Ratio Functions}
We begin with the following lemma.
\begin{lemma} \label{thm: pre-ratio}
	Suppose assumptions  \ref{ass: stationarity}, \ref{ass: policy class}, \ref{assumption: function classes}, \ref{assumption: ratio function} hold.  Let $\widehat \omega^\pi_N$ be the estimated ratio function with tuning parameter $\mu_{2N} \simeq \lambda_{2N} \simeq (1+VC(\Pi))(\log N)^{\frac{2+\alpha}{1+\alpha}}N^{-\frac{1}{1+\alpha}}$ defined in \eqref{ratio estimator}. 
	For any $m \geq 1$, there exists some constant such that with sufficiently large $N$, the following  holds with probability at least {$1-\frac{3+2m}{N}$}
	\begin{align*}
	& \norm{\h(\Hn) - \h(\Hpi)}_2^2 \leqconst \left(1+VC(\Pi)\right)\left[\log\left(\max(1/\delta, N)\right)\right]^{\frac{2+\alpha}{1+\alpha}} N^{-\frac{r_m+\frac{1}{1+\alpha}}{2}}
	\\
	& J_1^2(\Hn) \leqconst N^{\frac{ \frac{1}{1+\alpha}}{2}-\frac{r_m}{2}}, %=  C_{new}  (\delta) \iota^{\omega_{k+1}} p n^{1/(1+\alpha) - \beta_{new}}
	\end{align*}
	where $r_{m} = \frac{r_{m-1}+ 1/(1+\alpha)}{2} = \frac{1}{1+\alpha} -  \frac{(1-\alpha) 2^{-(m-1)}}{1+\alpha}$.
\end{lemma}
\textbf{Proof of Lemma \ref{thm: pre-ratio}}
%	Link the lemmas below to show the convergence rate of the ratio estimator. 
We start with 
\begin{align}
\norm{\h(\Hn) - \h(\Hpi)}^2 
\leq 2 \norm{\h(\Hn) - \hn(\Hn)}^2 + 2\norm{\hn(\Hn) - \h(\Hpi)}^2  \label{target}
\end{align}
The first term can be bounded by Lemma \ref{lemma: uniform l2 bound}. For sufficiently large $N$ and $\tau \leq \frac{1}{3}$, with probability at least $1-\delta$,  for all $\pi \in \Pi$,   we can have
%		\[
%		\norm{\g(\Hn) - \gn(\Hn)}_2^2 \leq C_1(\delta)pn^{-1/(1+\alpha)} + \inblue{K_1} n^{-1/(1+\alpha)} J_1^2(\Hn)
%		\]
\begin{align*}
& \norm{\h(\Hn) - \hn(\Hn)}_2^2 \\
& \leqconst  \mu_N (1+J_1^2(\Hn) + J_2^2(\h (\Hn))) + \frac{(VC(\Pi)+1)\left[\log(\max(N, 1/\delta))\right]^{\frac{1}{\tau}}}{N} + \frac{1 + VC(\Pi)}{N\mu_N^{\frac{\alpha}{1- \tau(2+\alpha)}}}\\
& \leqconst  \mu_N \left(1+J_1^2(\Hn)\right) + \frac{(VC(\Pi)+1)\left[\log(\max(N, 1/\delta))\right]^{\frac{1}{\tau}}}{N} + \frac{1 + VC(\Pi)}{N\mu_N^{\frac{\alpha}{1- \tau(2+\alpha)}}}.\\
\end{align*}

Now we discuss the second term. We apply Lemma \ref{lemma: upper level initial bound} with the same $\tau$ as above. Then for sufficiently large $N$, with the probability at least $1-2\delta$, for all $\pi \in \Pi$
\begin{align*}
& \norm{\hn(\Hn) - \h(\Hpi)}^2 + \lambda_N J_1^2(\Hn) 
\\
&\leqconst  \zeta_2(\delta, N, VC(\Pi), \mu_N, \lambda_N, \tau) + \Rem(\pi) +  \mu_N J_1(\Hn).
\end{align*}
Recall that $\Rem(\pi) =4\bigabs{\PN \h(S; A; \Hpi) [\Delta^\pi(S, A, S'; \Hn) - \Delta^\pi(S, A, S; \Hpi)] }$.
Here we define
\begin{align*}
&\zeta_2(\delta, N, VC(\Pi), \mu_N, \lambda_N, \tau) = (\mu_N+\lambda_N) (1+\sup_{\pi \in \Pi}J_1^2(\tilde H^\pi)) +\frac{(VC(\Pi)+1)\left[\log(\max(N, 1/\delta))\right]^\frac{1}{\tau}}{N} \\
& + \frac{1+VC(\Pi)}{N\mu_N^{\frac{\alpha}{1-\tau(2+\alpha)}}} +\frac{1}{N \lambda_N^{\frac{\alpha}{1 - \tau(2+ \alpha)}}} + \frac{\sqrt{\mu_N(VC(\Pi)+1)}\left[\log(\max(N, 1/\delta))\right]^{\frac{1}{2\tau}}}{\sqrt{N}}
\\
& + N^{-\frac{1-(2+\alpha)\tau}{1+\alpha-(2+\alpha)\tau}} + \frac{(VC(\Pi)+1)\log^{\frac{\alpha/\tau}{1+\alpha - \tau(2+\alpha)}}(\max(N,1/\delta))}{N\mu_N^{\alpha/(1+\alpha- (2+\alpha)\tau)}}  + \frac{\sqrt{1 + VC(\Pi)}}{\sqrt{N}\mu_N^{\frac{\alpha + \tau(2+\alpha)- 1}{2(1- \tau(2+\alpha))}}}.
\end{align*}
Letting
$$
\lambda_N \cong \mu_N  = \left(1+VC(\Pi)\right)(1+ \sup_{\pi \in \Pi}J_1^2(\tilde{H}^\pi))^{-\frac{1}{1+\alpha}}\left[\log\left(\max(1/\delta, N)\right)\right]^{\frac{2+\alpha}{1+\alpha}} N^{-\frac{1}{1+\alpha}},
$$
and 
$$
\tau = \frac{(1+\alpha)\log(B)}{\alpha \log(A) + (2+ \alpha)\log(B)},
$$
where
\begin{align*}
A & = N\left(1+ \sup_{\pi \in \Pi}J_1^2(\tilde{H}^\pi)\right)\\
B & = \log(\max(N, 1/\delta)).
\end{align*}
we can show that the first seven terms in $\zeta_2(\delta, N, VC(\Pi), \mu_N, \lambda_N, \tau)$ is proportionally less than or equal to
\begin{align*}
&  \left(1+VC(\Pi)\right)(1+ \sup_{\pi \in \Pi}J_1^2(\tilde{H}^\pi))^{\frac{\alpha}{1+\alpha}}\left[\log\left(\max(1/\delta, N)\right)\right]^{\frac{2+\alpha}{1+\alpha}} N^{-\frac{1}{1+\alpha}},
\end{align*}
which is similar to the derivation in the proof of Theorem \ref{thm: value}.
Now we discuss the last term of $\zeta_2(\delta, N, VC(\Pi), \mu_N, \lambda_N, \tau)$. As we know that
$$
\mu_N >  \bar{\mu}_N =
\left(1+VC(\Pi)\right)^{\frac{1}{1+\alpha }- 2\alpha \log(B) }(1+ \sup_{\pi \in \Pi}J_1^2(\tilde{H}^\pi))^{-\frac{1}{1+\alpha}}\left[\log\left(\max(1/\delta, N)\right)\right]^{\frac{2+\alpha}{1+\alpha}} N^{-\frac{1}{1+\alpha}}$$.
In addition, by the definition of $\bar \mu_N$, we can show that
$$
\bar{\mu}_N = \left[\frac{1+VC(\Pi)}{ A}\right]^{\frac{1-\tau(2+\alpha)}{1+\alpha -\tau(2+\alpha) }} \triangleq D^{\frac{1-\tau(2+\alpha)}{1+\alpha -\tau(2+\alpha) }}.
$$
Then we can show that
\begin{align*}
& \frac{\sqrt{1 + VC(\Pi)}}{\sqrt{N}\mu_N^{\frac{\alpha + \tau(2+\alpha)- 1}{2(1- \tau(2+\alpha))}}} \leq \sqrt{\frac{1 + VC(\Pi)}{N\bar \mu_N^{-1+ \frac{\alpha}{1-\tau(2+\alpha)}}}} = \sqrt{\frac{\bar \mu^2_N \left(1 + VC(\Pi)\right)}{N\times D}}\\
& = \left(1+VC(\Pi)\right)^{\frac{1}{1+\alpha } - 2\alpha\log(B)}(1+ \sup_{\pi \in \Pi}J_1^2(\tilde{H}^\pi))^{-\frac{1-\alpha}{2(1+\alpha)}}\left[\log\left(\max(1/\delta, N)\right)\right]^{\frac{2+\alpha}{1+\alpha}} N^{-\frac{1}{1+\alpha}}.
\end{align*}
Combining together, we can demonstrate that
$$
\zeta_2(\delta, N, VC(\Pi), \mu_N, \lambda_N, \tau) \leqconst \left(1+VC(\Pi)\right)^{\frac{1}{1+\alpha }}(1+ \sup_{\pi \in \Pi}J_1^2(\tilde{Q}^\pi))^{\frac{\alpha}{1+\alpha}}\left[\log\left(\max(1/\delta, N)\right)\right]^{\frac{2+\alpha}{1+\alpha}} N^{-\frac{1}{1+\alpha}}.
$$
As a result,  we obtain that for $N$ sufficiently large and the chosen $\mu_N$, with probability at least $1-3\delta$, for all $\pi \in \Pi$, 
\begin{align}
\lambda_N J_1^2(\Hn) &\leqconst \left(1+VC(\Pi)\right)(1+ \sup_{\pi \in \Pi}J_1^2(\tilde{H}^\pi))^{\frac{\alpha}{1+\alpha}}\left[\log\left(\max(1/\delta, N)\right)\right]^{\frac{2+\alpha}{1+\alpha}} N^{-\frac{1}{1+\alpha}} \notag \\
& +  \Rem(\pi) + \mu_NJ_1(\Hn) \label{J iter}
\end{align}
%and 
\begin{align}
\norm{\h(\Hn) - \h(\Hpi)}^2 &\leqconst \left(1+VC(\Pi)\right)(1+ \sup_{\pi \in \Pi}J_1^2(\tilde{H}^\pi))^{\frac{\alpha}{1+\alpha}}\left[\log\left(\max(1/\delta, N)\right)\right]^{\frac{2+\alpha}{1+\alpha}} N^{-\frac{1}{1+\alpha}} \notag \\
& +  \Rem(\pi) + \mu_NJ_1(\Hn) \label{err iter}
\end{align}

\paragraph{Initial Rate}
We derive an initial rate by bounding $\Rem(\pi)$ uniformly over $\pi \in \Pi$.  Let 
$$
f(H_1, H_2) (S, A, S') = 4 \h(S; A; H_2) [\Delta^\pi(S, A, S'; H_1) - \Delta^\pi(S, A, S'; H_2)].$$
We thus have $\Rem(\pi) =  \abs{\PN f(\Hn, \Hpi)}$. Note that under Assumption \ref{b3-2}, $\h(\Hpi) = e^\pi$. The orthogonality property (\ref{lemma: orthogonal}) then implies that  $P f(H, \Hpi) = 0$ for any $H \in \F$. We can bound $\Rem(\pi)$ by 
%	\begin{align*}
%	\Rem(\pi) = \sup_{\pi \in \Pi} \abs{\Pn f(\Hn, \Hpi)}
%	\end{align*}
\begin{align*}
\Rem (\pi) = \abs{\PN f(\Hn, \Hpi)} = J_1(\Hn-\Hpi) \frac{\abs{\PN f(\Hn, \Hpi) }}{J_1(\Hn - \Hpi)}   \leq  J_1(\Hn-\Hpi) \sup_{f \in \F_0} \abs{\PN f}
\end{align*}
where the function class $\F_0$ is given by
\begin{align*}
%	\F_0 & =  \{(Y_{q} - Y_{\Hpi}) \g(\Hpi)/J_1(q-H^{\pi}): H \in \F, \pi \in \Pi \} \\
\F_0 & =  \Bigdkh{(S, A, S') \mapsto \big[\frac{(H^{\pi} - H)(S, A)}{J_1(H^{\pi} - H)} 
	- \sum_{a'} \pi(a'|S')   \frac{(H^{\pi} - H)(S', a')}{J_1(H^{\pi} - H)}\big] \h(S, A; \Hpi): H \in \F, \pi \in \Pi }\\
&= \Bigdkh{D \mapsto \Big[H(S, A) 
	- \sum_{a'} \pi(a'|S')   H(S', a')\Big] \h(S, A): H \in \F, J_1(H) = 1, \pi \in \Pi }\\
&\subset \Bigdkh{D \mapsto \Big[H(S, A) 
	- \sum_{a'} \pi(a'|S')   H(S', a')\Big] h(S, A; \Hpi): H \in \F, J_1(H) \leq 1, \pi \in \Pi, J(h) \leq \sup_{\pi \in \Pi} J_2(e^\pi) }\\
& \triangleq \F_1.
\end{align*}
%	\begin{align*}
%	\F = & \Bigdkh{ (S, A, S') \mapsto \big[q(S, A) - \sum_{a'} \pi(a'|S')    q (S', a')\big] g(S, A): \\
%		& \hspace{30ex}H \in \F, J_1(q) \leq 1, \pi \in \Pi, J_2(g) \leq \sup_{\pi \in \Pi} J_2(e^\pi)}
%	\end{align*}
{Applying  Lemma \ref{lemma: donsker class prob bound} with  $M=1$ and $\sigma = 2G_{\max}F_{\max}$ implies that the following holds with probability at least $1-\delta - \frac{1}{N}$:}
\begin{align*}
\sup_{f \in \F_1}  \abs{\PN f}  & \leqconst \frac{\sqrt{VC(\Pi) +1 }}{\sqrt{N}}\log(\max(N, 1/\delta))
\end{align*}
As a result, combing with (\ref{J iter}) and with probability at least $1-4\delta - \frac{1}{N}$,  the following holds for all $\pi$:
\begin{align*}
\lambda_N J_1^2(\Hn) \leqconst& 
\left(1+VC(\Pi)\right)(1+ \sup_{\pi \in \Pi}J_1^2(H^\pi))^{\frac{\alpha}{1+\alpha}}\left[\log\left(\max(1/\delta, N)\right)\right]^{\frac{2+\alpha}{1+\alpha}} N^{-\frac{1}{1+\alpha}}\\
+  & \Rem(\pi) + \mu_NJ_1(\Hn) \\  
\leqconst& 
\left(1+VC(\Pi)\right)(1+ \sup_{\pi \in \Pi}J_1^2(H^\pi))^{\frac{\alpha}{1+\alpha}}\left[\log\left(\max(1/\delta, N)\right)\right]^{\frac{2+\alpha}{1+\alpha}} N^{-\frac{1}{1+\alpha}}\\
+  & \frac{\sqrt{VC(\Pi) + 1}}{\sqrt{N}}\log(\max(N, 1/\delta))J_1(\Hn - \tilde H^\pi) + \mu_NJ_1(\Hn) \\ 
\leqconst& 
\left(1+VC(\Pi)\right)(1+ \sup_{\pi \in \Pi}J_1^2(H^\pi))^{\frac{\alpha}{1+\alpha}}\left[\log\left(\max(1/\delta, N)\right)\right]^{\frac{2+\alpha}{1+\alpha}} N^{-\frac{1}{1+\alpha}}\\
+  & \frac{\sqrt{VC(\Pi) + 1}}{\sqrt{N}}\log(\max(N, 1/\delta))J_1(\Hn) \\
+& \frac{\sqrt{VC(\Pi) + 1}}{\sqrt{N}}\log(\max(N, 1/\delta))\sup_{\pi \in \Pi}J_1(\tilde H^\pi) + \mu_NJ_1(\Hn) 
\end{align*}
Dividing $\lambda_N$ on both sides and recalling that $\lambda_N \cong \mu_N$ give that 
%		& = U_1(\delta, n, p, \mu_n, \lambda_n) + U_2(\delta, n, p, \mu_n )
\begin{align*}
J_1^2(\Hn) \leqconst& 
1 + \sup_{\pi \in \Pi} J^2_1(\tilde H^\pi)+  \left\{ \frac{\sqrt{VC(\Pi) + 1}}{\sqrt{N}}\log(\max(N, 1/\delta)) / \lambda_N+ 1\right\} J_1(\Hn) \\
+& \frac{\sqrt{VC(\Pi) + 1}}{\sqrt{N}}\log(\max(N, 1/\delta))\sup_{\pi \in \Pi}J_1(\tilde H^\pi) /\lambda_N
\end{align*}
Let $x = J_1(\Hn)$ and the above inequality becomes $x^2 \leq a + b x$ for some $a, b > 0$. When $a \leq bx$, we have $x^2 \leq 2bx$, or $x^2 \leq 4b^2$. When $a > bx$, we have $x^2 \leq a+bx \leq 2a$. Thus  $x^2 \leq \max(4b^2, 2a) \leq 2a + 4b^2$.  Now we have
\begin{align*}
J_1^2(\Hn) \leqconst& 1+\sup_{\pi \in \Pi} J^2_1(\tilde H^\pi)+ \sqrt{\frac{VC(\Pi) + 1}{N}}\log(\max(N, 1/\delta))\sup_{\pi \in \Pi}J_1(\tilde H^\pi) /\lambda_N\\
& + \left\{ \frac{\sqrt{VC(\Pi) + 1}}{\sqrt{N}}\log(\max(N, 1/\delta)) / \lambda_N+ 1\right\}^2 \\
\leqconst  &  N^{\frac{1-\alpha}{1+\alpha}} (1+\sup_{\pi \in \Pi} J^2_1(H^\pi))^{\frac{2+\alpha}{1+\alpha}},
\end{align*}
where without loss of generality, we assume $\sup_{\pi \in \Pi} J_1(\tilde H^\pi) \geq 1$.
Now summarizing the previous probability bounds, we can show that w.p. $1-4\delta - \frac{1}{N}$ for all $\pi \in \Pi$:
\begin{align*}
\norm{\h(\Hn) - \h(\Hpi)}^2 &\leqconst \left(1+VC(\Pi)\right)(1+ \sup_{\pi \in \Pi}J_1^2(H^\pi))^{\frac{\alpha}{1+\alpha}}\left[\log\left(\max(1/\delta, N)\right)\right]^{\frac{2+\alpha}{1+\alpha}} N^{-\frac{1}{1+\alpha}} \notag \\
& + \mu_N\big[N^{\frac{1-\alpha}{1+\alpha}} (1+\sup_{\pi \in \Pi} J^2_1(\tilde H^\pi))^{\frac{2 + \alpha}{1+\alpha}}\big]^{\frac{1}{2}}\\
& + \sqrt{\frac{VC(\Pi)+1}{N}}\log(\max(N, 1/\delta))\sup_{\pi \in \Pi}J_1(\tilde H^\pi) \\
& + \sqrt{\frac{VC(\Pi)+1}{N}}\log(\max(N, 1/\delta))\times \big[N^{\frac{1-\alpha}{1+\alpha}} (1+\sup_{\pi \in \Pi} J^2_1(\tilde H^\pi))^{\frac{2 + \alpha}{1+\alpha}}\big]^{\frac{1}{2}}\\
& \leqconst \left(1+VC(\Pi)\right)(1+ \sup_{\pi \in \Pi}J_1^2(H^\pi))\left[\log\left(\max(1/\delta, N)\right)\right]^{\frac{2+\alpha}{1+\alpha}} N^{-\frac{\alpha}{1+\alpha}}
\end{align*}
Let $r_1 = \frac{\alpha}{1+\alpha}$.  We have shown that for $N$ sufficiently large, with probability at least $1-4\delta - \frac{1}{N}$,  the inequalities (\ref{J iter}), (\ref{err iter}) and the followings hold:
\begin{align*}
& \norm{\h(\Hn) - \h(\Hpi)}^2 \leqconst \left(1+VC(\Pi)\right)(1+ \sup_{\pi \in \Pi}J_1^2(H^\pi))\left[\log\left(\max(1/\delta, N)\right)\right]^{\frac{2+\alpha}{1+\alpha}} N^{-\frac{\alpha}{1+\alpha}}\\
& J_1^2(\Hn) \leqconst  N^{\frac{1}{1+\alpha} - r_1} (1+\sup_{\pi \in \Pi} J^2_1(H^\pi))^{\frac{2+\alpha}{1+\alpha}}
\end{align*}

\paragraph{Rate Improvement}
Let $r = r_1$. Denote by $E_N$ the event that the inequalities (\ref{J iter}) and (\ref{err iter}) hold and 
\begin{align*}
& \norm{\h(\Hn) - \h(\Hpi)}^2 \leqconst \left(1+VC(\Pi)\right)(1+ \sup_{\pi \in \Pi}J_1^2(H^\pi))\left[\log\left(\max(1/\delta, N)\right)\right]^{\frac{2+\alpha}{1+\alpha}} N^{-r}\\
& J_1^2(\Hn) 	\leqconst  N^{\frac{1}{1+\alpha} - r} (1+\sup_{\pi \in \Pi} J^2_1(H^\pi))^{\frac{2+\alpha}{1+\alpha}}.
\end{align*}
We have shown that $\Pr(E_N) \geq 1-4\delta - 1/N$.  Below we improve the rate by refining the bound of the remainder term, $\Rem(\pi)$.   First, note that under the event $E_N$, 
\begin{align}
P f^2(\Hn, \Hpi) & \leqconst G_{\max}^2 \norm{\h(\Hn) - \h(\Hpi)}^2  \label{key ineq}\\
& \leqconst \left(1+VC(\Pi)\right)(1+ \sup_{\pi \in \Pi}J_1^2(H^\pi))\left[\log\left(\max(1/\delta, N)\right)\right]^{\frac{2+\alpha}{1+\alpha}} N^{-r}\triangleq (I),
\end{align}
where we use Assumption \ref{assumption: ratio function} (c).
In addition,
%	see Derivation of Inequality (B.9) in JASA)
similarly,
%		\[
%		\inblue{K} = \frac{1 +  (1 + (1/T) \norm{\frac{d_{T+1}}{d_D}}_\infty) (1/p_{\min})}{\kappa^2} > 1
%		\]
\begin{align*}
J_1(\Hn - \Hpi) & \leqconst N^{\frac{1}{2(1+\alpha)} - r/2} (1+\sup_{\pi \in \Pi} J^2_1(H^\pi))^{\frac{2+\alpha}{2(1+\alpha)}} +\sup_{\pi \in \Pi} J_1(H^\pi) \\
& \leqconst N^{\frac{1}{2(1+\alpha)} - r/2} (1+\sup_{\pi \in \Pi} J^2_1(H^\pi))\triangleq (II).
\end{align*}
Then under the event $E_N$, we have 
\begin{align*}
\Rem (\pi) = \abs{\PN f(\Hn, \Hpi)} \leq   \sup_{f \in \F_0}\abs{\PN f}
\end{align*}
where $\F_0$ is given by 
\begin{align*}
\F_0 = & \Bigdkh{ f: (S, A, S') \mapsto  \big[H(S, A) - \sum_{a'} \pi(a'|S')    H(S', a')\big] h(S, A): \pi \in \Pi, h \in \G, H \in \F, \\
	& \hspace{10ex} J_1(H) \leqconst (II),  J_2(h) \leqconst \sup_{\pi \in \Pi} J_2(e^\pi), P f^2 \leqconst (I)}
\end{align*}
Apply Lemma \ref{lemma: donsker class prob bound} with $v = VC(\Pi) + 1$, $\sigma^2 = (I)$ and $M = (II)$, with probability at least $1-\delta - 1/N$, $$
\sup_{f \in \F_0}\abs{\PN f}  \leqconst  (VC(\Pi)+1)N^{-\frac{r+\frac{1}{1+\alpha}}{2}}\left(1+\sup_{\pi \in \Pi}J_1^2(\tilde H^\pi)\right)^{\frac{1+\alpha}{2}}\log^\frac{3}{2}(\max(\frac{1}{\delta}, N)),
$$
%	where
%	\begin{align*}
%	(III) &  = \sqrt{\frac{VC(\Pi) + 1}{N\log(N)}}(I)^\frac{1-\alpha}{2} + \frac{VC(\Pi)+1}{N}(I)^{-\alpha} + \sqrt{\frac{(I)}{N}} \\
%	 & +(II)^\alpha\sqrt{\frac{VC(\Pi+1) + 1}{N\log(N)}}(I)^\frac{1-\alpha}{2} + (II)^{2\alpha}\frac{VC(\Pi)+1}{N}(I)^{-2\alpha} + \frac{\log(N)}{N}\\
%	&\leqconst (VC(\Pi)+1)^{\frac{1}{1+\alpha}}N^{-\frac{r+\frac{1}{1+\alpha}}{2}}\left(1+\sup_{\pi \in \Pi}J_1^2(\tilde H^\pi)\right)^{\frac{1+\alpha}{2}}\log^\frac{3}{2}(\max(\frac{1}{\delta}, N))
%	\end{align*}
%		where $Z_1(\delta), Z_2(\delta)$ only depends on \peng{bla bla}. 
%\begin{align*}
%\tilde C_3(\delta) = C(\delta)\left[ (1 + \tilde{C}_2(\delta)^{\alpha}) \tilde{C}_1(\delta)^{\frac{1-\alpha}{2}} + \frac{1+\tilde{C}_2(\delta)^{2\alpha}}{\tilde{C}_1(\delta)^{\alpha}} + \sqrt{\tilde{C}_1(\delta)} + 2+\tilde{C}_2(\delta)^{\alpha/2} + \tilde{C}_2(\delta)^\alpha\right]
%\end{align*}
Combing with (\ref{J iter}), which holds under the event $E_N$, we have 
\begin{align*}
\lambda_N J_1^2(\Hn) &\leqconst \left(1+VC(\Pi)\right)(1+ \sup_{\pi \in \Pi}J_1^2(H^\pi))^{\frac{\alpha}{1+\alpha}}\left[\log\left(\max(1/\delta, N)\right)\right]^{\frac{2+\alpha}{1+\alpha}} N^{-\frac{1}{1+\alpha}} \notag \\
& +  (VC(\Pi)+1)N^{-\frac{r+\frac{1}{1+\alpha}}{2}}\left(1+\sup_{\pi \in \Pi}J_1^2(\tilde H^\pi)\right)^{\frac{1+\alpha}{2}}\log^\frac{3}{2}(\max(\frac{1}{\delta}, N))+ \mu_NJ_1(\Hn).
\end{align*}
Thus using the same argument as before gives 
\begin{align*}
J_1^2(\Hn) & \leqconst (1+\sup_{\pi \in \Pi}J_1^2(\tilde H^\pi))^{\frac{(1+\alpha)^2+2}{2(1+\alpha)}}N^{\frac{ \frac{1}{1+\alpha}}{2}-\frac{r}{2}}.
\end{align*}
Now using (\ref{err iter}) again with this inequality under the event $E_N$ gives
\begin{align}
\norm{\h(\Hn) - \h(\Hpi)}^2 &\leqconst \left(1+VC(\Pi)\right)(1+ \sup_{\pi \in \Pi}J_1^2(H^\pi))^{\frac{\alpha}{1+\alpha}}\left[\log\left(\max(1/\delta, N)\right)\right]^{\frac{2+\alpha}{1+\alpha}} N^{-\frac{1}{1+\alpha}} \notag \\
& +  \Rem(\pi) + \mu_NJ_1(\Hn) \notag\\
& \leqconst \left(1+VC(\Pi)\right)(1+ \sup_{\pi \in \Pi}J_1^2(H^\pi))^{\frac{1+\alpha}{2}}\left[\log\left(\max(1/\delta, N)\right)\right]^{\frac{2+\alpha}{1+\alpha}} N^{-\frac{r+\frac{1}{1+\alpha}}{2}} \notag
\end{align}
%		Let  
%		$$C^{[k+1]}(\delta) = \inblue{K_4} \bar C(\delta) (1 + \sqrt{C(\delta)})$$
%		$$\beta_{k+1} = \frac{\beta + 1/(1+\alpha)}{2} = \frac{1}{1+\alpha} -  \frac{(1-\alpha) 2^{-k}}{1+\alpha} $$ $$\omega_{k+1} = (1+\omega)/2 = 1 - 2^{-k}$$
Since $(1+\sup_{\pi \in \Pi}J_1^2(\tilde H^\pi)) < \infty$, we show that the following holds w.p. $1-(4+1)\delta - (1+1)/N$, for all $\pi \in \Pi$
\begin{align*}
& \norm{\h(\Hn) - \h(\Hpi)}_2^2 \leqconst \left(1+VC(\Pi)\right)\left[\log\left(\max(1/\delta, N)\right)\right]^{\frac{2+\alpha}{1+\alpha}} N^{-\frac{r+\frac{1}{1+\alpha}}{2}}
\\
& J_1^2(\Hn) \leqconst N^{\frac{ \frac{1}{1+\alpha}}{2}-\frac{r}{2}}. %=  C_{new}  (\delta) \iota^{\omega_{k+1}} p n^{1/(1+\alpha) - \beta_{new}}
\end{align*}
%\paragraph{Final Bound}
%\begin{align*}
%& \beta_{k} =  \frac{1}{1+\alpha} -  \frac{1-\alpha }{1+\alpha} 2^{-(k-1)}\\
%&\gamma_{k} = \frac{1-\alpha}{2(1+\alpha)} 2^{-(k-1)} \\
%& \omega_{k} = 1 - 2^{-(k-1)},  
%\end{align*}
Thus the convergence rate is improved to $r_{2} = \frac{r_1 + 1/(1+\alpha)}{2}$. The same procedure can be applied $m$ times. It is easy to verify that  for any $m \geq 2$,  
$r_{m} = \frac{r_{m-1}+ 1/(1+\alpha)}{2} = \frac{1}{1+\alpha} -  \frac{(1-\alpha) 2^{-(m-1)}}{1+\alpha} $, thus we obtain the desired result. 

%	\newpage
%		\peng{
%			} 

\textbf{Proof of Theorem \ref{thm: ratio} in the Main Text}

Recall the ratio estimator, $\widehat \omega^\pi_N$ in (\ref{ratio estimator}).  From Lemma \ref{thm: pre-ratio} and Lemma \ref{lemma: uniform l2 bound}, w.p. $1-(3+k)\delta - k/N$ for all $\pi \in \Pi$ we have
\begin{align*}
& \norm{\widehat e^\pi_N - e^\pi}^2  = \norm{\hn(\Hn) - \h(\Hn) + \h(\Hn) -  \h(\Hpi)}^2 \\
& \leq 2\norm{\hn(\Hn) - \h(\Hn)}^2 + 2\norm{\h(\Hn) -  \h(\Hpi)}^2 \\
& \leqconst N^{-r_k} (VC(\Pi)+ 1)(\log(\max(N, 1/\delta)))^{\frac{\alpha+2}{\alpha+1}}.
\end{align*}
In addition
\begin{align*}
\abs{\PN \hn(\Hn) - P e^\pi} & = \abs{\PN \hn(\Hn) - P \hn(\Hn) + P \hn(\Hn) - e^\pi} \\
& \leq  \abs{(\PN -P) \hn(\Hn)} + P \abs{ \hn(\Hn) - e^\pi}  \\
& \leq \abs{(\PN -P) \hn(\Hn)} + \norm{\hn(\Hn) - e^\pi}.
\end{align*}
%			We use the empirical process, Lemma \ref{lemma: donsker class prob bound}, to bound the first term.   \peng{At the minimum we should get the same rate as the second term. }
For the first term above, $\abs{(\PN -P) \hn(\Hn)}\leq  J_2(\hn(\Hn)) \frac{\PN \hn(\Hn)}{J_2(\hn(\Hn))} \leq  J_2(\hn(\Hn))  \sup_{h \in\G_1}  \abs{\PN h}$,
where $\G_1 = \{ h: h \in \G,  J_2(h) \leq 1 \}$.  Using Lemma \ref{lemma: uniform l2 bound} with proper chosen $\tau$ as before and Letting $N$ sufficiently large, with probability at least $1-\delta$ for all $\pi$, 
\begin{align*}
& J_2(\hn(\Hn))  \leqconst 1+J_1(\Hn) + J_2(\h (\Hn)) + \frac{\sqrt{VC(\Pi)+1}\left[\log(\max(N, 1/\delta))\right]^{\frac{1}{2\tau}}}{\sqrt{N\mu_N}} + \frac{\sqrt{1 + VC(\Pi)}}{\sqrt{N}\mu_N^{\frac{1 +\alpha - \tau(2+\alpha)}{2(1- \tau(2+\alpha))}}}\\
& \leqconst \left\{N^{\frac{ \frac{1}{1+\alpha}-r_k}{2}}(1+VC(\Pi))[\log(\max(N, 1/\delta))]^{\frac{2+\alpha}{1+\alpha}}\right\}^{\frac{1}{2}},
\end{align*}
Similar to the derivation of initial rate in the proof of Theorem \ref{thm: pre-ratio}, applying  Lemma \ref{lemma: donsker class prob bound} implies with probability at least $1-(4+k)\delta - \frac{k+1}{N}$,  $$\abs{\PN \hn(\Hn) - P e^\pi} \leq \bar C N^{-\frac{r_k }{2}}(1+VC(\Pi))[\log(\max(N, 1/\delta))]^{\frac{2+\alpha}{1+\alpha} + 1},$$
for some constant $\bar C$.
When $N$ is large enough such that $$
\bar{C}N^{-\frac{1+r_k - \frac{1}{1+\alpha}}{2}}(1+VC(\Pi))^{\frac{3}{2} }[\log(\max(N, 1/\delta))]^{\frac{2+\alpha}{2(1+\alpha)} + 1} < (1/2) P e^\pi,
$$
then $\abs{\PN \hn(\Hn) - P e^\pi} \leq (1/2) P e^\pi$. Thus we have $\PN \hn(\Hn) \geq (1/2) P e^\pi > 0$ and 
\begin{align*}
\norm{\widehat \omega^\pi_N - \omega^\pi} & \leq \frac{\norm{\widehat e^\pi_N - e^\pi}}{\abs{\PN \widehat e^\pi_N}} + \norm{e^\pi} \times \Bigabs{\frac{1}{\PN \widehat e^\pi_N} - \frac{1}{P e^\pi}} \\
& \leq \frac{2}{P e^\pi} \norm{\widehat e^\pi_N - e^\pi} +\frac{2 \norm{e^\pi} }{(Pe^\pi)^2} \times   \abs{\PN \hn(\Hn) - P e^\pi} 
\end{align*}
Recall that $\norm{\omega^\pi}^2 = \int \omega^\pi(s, a) d^\pi(s, a) > 1$. As a result, $P e^\pi  = \frac{1}{\norm{\omega^\pi}^2} = \norm{e^\pi}^2$. Finally we have 
\begin{align*}
\norm{\widehat \omega^\pi_N - \omega^\pi}  
& \leq 2 \norm{\omega^\pi}^2 \cdot \norm{\widehat e^\pi_N - e^\pi} + 2 \norm{\omega^\pi}^3 \cdot   \abs{\PN \hn(\Hn) - P e^\pi} \\ 
& \leq 2 \norm{\omega^\pi}^3  \left(\norm{\widehat e^\pi_N - e^\pi} + \abs{\PN \hn(\Hn) - P e^\pi}\right) \\
& \leqconst \sup_{\pi \in \Pi} \{\norm{\omega^\pi}^3\} N^{-r_k/2} (1+VC(\Pi))[\log(\max(N, 1/\delta))]^{\frac{2+\alpha}{2(1+\alpha)} + 1}		
\end{align*}

%	\subsection{Lower level}

\begin{lemma}[Lower level]
	\label{lemma: uniform l2 bound}
	Suppose Assumptions \ref{ass: stationarity}, \ref{ass: policy class}, \ref{assumption: function classes} and \ref{assumption: ratio function} hold. Then with sufficiently large $N$ and $0 \leq \tau \leq \frac{1}{3}$, $\Pr(E_N) \geq 1-\delta$, where the event $E_N$ is that  for all $(H, \pi) \in \F \times \Pi$, the followings hold
	%		with probability at least $1-\delta$, the following event,  $E_n$, hold  and $\alpha$ : 
	\begin{align*}
	& \norm{\hn (H) - \h (H)}_{}^2 \leqconst \mu_N (1+J_1^2(H) + J_2^2(\h (H))) + \frac{(VC(\Pi)+1)\left[\log(\max(N, 1/\delta))\right]^{\frac{1}{\tau}}}{N} + \frac{1 + VC(\Pi)}{N\mu_N^{\frac{\alpha}{1- \tau(2+\alpha)}}}\\
	&  J_2^2(\hn(H))\leqconst  1+J_1^2(H) + J_2^2(\h (H)) + \frac{(VC(\Pi)+1)\left[\log(\max(N, 1/\delta))\right]^{\frac{1}{\tau}}}{N\mu_N} + \frac{1 + VC(\Pi)}{N\mu_N^{\frac{1 +\alpha - \tau(2+\alpha)}{1- \tau(2+\alpha)}}}\\
	& \norm{\hn(H) - \h(H)}_{N}^2 \leqconst \mu_N (1+J_1^2(H) + J_2^2(\h (H))) + \frac{(VC(\Pi)+1)\left[\log(\max(N, 1/\delta))\right]^{\frac{1}{\tau}}}{N} + \frac{1 + VC(\Pi)}{N\mu_N^{\frac{\alpha}{1- \tau(2+\alpha)}}}.
	\end{align*}
\end{lemma}
The proof is similar to that of Lemma \ref{lemma: upper bound of eta, Q} so we omit details here.

\begin{lemma}[Decomposition] \label{lemma: second term} \label{lemma: upper level initial bound}
	Suppose Assumptions \ref{ass: stationarity}, \ref{ass: policy class}, \ref{assumption: function classes} and \ref{assumption: ratio function} hold. Then, the following hold with probability at least $1-{2\delta}$: for all policy $\pi \in \Pi$:
	\begin{align*}
	& \norm{\hn(\Hn) - \h(\Hpi)}^2 + \lambda_N J_1^2(\Hn) \\
	&\leqconst  (\mu_N+\lambda_N) (1+\sup_{\pi \in \Pi}J_1^2(\tilde H^\pi)) +\frac{(VC(\Pi)+1)\left[\log(\max(N, 1/\delta))\right]^\frac{1}{\tau}}{N} + \frac{1+VC(\Pi)}{N\mu_N^{\frac{\alpha}{1-\tau(2+\alpha)}}}\\
	& +\Rem(\pi) +  \mu_N J_1(\Hn) + 
	\frac{\sqrt{\mu_N(VC(\Pi)+1)}\left[\log(\max(N, 1/\delta))\right]^{\frac{1}{2\tau}}}{\sqrt{N}} + \frac{\sqrt{1 + VC(\Pi)}}{\sqrt{N}\mu_N^{\frac{\alpha + \tau(2+\alpha)- 1}{2(1- \tau(2+\alpha))}}}\\
	& + N^{-\frac{1-(2+\alpha)\tau}{1+\alpha-(2+\alpha)\tau}} + \frac{(VC(\Pi)+1)\log^{\frac{\alpha/\tau}{1+\alpha - \tau(2+\alpha)}}(\max(N,1/\delta))}{N\mu_N^{\alpha/(1+\alpha- (2+\alpha)\tau)}} +\frac{1}{N \lambda_N^{\frac{\alpha}{1 - \tau(2+ \alpha)}}}
	\end{align*}
	%	\begin{align*}
	%		\norm{\gn(\Hn) - \g(\Hpi)}^2 + \lambda_n J_1^2(\Hn)  \leq \peng{bla bla} +  \abs{\Pn  \g(\Hpi) Y_{\Hn} - Y_{\Hpi}] }
	%%	\norm{\gn(\Hn)}_{}^2 + \norm{\gn(\Hn)}_{n}^2& \leqconst n^{-\frac{1}{1+\alpha}} + \mu_n (1+J_1^2(\tilde Q^\pi))  +  \lambda_n J_1^2(\tilde Q^\pi)   + (n\lambda_n^{\alpha})^{-1}  \\
	%%	& + p (n \mu_n^{\alpha})^{-1}   + (n\mu_n^{\alpha/(1+\alpha)})^{-1}+  \log(1/\delta) (n^{-1} + (n\mu_n^{\alpha/(1+\alpha)})^{-1})
	%	\end{align*}
	where $\Rem(\pi) =4\bigabs{\PN \h(S; A; \Hpi) [\Delta^\pi(S, A, S'; \Hn) - \Delta^\pi(S, A, S; \Hpi)] }$
\end{lemma}
\textbf{Proof of Lemma \ref{lemma: upper level initial bound}}
For $h_1, h_2  \in \G$,  define the functionals $\f_1, \f_2$, 
\begin{align*}
& \f_1(h_1) (S, A) =  h^2_1(S, A)\\
& \f_2(h_1, h_2) (S, A) = 2 h_1(S, A)h_2(S, A) 
%	\\
%	& \f_3(g_1, g_2, g_3) (D) = (1/T) \sum_{t=1}^T g_3(S_t, A_t) (g_1 - g_2)(S_t, A_t)
%	f_2(g_1, g_2) (D) = (1/T) \sum_{t=1}^T (g_1(S_t, A_t)-g_2(S_t, A_t))^2
\end{align*}
With this definition, we have
\begin{align*}
& \norm{\hn(\Hn) - \g(\Hpi)}^2 + \lambda_N J_1^2(\Hn) 
= P\f_1(\hn(\Hn)-\h(\Hpi)) + \lambda_N J_1^2(\Hn) \\
&= 2 \times \big[\PN  \f_1(\hn(\Hn)) + \lambda_NJ_1^2(\Hn)\big] + P\f_1(\hn(\Hn)-\h(\Hpi))  \\
& \qquad \qquad + \lambda_N J_1^2(\Hn) - 2 \times \big[\PN  \f_1(\hn(\Hn)) + \lambda_N J_1^2(\Hn)\big] 
%\\& \leq 2 (\Pn  \f_1(\gn(\Hpi)) + \lambda_nJ_1^2(\Hpi)) + P\f_1(\gn(\Hn)-\g(\Hpi)) + \lambda_n J_1^2(\Hn) - 2 (\Pn  \f_1(\gn(\Hn))+ \lambda_nJ_1^2(\Hn))\\
\end{align*}	
Using the optimizing property of $\Hn$ in (\ref{qhat}), the first term can be bounded by the following inequality.
\begin{align*}
& \PN  \f_1(\hn(\Hn))+ \lambda_N J_1^2(\Hn) \leq \PN  \f_1(\hn(\Hpi)) + \lambda_N J_1^2(\Hpi)\\
& = \PN  \f_1 (\hn(\Hpi) - \h(\Hpi)) +  \PN \f_1(\h(\Hpi))+  \PN \f_2(\hn(\Hpi) - \h(\Hpi), \h(\Hpi))  + \lambda_NJ_1^2(\Hpi) \\
& = \PN  \f_1(\hn(\Hpi) - \h(\Hpi)) +  \lambda_N J_1^2(\Hpi) +  (1/2) \PN \f_2(2\hn(\Hpi) - \h(\Hpi), \h(\Hpi))
\end{align*}
so that
\begin{align*}
&  \norm{\hn(\Hn) - \h(\Hpi)}^2 + \lambda_N J_1^2(\Hn) \\
&  \leq 2 \left[\PN  \f_1(\hn(\Hpi) - \h(\Hpi)) +  \lambda_N J_1^2(\Hpi) + (1/2) \PN \f_2(2\hn(\Hpi) - \h(\Hpi), \h(\Hpi))  \right]  \\
& \qquad + P\f_1(\hn(\Hn)-\h(\Hpi)) + \lambda_N J_1^2(\Hn) - 2 (\PN  \f_1(\hn(\Hn))+ \lambda_N J_1^2(\Hn)) \\
& = 2 \left[\Pn  \f_1(\gn(\Hpi) - \h(\Hpi)) +  \lambda_N J_1^2(\Hpi) \right]   + P\f_1(\hn(\Hn)-\h(\Hpi)) - \lambda_N J_1^2(\hn) \\
& \qquad - 2 \PN [ \f_1(\hn(\Hn)) +  (1/2)  \f_2(\h(\Hpi) - 2\hn(\Hpi), \h(\Hpi))] \\
& = 2 \left[\PN  \f_1(\hn(\Hpi) - \h(\Hpi)) +  \lambda_NJ_1^2(\Hpi) \right]  + P\f_1(\hn(\Hn)-\g(\Hpi)) - \lambda_N J_1^2(\Hn) \\
&	\qquad - 2 \PN [ \f_1(\hn(\Hn)-\h(\Hpi)) + \f_2 (\hn(\Hn)- \hn(\Hpi) , \h(\Hpi))]\\
& \leq 2 \left[\PN  \f_1(\hn(\Hpi) - \h(\Hpi)) +  \lambda_N J_1^2(\Hpi) \right]  \\
& \qquad + P\f_1(\hn(\Hn)-\h(\Hpi)) - \lambda_N J_1^2(\Hn) - 2 \Pn  \f_1(\hn(\Hn)-\h(\Hpi))\\
&	 \qquad + 2\abs{ \PN \f_2 (\hn(\Hn)- \hn(\Hpi) \h(\Hpi))}.
\end{align*}
Then we decompose the error into three components.
\begin{align*}
\norm{\hn(\Hn) - \h(\Hpi)}^2 + \lambda_N J_1^2(\Hn) = I_{1}(\pi) +  I_2(\pi) + I_3(\pi)
\end{align*}
where
\begin{align*}
& I_1 (\pi) = 2 \left[\PN  \f_1(\hn(\Hpi) - \h(\Hpi)) +  \lambda_N J_1^2(\Hpi) \right] \\
& I_2(\pi) = P\f_1(\hn(\Hn)-\h(\Hpi)) - \lambda_N J_1^2(\Hn) - 2 \PN  \f_1(\hn(\Hn)-\hn(\Hpi))\\
& I_3(\pi) = 2\abs{ \PN \f_2 (\hn(\Hn)- \hn(\Hpi) , \h(\Hpi))}
\end{align*}
Below we provide a bound for each of the three terms. By Lemma \ref{lemma: uniform l2 bound} with $\tau_1 \leq \frac{1}{3}$ and sufficiently large $N$, three inequalities in Lemma \ref{lemma: uniform l2 bound} hold. Denote such event as $E_N$. 
%		\begin{align*}
%		E_n = \{ bla bla  \}
%		\end{align*}
\paragraph{Bounding $I_1(\pi)$} Under the event $E_N$, we can have
\begin{align*}
I_1(\pi) & =  2\norm{\hn (\Hpi) - \h (\Hpi)}_{}^2 + \lambda_N J_1^2(\Hpi) \\
&	\leqconst \Big(\mu_N (1+J_1^2(\Hpi) + J_2^2(\h (\Hpi)))\Big)  +\frac{(VC(\Pi)+1)\left[\log(\max(N, 1/\delta))\right]^\frac{1}{\tau_1}}{N} \\
& + \frac{1+VC(\Pi)}{N\mu_N^{\frac{\alpha}{1-\tau(2+\alpha)}}}+ \lambda_N J_1^2(\Hpi) \\
&\leqconst	(\mu_N+\lambda_N) (1+\sup_{\pi \in \Pi}J_1^2(\tilde H^\pi)) +\frac{(VC(\Pi)+1)\left[\log(\max(N, 1/\delta))\right]^\frac{1}{\tau_1}}{N} + \frac{1+VC(\Pi)}{N\mu_N^{\frac{\alpha}{1-\tau(2+\alpha)}}}.
\end{align*}

\paragraph{Bounding $I_3(\pi)$}

%	where we define $f_2(g_1, g_2, g_3) = 2 g_3(g_1 - g_2)$. 
%\[
%f_2 (\gn(\Hn), \gn(\Hpi) , \g(\Hpi)) =  2 \g(\Hpi) [ \gn(\Hn)     - \gn(\Hpi)]
%\]
%	Thus we have 
%	\begin{align*}
%	& \norm{\gn(\Hn) - \g(\Hpi)}^2 + \lambda_n J_1^2(\Hn) \\
%	& \leq 2 \left[\Pn  \f_1(\gn(\Hpi) - \g(\Hpi)) +  \lambda_nJ_1^2(\Hpi) \right] + P\f_1(\gn(\Hn)-\g(\Hpi)) - \lambda_n J_1^2(\Hn) - 2 \Pn  (f_1 +f_2)  
%	\end{align*}
%Let $\tilde f_2 = f_2 - P f_2$ and $f = f_1 + \tilde f_2$, then we have $\EE[f^2] \leqconst \EE[f_1]  = \EE[f]$.  The above is then 
%\begin{align*}
%& \norm{\gn(\Hn) - \g(\Hpi)}^2 + \lambda_n J_1^2(\Hn) \\
% & \leq 2 \left[\Pn  \f_1(\gn(\Hpi) - \g(\Hpi)) +  \lambda_nJ_1^2(\Hpi) \right] + Pf(\gn(\Hn), \g(\Hpi)) - \lambda_n J_1^2(\Hn) - 2 \Pn f - 2 Pf_2
%\end{align*}
%Now we introduce $Y_{\widehat q}$ (again this would be problematic)
%\begin{align*}
%(1/2) f_2 (\gn(\Hn), \gn(\Hpi) , \g(\Hpi) = \Pn \g(\Hpi) [ \gn(\Hn)     - \gn(\Hpi)] = \Pn \g(\Hpi) [ Y_{\Hn} - Y_{\Hpi}] + \Rem
%\end{align*}
%\peng{If there is any hope that $\Pn f_2$ is of the optimal order -- we must be able to show $P f_2$ is of the same order. So we are definitely making the life easier.  }
%For $P f_2$. We have for any $q$:
%\begin{align*}
%& \innerprod{\gn(q) - \g(H)}{\g(H)}_2 \\ 
%& =  \innerprod{\gn(q) - Y_{q} + Y_{q} - \g(H)}{\g(\Hpi)}_2 \\ 
%& =  \innerprod{\gn(q) - Y_{q}}{\g(\Hpi)}_2 \\
%& = P((\gn(q) - Y_{q})\g(\Hpi)) \\
%& = (P-\Pn) [(\gn(q) - Y_{q})\g(\Hpi)] + \Pn  [(\gn(q) - Y_{q})\g(\Hpi)]
%\end{align*}

%	\peng{$Y_q, \Pn $...} 
Using the optimizing property of $\hn(H)$ and Assumption \ref{b3-2} that $\h(\Hpi) = e^\pi \in \G$,  the followings holds for all $H \in \F, \pi \in \Pi$, 
%	\begin{align*}
%	\Pn  [(1/T) \sum_{t=1}^T (\gn(q)(S_t, A_t) - (1-H(S_t, A_t) + \sum_{a'} \pi(a'|S_{t+1}) H(S_{t+1}, a')))\g(\Hpi)] = - \mu_n J_2( \gn(q), \g(\Hpi)) 
%	\end{align*}
\begin{align*}
& \mu_N J_2( \hn(H), \h(\Hpi))  \\
& = \Pn  [(1/T) \sum_{t=1}^T \big(1-H(S_t, A_t) + \sum_{a'} \pi(a'|S_{t+1}) H(S_{t+1}, a') - \hn(S_t, A_t; H)\big)\h(S_t, A_t; \Hpi)] \\
& = \Pn  [(1/T) \sum_{t=1}^T (\Delta^\pi(S_t, A_t, S_{t+1}; H) - \hn(S_t, A_t; H))\h(S_t, A_t; \Hpi)]
\end{align*}
Thus we have 
\begin{align*}
& (1/2) \PN \f_2 (\hn(\Hn)- \hn(\Hpi) , \h(\Hpi)) \\
& =   \Pn (1/T) \sum_{t=1}^T  \h(S_t, A_t; \Hpi) [ \hn(S_t, A_t; \Hn)     - \hn(S_t, A_t; \Hpi)] \\
& =  \Pn (1/T) \sum_{t=1}^T  \h(S_t, A_t; \Hpi) [ \hn(S_t, A_t; \Hn)    - \Delta^\pi(S_t, A_t, S_{t+1}; \Hn) + \Delta^\pi(S_t, A_t, S_{t+1}; \Hn) \\
& \qquad\qquad  - \Delta^\pi(S_t, A_t, S_{t+1}; \Hpi) + \Delta^\pi(S_t, A_t, S_{t+1}; \Hpi) - \hn(S_t, A_t; \Hpi)] \\
%	& =   \Pn  \g(\Hpi) [ \gn(\Hn)   - Y_{\Hn} + Y_{\Hn} - Y_{\Hpi} + Y_{\Hpi}- \gn(\Hpi)] \\
& =   \Pn (1/T) \sum_{t=1}^T \h(S_t; A_t; \Hpi) [\Delta^\pi(S_t, A_t, S_{t+1}; \Hn) - \Delta^\pi(S_t, A_t, S_{t+1}; \Hpi)] \\
& \qquad\qquad  + \mu_N J_2( \hn(\Hpi), \h(\Hpi)) - \mu_N J_2(\hn(\Hn), \h(\Hpi))
%	& =   \Pn  \g(\Hpi) [Y_{\Hn} - Y_{\Hpi}] + \peng{\mu_n J_2( \gn(\Hpi)  - \gn(\Hn), \g(\Hpi))}
%& = \Pn (Y_{\Hn} - Y_{\Hpi}) \g(\Hpi)  - 2\mu_n J_2( \gn(q), \g(\Hpi))  \\& = \Pn f_3 - 2\mu_n J_2( \gn(q), \g(\Hpi))
%\Pn f_2 & =  2 \Pn  \g(\Hpi) [ \gn(\Hn)     - \gn(\Hpi)] \\
%& = \Pn (Y_{\Hn} - Y_{\Hpi}) \g(\Hpi)  - 2\mu_n J_2( \gn(q), \g(\Hpi))  \\& = \Pn f_3 - 2\mu_n J_2( \gn(q), \g(\Hpi))
\end{align*}
In addition, under the event $E_N$, we have
%		\[
%		J_2(\gn(q))\leqconst  1 + J_1(q) +  J_2(\g (q)) + \sqrt{\frac{p}{n \mu_n^{\alpha+1}} }+ \sqrt{\frac{1}{n \mu_n}} +  \sqrt{\frac{\log(1/\delta)}{n\mu_n}}\\
%		\]
\begin{align*}
& \abs{\mu_N J_2( \hn(\Hpi), \h(\Hpi)) - \mu_N J_2(\hn(\Hn), \h(\Hpi))}  \leq \mu_N J_2(e^\pi) \left(J_2(\hn(\Hpi) ) + J_2(\hn(\Hn))\right) \\
%	& \leq \mu_n (\sup_{\pi \in \Pi} J_2(e^\pi) ) \left(J_2(\gn(\Hpi)) + J_2(\gn(\Hn))\right) \\
& \leqconst  \mu_NJ_2(e^\pi)\Big(1 + J_1(\tilde H^\pi) + J_1(\Hn) + \frac{\sqrt{VC(\Pi)+1}\left[\log(\max(N, 1/\delta))\right]^{\frac{1}{2\tau}}}{\sqrt{N\mu_N}} + \frac{\sqrt{1 + VC(\Pi)}}{\sqrt{N}\mu_N^{\frac{1 +\alpha - \tau(2+\alpha)}{2(1- \tau(2+\alpha))}}}\Big)\\
& \leqconst  \mu_N(1 + \sup_{\pi \in \Pi} J_2(H^\pi))\Big(1 + \sup_{\pi \in \Pi}J_1(\tilde H^\pi) + J_1(\Hn) + \frac{\sqrt{VC(\Pi)+1}\left[\log(\max(N, 1/\delta))\right]^{\frac{1}{2\tau}}}{\sqrt{N\mu_N}} + \frac{\sqrt{1 + VC(\Pi)}}{\sqrt{N}\mu_N^{\frac{1 +\alpha - \tau(2+\alpha)}{2(1- \tau(2+\alpha))}}}\Big)\\
& \leqconst  \mu_N J_1(\Hn) + 
\mu_N\big(1 + \frac{\sqrt{VC(\Pi)+1}\left[\log(\max(N, 1/\delta))\right]^{\frac{1}{2\tau}}}{\sqrt{N\mu_N}} + \frac{\sqrt{1 + VC(\Pi)}}{\sqrt{N}\mu_N^{\frac{1 +\alpha - \tau(2+\alpha)}{2(1- \tau(2+\alpha))}}}\Big),
\end{align*}
where the last equality holds by the assumption that $\sup_{\pi \in \Pi} J_2(\tilde H^\pi) < \infty$.
Thus we have
\begin{align*}
I_{3}(\pi) & =  2\abs{ \PN \f_2 (\hn(\Hn)- \hn(\Hpi) , \h(\Hpi))} \\
& \leqconst 4\bigabs{\PN \h(S; A; \Hpi) [\Delta^\pi(S, A, S'; \Hn) - \Delta^\pi(S, A, S; \Hpi)] }  \\
& \qquad +  \mu_N J_1(\Hn) + 
\mu_N\big(1 + \frac{\sqrt{VC(\Pi)+1}\left[\log(\max(N, 1/\delta))\right]^{\frac{1}{2\tau}}}{\sqrt{N\mu_N}} + \frac{\sqrt{1 + VC(\Pi)}}{\sqrt{N}\mu_N^{\frac{1 +\alpha - \tau(2+\alpha)}{2(1- \tau(2+\alpha))}}}\Big)\\
& = \Rem(\pi) +  \mu_N J_1(\Hn) + 
\mu_N\big(1 + \frac{\sqrt{VC(\Pi)+1}\left[\log(\max(N, 1/\delta))\right]^{\frac{1}{2\tau}}}{\sqrt{N\mu_N}} + \frac{\sqrt{1 + VC(\Pi)}}{\sqrt{N}\mu_N^{\frac{1 +\alpha - \tau(2+\alpha)}{2(1- \tau(2+\alpha))}}}\Big),
\end{align*}
where we let $\Rem(\pi) =4\bigabs{\PN \h(S; A; \Hpi) [\Delta^\pi(S, A, S'; \Hn) - \Delta^\pi(S, A, S; \Hpi)] }$.
\paragraph{Bounding $I_2(\pi)$}
%		We can show w.p. $1-2\delta$, for all policy $\pi \in \Pi$, 
%	$$P\f_1(\gn(\Hn)-\g(\Hpi)) - \lambda_n J_1^2(\Hn) - 2 \Pn \f_1(\gn(\Hn)-\g(\Hpi))= O_P(n^{-1/(1+\alpha)})$$

%		\begin{align*}
%		I_2(\pi) 
%		%	 & \leq c_2n^{-1} + c_5 n^{-\frac{1}{1+\alpha}}+ c_5(n\lambda_n^{\alpha})^{-1}  + c_5 \left[\frac{p^{\frac{\alpha}{1+\alpha}}}{n \mu_n^{\alpha}} + \frac{1+\log^{\frac{\alpha}{1+\alpha}}(1/\delta)}{n\mu_n^{\alpha/(1+\alpha)}}\right] + \log(120/\delta) c_6n^{-1}\\
%		&  \leq \inblue{K_2}\Big[ \frac{1}{n^{1/(1+\alpha)}} + \frac{1}{n\lambda_n^\alpha} + \frac{p^{\frac{\alpha}{1+\alpha}}}{n \mu_n^{\alpha}} + \frac{1+\log^{\frac{\alpha}{1+\alpha}}(1/\delta)}{n\mu_n^{\alpha/(1+\alpha)}} + \frac{1+\log(1/\delta)}{n} \Big]
%		\end{align*}
%
%		\peng{A long proof here.}

For the second term, 
\begin{align*}
I_2(\pi) & = P\f_1(\hn(\Hn)-\h(\Hpi)) - \lambda_N J_1^2(\Hn) - 2 \PN  \f_1(\hn(\Hn)-\h(\Hpi)) \\
& = 2(P - \PN) \f_1(\hn(\Hn)-\h(\Hpi)) - \lambda_N J_1^2(\Hn) - P \f_1(\hn(\Hn)-\h(\Hpi))
\end{align*}
Let $\zeta(N, \mu_N, \delta, \tau_1) = 1 + \frac{\sqrt{VC(\Pi)+1}\left[\log(\max(N, 1/\delta))\right]^{\frac{1}{2\tau}}}{\sqrt{N\mu_N}} + \frac{\sqrt{1 + VC(\Pi)}}{\sqrt{N}\mu_N^{\frac{1 +\alpha - \tau(2+\alpha)}{2(1- \tau(2+\alpha))}}}$. Under $E_N$, we have 
\begin{align*}
& J_2(\hn(\Hn)-\h(\Hpi))  \\
& \leq J_2(\hn(\Hn))  + J_2(\h(\Hpi)) \\ 
& \leqconst (1 + J_1(\Hn) +  J_2(\h (\Hn)) + \zeta(N, \mu_N, \delta, \tau_1) + J_1(\tilde H^\pi) \\
& \leqconst J_1(\Hn) +  \zeta(N, \mu_N, \delta, \tau_1) + \sup_{\pi \in \Pi} J_1(\tilde H^\pi) 
\end{align*} 
Now we have $\Pr(\exists \pi \in \Pi,  I_2(\pi) > t) \leq \Pr(\{\exists \, \pi \in \Pi,  I_2(\pi) > t\} \jiao E_N) + \delta$ and we bound the first term using peeling device on $\lambda_n J_1^2(\Hn)$ in $I_2(\pi)$: 
\begin{align*}
& \Pr(\{\exists \pi \in \Pi,  I_2(\pi) > t\} \jiao E_n)  \\
& = \sum_{l=0}^\infty \Pr\big(\{\exists \pi \in \Pi,  I_2(\pi) > t, ~\lambda_N J_1^2(\Hn) \in [2^{l}t \indicator{t\neq 0}, 2^{l+1} t)\} \jiao E_N \big)\\
& \leq \sum_{l=0}^\infty \Pr\big(\exists \pi \in \Pi, ~ 2(P-\PN) \f_1(\hn(\Hn)-\h(\Hpi)) > P \f_1(\hn(\Hn)-\h(\Hpi)) + \lambda_N J_1^2(\Hn) + t, \\
& \hspace{12ex} \lambda_N J_1^2(\Hn) \in [2^{l}t \indicator{t\neq 0}, 2^{l+1} t), J_2(\hn(\Hn)-\h(\Hpi)) \leq c_1 (J_1(\Hn) + \zeta(N, \mu_N, \delta, \tau_1))\big)\\
& \leq  \sum_{l=0}^\infty \Pr\big(\exists \pi \in \Pi, ~2(P-\PN) \f_1(\hn(\Hn)-\h(\Hpi)) > P \f_1(\hn(\Hn)-\h(\Hpi)) + 2^{l}t \indicator{t\neq 0} + t, \\
& \hspace{12ex} \lambda_N J_1^2(\Hn) \leq 2^{l+1} t, J_2(\hn(\Hn)-\h(\Hpi)) \leq c_1 (\sqrt{(2^{l+1} t)/\lambda_n} + \zeta(N, \mu_N, \delta, \tau_1))\big)\\
& \leq  \sum_{l=0}^\infty \Pr\big(\exists \pi \in \Pi, ~2(P-\PN) \f_1(\hn(\Hn)-\h(\Hpi)) > P \f_1(\hn(\Hn)-\h(\Hpi)) + 2^{l}t, \\
& \hspace{12ex} J_2(\hn(\Hn)-\h(\Hpi)) \leq c_1 (\sqrt{ (2^{l+1} t)/\lambda_N} + \zeta(N, \mu_N, \delta, \tau_1))\big)\\
& \leq \sum_{l=0}^\infty \Pr\left( \sup_{f \in \F_l}  \frac{(P-\PN) f(S, A)}{P f(S, A) + 2^l t} > \frac{1}{2} \right),
\end{align*}
where $c_3$ is some constant, and $\F_l = \{ \f_1(h): J_2(h) \leq c_1 ( \sqrt{(2^{l+1} t)/\lambda_N} + \zeta(N, \mu_N, \delta, \tau_1)), h \in \G  \}$. It is easy to see that  $\abs{f(g)(S, A)} \leq G_{\max}^2 \triangleq K_1$.

Similar to Lemma \ref{lemma: second term}, we bound each term of the above probabilities by using the independent block technique. 
For each $l \geq 0$, we will use an  independent block sequence $(x_{N, l}, v_{N, l})$ with the residual $R_{ l}$. By controlling the size of these blocks, we can  optimize the bound. We let
$$
x_{N, l} = \lfloor x'_{N, l} \rfloor \epc \mbox{and } \epc  v_{N, l} = \lfloor \frac{N}{2x_{N, l}} \rfloor,
$$
where $x'_{N, l} = (Nt)^\tau (2^l)^p$ and $v'_{N, l} = \frac{N}{2x'_{N, l}}$ with some positive constants $\tau$ and $p$. Let $\tau \leq p \leq \frac{1}{2 + \alpha} \leq \frac{1}{2}$ and $N$ satisfies the following constraint:
\begin{align}\label{sample size constraint for Q 1}
N \geq c_1 \triangleq 4 \times 8^2 \times K_1 \geq 4^{\frac{p}{1-p}}8^{\frac{1}{1-p}}.
\end{align} 
By the definition of $x'_{N,l}$ and assuming $t \geq \frac{1}{N}$, $x_{N,l}\geq 1$. Then we consider two cases. The first case is  any $l$ such that $x'_{N,l} \geq \frac{N}{8}$. In such case, based on the assumption over $\tau$ and $p$, we can show that $x'_{N, l} \leq (Nt2^l)^p$, which further implies that $(Nt2^l) \geq 4NK_1$ by the sample constraint and $p \leq \frac{1}{2+\alpha}$. Then we can show that for this case,
$$
\frac{(P-\PN) \left\{f(S, A)\right\}}{P \left\{f(S, A)\right\} + 2^l t} \leq  \frac{2K_1}{2^lt} \leq  \frac{1}{2},
$$
for sufficiently large $N$.
Thus such terms does not contribute to the probability bound.

The second case we consider is any $l$ such that $x'_{N, l} < \frac{N}{8}$. We again apply the relative deviation concentration inequality for the exponential $\boldsymbol{\beta}$-mixing stationary process given in Theorem 4 of \cite{farahmand2012regularized}, which combined results in \cite{yu1994rates} and Theorem 19.3 in \cite{gyorfi2006distribution}. It then suffices to verify conditions (C1)-(C5) in Theorem 4 of \cite{farahmand2012regularized} with $\F = \F_l$, $\epsilon = 1/2$ and $\eta = 2^{l} t$ to get an exponential inequality for each term in the summation. The conditions (C1) has been verified. For (C2), we have $P f^2(g) \leq G_{\max}^2 P f(g)$ and thus (A2) holds by choosing $K_2 = G_{\max}^2$

%$K_1 ＝ 6G_{\max}(2R_{\max}+ 2Q_{\max} + 3G_{\max})$ and  $K_2$
%6G_{\max}(2R_{\max}+ 2Q_{\max} + 3G_{\max})$. 
%There exists some $K_1, K_2$ (depending on $\pi_{\min}$, $R_{\max}, Q_{\max}$) such that $\norm{f}_\infty \leq K_1$ and $\EE\left[f(\D)^2\right] \leq K_2 \EE\left[f(\D)\right]$. 

To verify the condition  (C3), without loss of generality, we assume $K_1 \geq 1$. Otherwise, let $K_1 = \max(1, K_1)$. Then we know that $2K_1x_{N, l} \geq \sqrt{2K_1x_{N, l}}$ since $x_{N, l} \geq 1$. We need to have $\sqrt{N} \epsilon\sqrt{1-\epsilon} \sqrt{\eta} \geq 1152K_1x_{N, l}$, or suffice to have $\sqrt{N} \epsilon\sqrt{1-\epsilon} \sqrt{\eta} \geq 1152K_1x'_{N, l}$. Recall that $\epsilon = 1/2$ and $\eta = 2^{l} t$. So it is enough to show that
$$
\sqrt{N}\frac{\sqrt{2}}{4}  \sqrt{2^lt} \geq 1152K_1(Nt2^l)^p.
$$
We can check that if $t \geq \frac{2304\sqrt{2}K_1}{N}$, the above inequality holds for every $l \geq$ since $p \leq \frac{1}{2+\alpha}$.

Next we verify (C4) that $\frac{|R_l|}{N} \leq \frac{\epsilon\eta}{6K_1}$.  Recall that $|R_l| \leq 2 x_{N, l} \leq 2x'_{N, l} = (Nt)^\tau(2^l)^p$. So if $t \geq \frac{c_2}{n}$ for some positive constant $c_2$, we can have
$$
\frac{\epsilon\eta}{6K_1} = \frac{2^lt}{12K_1} \geq \frac{2(Nt)^\tau (2^l)^p}{N} = \frac{2x'_{N,l}}{N}\geq \frac{|R_l|}{N}.
$$
In addition, $|R_l| \leq 2 x'_{N, l} < \frac{N}{2}$.

We now verify the final condition (C5).  First, we obtain an upper bound $\mathcal{N}_{}(u, \F_l; \norm{\cdot}_\infty)$  for all possible realization of $(S, A)$.  For any $g_1, g_2 \in \G$, 
\[
\PP_N\left[f(g_1)(S, A) - f(g_2)(S, A)\right]^2 \leq 4 G_{\max}^2 \norm{g_1 - g_2}_{N}^2.
\]
Thus applying Assumption \ref{assumption: function classes} implies that for some constant $c_3$, the metric entropy for each $l$ is bounded  by
\begin{align*}
& \log {\mathcal{N}}(u, \F_l, \norm{\cdot}_\infty) \\
& \leq \log {\mathcal{N}}\left(\frac{u}{2G_{\max}}, \{h: J_2(h) \leq c_1 ( \sqrt{(2^l t)/\lambda_N }+ \zeta(N, \mu_N, \delta, \tau_1)), h \in \G  \}, \norm{\cdot}_\infty \right)\\
& \leq C_3\left[\frac{c_1 \left\{ \sqrt{(2^l t)/\lambda_N} + \zeta(N, \mu_N, \delta, \tau_1)\right\}}{u/(2G_{\max})}\right]^{2\alpha} \leq c_3 \left\{ \left(\frac{2^l t}{\lambda_N}\right)^\alpha + \zeta(N, \mu_N, \delta, \tau_1)^{2\alpha}\right\} u^{-2\alpha},
\end{align*} 
for some positive constant $c_3$.

Now we see the condition (C5) is satisfied if the following inequality holds for all $x \geq (2^ltx_{N, l})/8$ such that
\begin{align*}
\frac{\sqrt{v_{N, l}} (1/2)^2 x}{96x_{N, l}\sqrt{2} \max(K_1, 2K_2)} 
&\geq \int_{0}^{\sqrt{x}} \sqrt{c_3} \left\{ \left(\frac{2^l t}{\lambda_N}\right)^\alpha + \zeta(N, \mu_N, \delta, \tau_1)^{2\alpha}\right\}^{1/2}  \left(\frac{u}{2x_{N, l}}\right)^{-\alpha} du \\
&= x_{N, l}^{\alpha} x^{\frac{1-\alpha}{2}} \sqrt{2^\alpha c_3} \left( \left(\frac{2^l t}{\lambda_N}\right)^\alpha + \zeta(N, \mu_N, \delta, \tau_1)^{2\alpha}\right)^{1/2}.
\end{align*}
It is sufficient to let the following inequality hold:
\begin{align*}
\frac{\sqrt{v_{N, l}}  }{384x_{N, l}\sqrt{2} \max(K_1, 2K_2)}  x^{\frac{1+\alpha}{2}} \geq  \sqrt{c_3'}x_{N, l}^{\alpha}  \left\{ \left(\frac{2^l t}{\lambda_n}\right)^\alpha + \zeta(N, \mu_N, \delta, \tau_1)^{2\alpha} \right\}^{1/2},
\end{align*}
for some constant $c_3'$.
Using the inequality that $(a+b)^{1/2} \leq \sqrt{a} + \sqrt{b}$ and the fact that LHS is increasing function of $x$, it's enough to ensure that the following two inequalities hold:
\begin{align*}
& \frac{\sqrt{v_{N, l}}  }{384x_{N, l}\sqrt{2} \max(K_1, 2K_2)}  (x_{N, l}2^l t/8)^{\frac{1+\alpha}{2}} \geq  \sqrt{c_3'} x_{N, l}^{\alpha} \left(\frac{2^l t}{\lambda_N}\right)^{\alpha/2}\\
& \frac{\sqrt{v_{N, l}}  }{384x_{N, l}\sqrt{2} \max(K_1, 2K_2)}  (x_{N, l}2^l t/8)^{\frac{1+\alpha}{2}} \geq \sqrt{c_3'}x_{N, l}^{\alpha}  \zeta(N, \mu_N, \delta, \tau_1)^\alpha.
\end{align*}
By the definition of $v_{N, l}$ and $x_{N, l}$, after some algebra, we can see that the first inequality holds if
$$
t \geq c_5 \frac{1}{N \lambda_N^{\frac{\alpha}{1 - \tau(2+ \alpha)}}}.
$$
The second inequality holds if
$t$ satisfies
\begin{align*}
t & \geq  c_6 N^{-\frac{1-(2+\alpha)\tau}{1+\alpha- (2+\alpha)\tau}}\zeta(N, \mu_N, \delta, \tau_1)^{\frac{2\alpha}{1+\alpha - (2+\alpha)\tau}}.
\end{align*}
Choosing $\tau = \tau_1 \leq 1/3$, we can obtain that
\begin{align*}
&N^{-\frac{1-(2+\alpha)\tau}{1+\alpha-\tau(2+\alpha)}}\zeta(N, \mu_N, \delta, \tau_1)^{\frac{2\alpha}{1+\alpha - (2+\alpha)\tau}}\\ 
& = N^{-\frac{1-(2+\alpha)\tau}{1+\alpha-\tau(2+\alpha)}} \left[1+ \frac{(VC(\Pi)+1)\left[\log\left(\max(1/\delta, N)\right)\right]^{\frac{1}{\tau_1}}}{N\mu_N} + \frac{1+VC(\Pi)}{N \mu_N^{\frac{1 - \tau_1(2+ \alpha) + \alpha}{1 - \tau_1(2+ \alpha)}}} \right]^{\frac{\alpha}{1+\alpha - (2 +\alpha)\tau}} \\
& \leqconst N^{-\frac{1-(2+\alpha)\tau}{1+\alpha-\tau(2+\alpha)}} + \frac{1+VC(\Pi)}{N \mu_N^{\alpha/(1-\tau(2+\alpha))}} + \frac{(VC(\Pi)+1)\log^{\frac{\alpha/\tau}{1+\alpha - \tau(2+\alpha)}}(\max(N,1/\delta))}{N\mu_N^{\alpha/(1+\alpha- (2+\alpha)\tau)}}. 
%& \leqconst N^{-\frac{1-(2+\alpha)\tau}{1+\alpha-\tau(2+\alpha)}} + \frac{\log^{\frac{\alpha/\tau}{1+\alpha - \tau(2+\alpha)} }(\max(N,1/\delta))+ VC(\Pi)}{N \mu_N^{\alpha/(1-\tau(2+\alpha))}}. \mbox{\QZL{Don't proceed to this step}}
\end{align*}

Putting all together, all conditions (C1) to (C5) would be satisfied for all $l \geq 0$ when
\begin{align*}
t \geq & \frac{c_2}{N}  +c'_5 \left\{N^{-\frac{1-(2+\alpha)\tau}{1+\alpha-(2+\alpha)\tau}} + \frac{1+VC(\Pi)}{N \mu_N^{\alpha/(1-\tau(2+\alpha))}} + \frac{(VC(\Pi)+1)\log^{\frac{\alpha/\tau}{1+\alpha - \tau(2+\alpha)}}(\max(N,1/\delta))}{N\mu_N^{\alpha/(1+\alpha- (2+\alpha)\tau)}}\right\} +c_5 \frac{1}{N \lambda_N^{\frac{\alpha}{1 - \tau(2+ \alpha)}}},
\end{align*}
for some constant $c_5'$.

Applying Theorem 4 in \cite{farahmand2012regularized} with $\F = \F_l$, $\epsilon = 1/2$ and $\eta = 2^{l} t$, for sufficiently large $N$, we can obtain that 
\begin{align*}
& \Pr\left(\left\{\exists \pi \in \Pi,  I_2(\pi) > t\right\}  \jiao E_N\right) \\
& \leq \sum_{l=0}^\infty \Pr \left[\sup_{h \in \F_l}  \frac{(P-\PN) \left\{h(S, A)\right\}}{P \left\{h(S, A)\right\} + 2^l t} > \frac{1}{2} \right] \\
& \leq \sum_{l=0}^\infty 120 \exp\left\{-c_6\frac{v_{N,l }^2t2^l}{N}\right\} + 2\beta_{x_{N,l}}v_{N, l} \\
& \leq \sum_{l=0}^\infty 120 \exp\left\{-c_6\frac{v_{N,l }^2t2^l}{N}\right\} + 2\beta_0\exp\left(-\beta_1 x_{N, l} + \log v_{N, l}\right) ,
\end{align*}
where the last inequality is based on Assumption \ref{ass: stationarity}.  When $t \geq \frac{\left(4/\beta_1\log(N)\right)^{1/\tau}}{N}$, we have $\log v_{N, l} \leq \frac{1}{2}\beta_1 x_{N, l}$. This will further imply that $2\beta_{x_{N,l}}v_{N, l} \leq 2\beta_0\exp\left(-\beta_1 x_{N, l}/2 \right)$. Then we will have 
\begin{align*}
& \Pr\left(\left\{\exists \pi \in \Pi,  I_2(\pi) > t\right\}  \jiao E_N\right)\\
& \leq \sum_{l=0}^\infty 120 \exp\left\{-c_6\frac{v_{N,l }^2t2^l}{N}\right\} + 2\beta_0\exp\left(-\beta_1 x_{N, l} + \log v_{N, l}\right)\\
& \leq \sum_{l=0}^\infty 120\exp\left(-c_7(Nt)^{1-2\tau}(2l)^{1-2p} \right)+ 2\beta_0\exp\left(-\beta_1 (Nt)^{\tau}(2^l)^p\right)\\
& \leq c_8 \exp\left(-c_9(Nt)^{1-2\tau}\right) + c_{10}\exp\left(-c_{11} (Nt)^{\tau}\right).
\end{align*}
As long as $t$ satisfies all the above constraints, then 
$$
I_2(\eta, Q) \leq \frac{1}{N}\left\{\left(\frac{\log(\frac{2c_8}{\delta})}{c_9}\right)^{\frac{1}{1-2\tau}} + \left(\frac{\log(\frac{2c_{10}}{\delta})}{c_{11}}\right)^{\frac{1}{\tau}} \right\},
$$
with probability at least $1-2\delta$.
Collecting all the conditions on $t$ and combining with the bound of $I_1(\eta, Q)$, we have shown that for sufficiently large $N$ and $\tau  = \tau_1 \leq \frac{1}{3}$, with probability at least $1-2\delta$, the following holds for all  $\pi \in \Pi$:
\begin{align*}
& \norm{\hn(\Hn) - \g(\Hpi)}^2 + \lambda_N J_1^2(\Hn) = I_1(\pi) + I_2(\pi) + I_3(\pi)  \\
& \leqconst  (\mu_N+\lambda_N) (1+\sup_{\pi \in \Pi}J_1^2(\tilde H^\pi)) +\frac{(VC(\Pi)+1)\left[\log(\max(N, 1/\delta))\right]^\frac{1}{\tau_1}}{N} + \frac{1+VC(\Pi)}{N\mu_N^{\frac{\alpha}{1-\tau(2+\alpha)}}}\\
& +\Rem(\pi) +  \mu_N J_1(\Hn) + 
\mu_N\big(1 + \frac{\sqrt{VC(\Pi)+1}\left[\log(\max(N, 1/\delta))\right]^{\frac{1}{2\tau}}}{\sqrt{N\mu_N}} + \frac{\sqrt{1 + VC(\Pi)}}{\sqrt{N}\mu_N^{\frac{1 +\alpha - \tau(2+\alpha)}{2(1- \tau(2+\alpha))}}}\Big)\\
& + N^{-1}  + N^{-\frac{1-(2+\alpha)\tau}{1+\alpha-(2+\alpha)\tau}} + \frac{1+VC(\Pi)}{N \mu_N^{\alpha/(1-\tau(2+\alpha))}} + \frac{(VC(\Pi)+1)\log^{\frac{\alpha/\tau}{1+\alpha - \tau(2+\alpha)}}(\max(N,1/\delta))}{N\mu_N^{\alpha/(1+\alpha- (2+\alpha)\tau)}} +\frac{1}{N \lambda_N^{\frac{\alpha}{1 - \tau(2+ \alpha)}}}\\
&+ \frac{1}{N}\left\{\left(\frac{\log(\frac{2c_8}{\delta})}{c_9}\right)^{\frac{1}{1-2\tau}} + \left(\frac{\log(\frac{2c_{10}}{\delta})}{c_{11}}\right)^{\frac{1}{\tau}} \right\} \\
& \leqconst  (\mu_N+\lambda_N) (1+\sup_{\pi \in \Pi}J_1^2(\tilde H^\pi)) +\frac{(VC(\Pi)+1)\left[\log(\max(N, 1/\delta))\right]^\frac{1}{\tau}}{N} + \frac{1+VC(\Pi)}{N\mu_N^{\frac{\alpha}{1-\tau(2+\alpha)}}}\\
& +\Rem(\pi) +  \mu_N J_1(\Hn) + 
\frac{\sqrt{(VC(\Pi)+1)\mu_N}\left[\log(\max(N, 1/\delta))\right]^{\frac{1}{2\tau}}}{\sqrt{N}} + \frac{\sqrt{1 + VC(\Pi)}}{\sqrt{N}\mu_N^{\frac{\alpha + \tau(2+\alpha)- 1}{2(1- \tau(2+\alpha))}}}\\
& + N^{-\frac{1-(2+\alpha)\tau}{1+\alpha-(2+\alpha)\tau}} + \frac{(VC(\Pi)+1)\log^{\frac{\alpha/\tau}{1+\alpha - \tau(2+\alpha)}}(\max(N,1/\delta))}{N\mu_N^{\alpha/(1+\alpha- (2+\alpha)\tau)}} +\frac{1}{N \lambda_N^{\frac{\alpha}{1 - \tau(2+ \alpha)}}}
\end{align*}
which concludes our proof.

\begin{lemma}[Orthogonality] \label{lemma: orthogonal}
	The function, $e^\pi(\cdot, \cdot)$, satisfies the orthogonality property, i.e., for any state-action function $H(\cdot, \cdot)$, 
	\begin{align}
	\EE\Big[\sum_{t=1}^{T} e^\pi(S_t, A_t) \big(H(S_t, A_t) - \EE\big[\sum_{a'} \pi(a'|S_{t+1})  H(S_{t+1}, a')|S_t, A_t\big]\big) \Big] = 0 
	\end{align}
	As a result, $H^{\pi} \in \argmin_{q }\EE[\frac{1}{T}\sum_{t=1}^T( \EE[\epsilon^\pi(Z_t; H)|S_t, A_t])^2] $ and it is unique up to a constant shift.

\end{lemma}
The proof is straightforward, so we omit here.

%\begin{lemma}\label{lem: VQPI}
%	Let $\V(\Pi, \F_M) = \{ s\mapsto \sum_{a} \pi(a|s) Q(s, a): Q \in \F_M, \pi \in \Pi \}$. For $M \geq 1$, $\log N(\epsilon, \V(\Pi, \F_M), \norm{\cdot}_\infty) \leq \tilde C_1  (M/\epsilon)^{2\alpha}$ where $\tilde C_1  = C_1  + p (\text{diam}(\Theta)L_\Theta F_{\max})^{2\alpha}/(2\alpha)$ with the constant $C_1$   specified in Assumption \ref{b6}.
%\end{lemma}

\begin{lemma} \label{lemma: donsker class prob bound}
	Let $Z_i, i = 1, \cdots, N$ be an exponential $\beta$-mixing stationary sequences and $\F$ be a family of point-wise measurable real-valued functions such that $\norm{f}_\infty \leq F < \infty$ and $\EE[f(Z_1)] = 0$ for all $f \in \F$. In addition, 	$\EE[f(Z_1)^2] \leq \sigma^2$ and $J_1(f) \leq M$. Then under the entropy condition that
	\begin{align*}
	& \log {N}(\epsilon, \F, \norm{\cdot}_\infty) \leqconst v \left(\frac{M}{\epsilon}\right)^{2\alpha}
	\end{align*}
	for some positive constant $v$, then
	with probability at least $1 - \frac{1}{N} - \delta$,
	\begin{align*}
	\sup_{f \in \F} | \sum_{t = 1}^{N}f(Z_t)| &\leqconst \sqrt{vN \log(N)}\sigma^{1-\alpha} + \log(N)v\sigma^{-2\alpha} \\
	&+ M^\alpha\sqrt{vN\log(N)}\sigma^{1-\alpha}+ \log(N)M^{2\alpha}v\sigma^{-2\alpha} \\
	& + \log^\frac{3}{2}(\max(N, 1/\delta))\left\{(v\frac{N}{\log(N)})^{\frac{1}{4}}\sigma^{\frac{1-\alpha}{2}}+\sqrt{v}\sigma^{-\alpha}+M^\frac{\alpha}{2}((v\frac{N}{\log(N)})^{\frac{1}{4}})+M^\alpha\sqrt{v}\sigma^{-\alpha} \right\} \\[0.1in]
	&+ \sigma\sqrt{N} \log(\max(N, 1/\delta)) + \log^2(\max(N, 1/\delta))
	\end{align*}

	%	For $M, \sigma > 0$, let	
	%	\begin{align*}
	%	\F^*= & \Bigdkh{ f: D \mapsto (1/T) \sum_{t=1}^T \big[H(S_t, A_t) - \sum_{a'} \pi(a'|S_{t+1})    q (S_{t+1}, a')\big] g(S_t, A_t): \\
	%		& \qquad \qquad \pi \in \Pi, g \in \G, H \in \F, J_1(q) \leq M, P f^2 \leq \sigma^2, J_2(g) \leq \sup_{\pi \in \Pi} J_2(e^\pi)}
	%	\end{align*}
	%	Under Assumption \peng{blabla}, the following holds  with probability at least $1-\delta$, 
	%	\begin{align*}
	%	\sup_{f \in \F^*}\abs{\Gn f} & \leq C(\delta) \gamma(n, p, M, \sigma) 
	%	\end{align*}
	%	where
	%	\begin{align*}
	%	C(\delta) = K_1  + \frac{4+ 32G_{\max} F_{\max} K_1}{ (8/3) G_{\max}F_{\max}} + (8/3) G_{\max}F_{\max} \log(1/\delta)
	%	\end{align*}
	%	\begin{align*}
	%	\gamma(n, p, M, \sigma)  & = (\sqrt{p}\sigma^{1-\alpha} + p\sigma^{-2\alpha} n^{-1/2} + M^{\alpha}  \sqrt{p}\sigma^{1-\alpha}   + M^{2\alpha} p\sigma^{-2\alpha} n^{-1/2} + \sigma + n^{-1/2} \\
	%	&\qquad + p^{1/4} n^{-1/4} \sigma^{(1-\alpha)/2} + \sqrt{p} \sigma^{-\alpha}n^{-1/2} + M^{\alpha/2} p^{1/4} n^{-1/4} \sigma^{(1-\alpha)/2} + M^{\alpha}\sqrt{p} \sigma^{-\alpha}n^{-1/2} ) 
	%	\end{align*}
\end{lemma}
\textbf{Proof of Lemma \ref{lemma: donsker class prob bound}}
We apply the Berbee's coupling lemma \cite{dedecker2002maximal}, which can approximate $\sup_{f \in \F} | \sum_{t = 1}^{N}f(Z_t)|$. By Lemma 4.1 of \cite{dedecker2002maximal}, we can construct a sequence of random variables $\left\{Z^{0}_t\right\}_{t=1}^N$ such that $\sup_{f \in \F} | \sum_{t = 1}^{N}f(Z_t)| = \sup_{f \in \F} | \sum_{t = 1}^{N}f(Z^0_t)|$, and that the sequence $\left\{Z^0_{2kx_N+j} \right\}_{j=1}^{x_N}$ for $k= 0, \cdots, (v_N-1)$ is i.i.d and so is $\left\{Z^0_{(2k+1)x_N+j} \right\}_{j=1}^{a_N}$ $k= 0, \cdots, (v_N-1)$ are i.i.d with probability at least $1 - \frac{N\boldsymbol{\beta}(x_N)}{x_N}$. Here we assume we can divide the index $\left\{1, \cdots, N\right\}$ into $2v_N$ block with equal length $x_N$. Denote the remainder index set as $R_N$ and without loss of generality assume that $|R_N| \leq x_N$. Then we can show that
\begin{align*}
\sup_{f \in \F} | \sum_{t = 1}^{N}f(Z_t)| &\leq \sum_{j = 1}^{x_N}\sup_{f \in \F} | \sum_{k = 0}^{2v_N-1}f(Z^0_{kx_N+j})| + \sum_{f \in \F} |\sum_{j \in R_N}f(Z^0_j)|\\[0.1in]
& \leq \sum_{j = 1}^{2x_N}\sup_{f \in \F} | \sum_{k = 0}^{v_N-1}f(Z^0_{2kx_N+j})| + |R_N| M\\[0.1in]
& \leq \sum_{j = 1}^{2x_N}\sup_{f \in \F} | \sum_{k = 0}^{v_N-1}f(Z^0_{2kx_N+j})| + x_NM,
\end{align*}
where the last inequality holds because $|R_N| \leq x_N$.

As we know $Z^0_{2kx_N+j}$ is i.i.d. for $k= 0, \cdots, v_N-1$. Then we can first apply Talagrand inequality to show that for any $t>0$,
with probability at least $1-\exp(t)$, we have
\begin{align*}
\sup_{f \in \F} \abs{\sum_{k = 0}^{v_N-1}f(Z^0_{2kx_N+j})}&\leq  \EE\sup_{f \in \F}\abs{\sum_{k = 0}^{v_N-1}f(Z^0_{2kx_N+j})} \\[0.1in]
&+\sqrt{4Ft\EE\sup_{f \in \F}\abs{\sum_{k = 0}^{v_N-1}f(Z^0_{2kx_N+j})} + 2v_N\sigma^2t} + \frac{Ft}{3},
\end{align*}
which can further imply that with probability at least $1-\delta$
\begin{align*}
\sum_{j=1}^{2x_N}\sup_{f \in \F} \abs{\sum_{k = 0}^{v_N-1}f(Z^0_{2kx_N+j})}&\leq  \sum_{j=1}^{2x_N}\EE\sup_{f \in \F}\abs{\sum_{k = 0}^{v_N-1}f(Z^0_{2kx_N+j})}\\[0.1in]
& +\sum_{j=1}^{2x_N}2 \sqrt{F\log(\frac{2x_N}{\delta})}\sqrt{\EE\sup_{f \in \F}\abs{\sum_{k = 0}^{v_N-1}f(Z^0_{2kx_N+j})}}\\[0.1in]
& + 2x_N\sigma\sqrt{2v_N\log(\frac{2x_N}{\delta})} + 2x_N\frac{F\log(\frac{2x_N}{\delta})}{3}
\end{align*}
Then applying maximal inequality with uniform entropy condition, we can show that
$$
\EE\sup_{f \in \F}\abs{\sum_{k = 0}^{v_N-1}f(Z^0_{2kx_N+j})} \leqconst \sqrt{v v_N}\sigma^{1-\alpha} + v\sigma^{-2\alpha} + M^\alpha\sqrt{vv_N}\sigma^{1-\alpha}+ M^{2\alpha}v\sigma^{-2\alpha}.
$$
Summarizing together and choosing $x_N = 2\log(N)$, we obtain that with probability at least $1 - \frac{1}{N} - \delta$,
\begin{align*}
\sup_{f \in \F} | \sum_{t = 1}^{N}f(Z_t)| &\leqconst \sqrt{vN \log(N)}\sigma^{1-\alpha} + \log(N)v\sigma^{-2\alpha} \\
&+ M^\alpha\sqrt{vN\log(N)}\sigma^{1-\alpha}+ \log(N)M^{2\alpha}v\sigma^{-2\alpha} \\
& + \log^\frac{3}{2}(\max(N, 1/\delta))\left\{(v\frac{N}{\log(N)})^{\frac{1}{4}}\sigma^{\frac{1-\alpha}{2}}+\sqrt{v}\sigma^{-\alpha}+M^\frac{\alpha}{2}(v\frac{N}{\log(N)})^{\frac{1}{4}}+M^\alpha\sqrt{v}\sigma^{-\alpha} \right\} \\[0.1in]
&+ \sigma\sqrt{N} \log(\max(N, 1/\delta)) + \log^2(\max(N, 1/\delta))
\end{align*}
which concludes our proof by dividing $N$ at both sides. In particular, when $M, \sigma$ are all constants, we can show with probability at least $1 - \frac{1}{N} - \delta$,
$$
\sup_{f \in \F} | \PP_Nf(Z)| \leqconst \sqrt{\frac{v}{N}}\log(\max(N, \frac{1}{\delta})).
$$
\subsection{Regret Bound}
\textbf{Proof of Theorem \ref{thm: policy learning}}
Let $\pi^*$ is the in-class optimal policy and assume $\pi^* \in \Pi$ and denote $\widehat \pi_N = \widehat{\pi}_N^c$ to indicate its dependency on $N$. We bound the regret by
\begin{align*}
\Regret (\widehat \pi_N) & = \sup_{\pi \in \Pi, |\beta|\leq R_{\max}} M(\beta, \pi) -  \sup_{|\beta|\leq R_{\max}}M(\beta, \widehat \pi_N) \leq M(\beta^\ast, \pi^*) - M(\widehat \beta, \widehat \pi_N) \\
& = (\widehat M_N - M)(\widehat \beta, \widehat \pi_N) - (\widehat M_N - M)(\beta^\ast, \pi^*) +  \widehat M_N(\beta^\ast, \pi^*) - \widehat M_N(\widehat \pi_N) \\
& \leq (\widehat M_N - M)(\widehat \beta, \widehat \pi_N) - (\widehat M_N - M)(\beta^\ast, \pi^*)
\end{align*}
Define
$$
\phi_\beta^\pi(S, A, S') = 	
\omega^\pi(S, A) [\beta - \frac{1}{1-c}(\beta - R)_+ +  \sum_{a'} \pi(a'|S') Q^\pi(S', a') - Q^\pi(S, A) -  M(\beta, \pi) ].
$$
Define the remainder term $\Rem_N(\beta, \pi) =(\Mn (\beta, \pi) - M(\beta, \pi)) - \PN \phi^{\pi}_\beta$. Letting $\bar B = [-R_{\max}, R_{\max}]$, We then have
\begin{align*}
\Regret (\widehat \pi_N)  & \leq \PN (\phi^{\widehat \pi_N}_{\widehat \beta} - \phi^{\pi^*}_{\beta^\ast}) + (\Rem_N(\widehat \beta, \widehat \pi_N) - \Rem_N(\beta^\ast, \pi^*)) \\
& \leq \sup_{\pi \in \Pi, \beta \in \bar B} ~ \PN (\phi^{\pi}_\beta - \phi^{\pi^*}_\beta) + 2 \sup_{\pi \in \Pi, \beta \in \bar B} \abs{\Rem_N(\beta, \pi)}.
\end{align*}

%		Below we show that (i) $\sup_{\pi \in \Pi} \Pn (\phi^{\pi} - \phi_{\pi^*})  = O_P(n^{-1/2})$ and (ii) $\sup_{\pi \in \Pi} \abs{\Rem_n(\pi)}  = o_P(n^{-1/2})$ and thus $\Regret(\widehat \pi_n) = O_P(n^{-1/2})$.

\paragraph{(i) Leading Term} 
For any $(s, a, s', r)$, we have
\begin{align*}
& \abs{\omega^{\pi_1}(s, a) (\beta_1 - \frac{1}{1-c}(\beta_1-\RR(s))_+ + U^{\pi_1}(s, a, s') - M(\beta_1,\pi_{1})) \\
	&- \omega^{\pi_{2}}(s, a) (\beta_2 - \frac{1}{1-c}(\beta_2-\RR(s))_+ + U^{\pi_2}(s, a, s') - M(\beta_2, \pi_2))}\\
& \leq 	\frac{2-c}{1-c}\sup_{\pi \in \Pi}\| \omega^\pi \|_\infty | \beta_1 - \beta_2| + 2 \left(\frac{2}{1-c}R_{\max}+F_{\max} \right)|\omega^{\pi_1}(s, a) - \omega^{\pi_2}(s, a)|\\
& +  \sup_{\pi \in \Pi} \| \omega^\pi \|_\infty |U^{\pi_1}(s, a, s') - U^{\pi_2}(s, a, s')| + \sup_{\pi \in \Pi} \| \omega^\pi \|_\infty |M(\beta_1, \pi_1) - M(\beta_2, \pi_2)|	\end{align*}
By our assumption, we know
\begin{align*}
&\abs{\omega^{\pi_1}(s, a) - \omega^{\pi_2}(s, a)} 
%		\leq \frac{1}{d_D(s, a)} (C_d \norm{\theta_1 - \theta_2}_2 + L_{\Pi} \norm{\theta_1 - \theta_2}) 
\leqconst d_{\Pi}(\pi_1, \pi_2)\\
&	\abs{M(\beta_1, \pi_1) - M(\beta_2, \pi_2)}  
%		R_{\max} \abs{\A} L_{\Pi} \norm{\theta_1 - \theta_2}_{2} + R_{\max}  C_d\norm{\theta_1 - \theta_2}_2 =\\
\leqconst d_{\Pi}(\pi_1, \pi_2) + |\beta_1 - \beta_2| \\
&\abs{U^{\pi_1}(s, a, s') - U^{\pi_2}(s, a, s')} 
%		\\&\leq \abs{Q^{\pi_{\theta_1}}(s, a) - Q^{\pi_{\theta_2}}(s, a)} + \abs{\sum_{a'} \pi_{\theta_1}(a'|s') Q^{\pi_{\theta_1}}(s', a') - \pi_{\theta_2}(a'|s') Q^{\pi_{\theta_2}}(s', a')} \\
%	&	\leq \abs{\V(\pi_{\theta_1}) - \V(\pi_{\theta_2})} + \norm{V^{\pi_{\theta_1}} - V^{\pi_{\theta_2}}}_{\infty} \left(1 + L_{\Pi} \abs{\A} + 1\right) \\ & 
\leqconst d_{\Pi}(\pi_1, \pi_2) + |\beta_1 - \beta_2|.
\end{align*}
Then we have $$\abs{\phi^{\pi_1}_{\beta_1}(s, a,s') - \phi^{\pi_2}_{\beta_2}(s, a,s')}\leqconst d_\Pi(\pi_1, \pi_2) + |\beta_1 - \beta_2|$$	
On the other hand, 
\begin{align*}
\abs{\phi^{\pi}_\beta(s, a,s')} \leq \phi_{\max} := 2 (\frac{1}{1-c}R_{\max} + F_{\max}) \cdot \sup_{\pi \in \Pi} \norm{\omega^\pi}_\infty < \infty
\end{align*}
We will apply the maximal inequality with the bracketing number. This only requires a slight modification of Lemma \ref{lemma: donsker class prob bound}. We can show that 
\begin{align*}
&  \sup_{\pi \in \Pi, \beta \in \bar B} ~ \PN (\phi_\beta^{\pi} - \phi_{\beta^\ast}^{\pi^*})  \leqconst \log(\max(N, 1/\delta))\sqrt{\frac{\Sigma}{N}}J_{[]}(\phi_{\max}, \F^*, L_2),
\end{align*}
with probability $(1-\delta - \frac{1}{\log(N)})$,
where $\F^* = \{ \phi_\beta^{\pi} - \phi_{\beta^\ast}^{\pi^*}: \pi \in \Pi \}$, $\Sigma = \underset{\pi \in \Pi,|\beta| \leq R_{\max}}{\sup} \EE\left[\psi^2(Z;U^{\pi,\beta}, \omega^\pi) \right]$ and the bracketing entropy $J_{[]}(\phi_{\max}, \F^*, L_2) = \int_{0}^{\phi_{\max}} \sqrt{ \log N_{[]}(\epsilon, \F^*, L_2)} d\epsilon$. Using the Lipschitz property gives  
\begin{align*}
J_{[]}(\phi_{\max}, \F^*, L_2) & \leqconst  \int_{0}^{\phi_{\max}} \sqrt{ \log N((R_{\max})^{-1} \epsilon, \bar B, \norm{\cdot}_2) + \log N(\epsilon, \Pi, d_{\Pi}(\bullet))} d\epsilon\\
& \leqconst \sqrt{VC(\Pi) + 1}.
\end{align*}

\paragraph{(ii) Remainder Term} 
%		\ \peng{rewrite the proofs for choosing $k$ s.t. $\beta_k$ is close to the target rate. And get the exact rate...}
For the ease of notation, define
\[
f(\omega, U, \beta, \pi):  (S, A, S') \mapsto   \omega(S, A) (\beta - \frac{1}{1-c}(\beta - R)_+ + U(S, A, S) - M(\beta, \pi))
\]
Note that we have $\phi^{\pi}_\beta = f(\omega^\pi, U^\pi, \beta, \pi)$. Let $\widehat \phi^\pi_\beta= f(\widehat \omega^\pi_N, \UpiN,\beta, \pi)$ be a ``plug-in'' estimator of $\phi^\pi$. Since the ratio estimator satisfies $\PN \omega^\pi_N(S, A) = 1$ by construction, we have
\begin{align*}
& \Rem_N(\beta, \pi)  = \Mn (\beta, \pi) - M(\beta, \pi) - \PN \phi^{\pi}_\beta \\
& = (\PN - P)(\widehat \phi^\pi_\beta - \phi^\pi_\beta) + P (\widehat \phi^\pi_\beta - \phi^\pi_\beta)
\end{align*}
This implies that $$\sup_{\pi \in \Pi, \beta \in \bar B} \abs{\Rem_N(\beta, \pi)} \leq \sup_{\pi \in \Pi, \beta \in \bar B} \abs{ P (\widehat \phi^\pi_\beta - \phi^\pi_\beta)} + \sup_{\pi \in \Pi, \beta \in \bar B} \abs{(\PN - P)(\widehat \phi^\pi_\beta - \phi^\pi_\beta) } $$.

Consider the first term. The doubly-robustness structure of the estimating equation, implies that
\begin{align*}
& P(\widehat \phi^\pi_\beta  - \phi^\pi_\beta)  = P ( f(\widehat \omega^\pi_N, \UpiN, \pi) -  f(\omega^\pi, U^\pi, \pi) ) \\
& = P [f(\widehat \omega^\pi_N, \UpiN, \pi) - f(\widehat \omega^\pi_N, \Upi,\beta, \pi)  +  f(\widehat \omega^\pi_N, \Upi, \beta,\pi)   - f(\omega^\pi, U^\pi,\beta, \pi) ] \\
& = P [f(\widehat \omega^\pi_N, \Upi,\beta, \pi)   - f(\omega^\pi, U^\pi, \beta, \pi)]   + \Big(P [f(\widehat \omega^\pi_N, \UpiN, \pi) - f(\widehat \omega^\pi_N, \Upi, \beta, \pi) ] \\
& \qquad - P [f( \omega^\pi, \Upin, \pi) - f(\omega^\pi, \Upi, \beta,\pi)] \Big)+ P [f( \omega^\pi, \UpiN, \beta,\pi) - f(\omega^\pi, \Upi, \beta,\pi)] \\
& = \EE\Big[ (1/T) \sum_{t=1}^T (\widehat \omega^\pi_N - \omega^\pi)(S_t, A_t) (R_{t+1} + \Upi(S_t, A_t, S_{t+1}) - M(\beta, \pi) \Big]  \\
& \qquad + \EE\Big[ (1/T) \sum_{t=1}^T (\widehat \omega^\pi_n - \omega^\pi)(S_t, A_t) \cdot (\widehat U^\pi_n-U^\pi)(S_t, A_t, S_{t+1}) \Big]\\
& \qquad + \EE\Big[ (1/T) \sum_{t=1}^T \omega^\pi(S_t, A_t) (\UpiN - U^\pi)(S_{t}, A_{t}, S_{t+1}) \Big]\\
& =\EE\Big[  (1/T) \sum_{t=1}^T   (\widehat \omega^\pi_N - \omega^\pi)(S_t, A_t) \cdot (\UpiN -U^\pi)(S_t, A_t, S_{t+1})\Big]
%  & \leq \norm{\widehat w^{-k} - \omega^\pi}_2 \cdot \norm{\Upin-U^\pi}_2
\end{align*}
where the last equality holds by noting  $\sum_{s, a} \EE[(\UpiN - \Upi)(S_{t}, A_t, S_{t+1})|S_t=s, A_t=a] d^\pi(s, a) = 0$. 
Furthermore, applying Cauchy inequality twice gives
\begin{align*}
& \abs{P( \widehat \phi^\pi_N  - \phi^\pi)} =  \abs{(1/T) \sum_{t=1}^T	\EE[ (\widehat \omega^\pi_N- \omega^\pi)(S_t, A_t) \cdot(\UpiN-U^\pi)(S_t, A_t, S_{t+1}) ]} \\
& \leq (1/T) \sum_{t=1}^T \sqrt{\EE\big[ (\widehat \omega^\pi_N - \omega^\pi)^2(S_t, A_t)\big]} \cdot \sqrt{\EE[(\UpiN-U^\pi)^2(S_t, A_t, S_{t+1})]}\\
& \leq \sqrt{(1/T) \sum_{t=1}^T \EE[ (\widehat \omega^\pi_N- \omega^\pi)^2(S_t, A_t)]} \cdot \sqrt{(1/T) \sum_{t=1}^T \EE[(\UpiN-U^\pi)^2(S_t, A_t, S_{t+1})]} \\
& = \norm{\widehat \omega^\pi_N - \omega^\pi} \cdot \norm{\widehat U^\pi_N - U^\pi} 
\end{align*}	
Using Theorem \ref{thm: ratio} and Theorem \ref{thm: value}, we can show that
\begin{align*}
\sup_{\pi \in \Pi, \beta \in \bar B} \abs{P( \widehat \phi^\pi_\beta  - \phi^\pi_\beta)} & \leq  \sup_{\pi \in \Pi, \beta \in \bar B}  \dkh{\norm{\widehat \omega^\pi_N - \omega^\pi} \cdot \norm{\widehat U^\pi_N - U^\pi}} \\
& \leq (\sup_{\pi \in \Pi}  \norm{\widehat \omega^\pi_N - \omega^\pi}) ( \sup_{\pi \in \Pi, \beta \in \bar{B}} \norm{\widehat U^\pi_N - U^\pi} ) \\
& \leqconst N^{-r_k/2} (1+VC(\Pi))[\log(\max(N, 1/\delta))]^{\frac{1+\alpha/2}{1+\alpha} + 1/2}\\
& \times	\left[\log\left(\max(1/\delta, N)\right)\right]^{\frac{2+\alpha}{2(1+\alpha)}} N^{-\frac{1}{2(1+\alpha)}}\\
& \leqconst N^{-(r_k+\frac{1}{1+\alpha})/2} (1+VC(\Pi))[\log(\max(N, 1/\delta))]^{\frac{2+\alpha}{1+\alpha} + 1/2},
\end{align*}
with probability at least $1 - (k+6)\delta - k/N$.
As long as $k\geq 3$, choosing $\delta = \frac{1}{N}$, the above term decays faster than $\frac{1}{\sqrt{N}}$.
Now we consider the second term. 	
%			Suppose $n$ is large enough such that $C_3^{[k]}(\delta) \iota^{\omega_k/2}\sqrt{p}n^{-(1-1/(1+\alpha)+ \beta_k)/2} < (1/2) P e^\pi$, then $\abs{\Pn \gn(\qn) - P e^\pi} \leq (1/2) P e^\pi$. Thus we have $\Pn \gn(\qn) \geq (1/2) P e^\pi > 0$
%		\begin{align*}
%		\abs{(\Pn - P)(\widehat \phi^\pi_n - \phi^\pi) }  \\
%		\leq \abs{(\Pn - P)(\widehat \phi^\pi_n - \phi^\pi) }
%		\end{align*}
Define 
$$
(I) \triangleq N^{-r_k} (1+VC(\Pi))[\log(\max(N, 1/\delta))]^{\frac{2+\alpha}{1+\alpha} + 1},
$$
and 
$$
(II) \triangleq \left(1+VC(\Pi)\right)\left[\log\left(\max(1/\delta, N)\right)\right]^{\frac{2+\alpha}{1+\alpha}} N^{-\frac{1}{1+\alpha}},
$$
As we know $\Pr(E_N) \geq 1- (6+k)/N- k/\log(N)$, where $E_N = E_{N, 1} \jiao E_{N, 2}$ and 
\begin{align*}
& E_{N, 1} = \{ \norm{\widehat \omega^\pi_N - w^\pi}^2 \leqconst (I), J_2(\widehat e^\pi_N) \leqconst N^{\frac{ \frac{1}{1+\alpha}}{2}-\frac{r_k}{2}}, \forall \pi \in \Pi\}\\
& E_{N, 2} = \{ \norm{\UpiN - \Upi}^2 \leqconst  (II), J_1(\widehat Q_N^\pi) \leqconst 1, \forall \pi \in \Pi\}.
\end{align*}
By the previous argument, we know that for $N$ sufficiently large, we have
$$
\PN\hn \geq \frac{1}{2}P e^\pi.
$$
Therefore
$$
J_2(\omega_N^\pi) \leqconst N^{\frac{ \frac{1}{1+\alpha}}{2}-\frac{r_k}{2}}.
$$
Then 
under this event $E_N$, we have
\begin{align*}
\sup_{\pi \in \Pi, \beta \in \Pi} \bigdkh{\abs{(\PN - P)( \widehat \phi^\pi_\beta  - \phi^\pi_\beta) }} \leq \sup_{f \in \F^*, Pf^2 \leqconst \zeta(N)}  \abs{ (\PN - P) f} 
\end{align*}
where $\zeta(N) = (I)$ and $$\F^* = \{ f: (S, A, S') \mapsto  g(S, A, S') - \phi^\pi_\beta(S, A, S') \given \beta \in \bar B,   \pi \in \Pi, g \in \G^*\}$$
\begin{align*}
\G^* = \{& g: (S, A, S')   \mapsto  w(s, a) (\beta - \frac{1}{1-c}\left(\beta - \RR(s)\right)_+ +\sum_{a' \in \cal A}\pi(a'|s')Q(s', a') - \eta) \given \\
&    {	J_2(w) \leqconst N^{\frac{ \frac{1}{1+\alpha}}{2}-\frac{r_k}{2}}}, J(Q) \leqconst 1, \pi \in \Pi, \eta, \beta \in [-R_{\max}, R_{\max}]\}
\end{align*}
One can show that
$$
\log(N(\epsilon, \F^*, \norm{\cdot}_\infty)) \leqconst (VC(\Pi) + 1)\left(\frac{M}{\epsilon}\right)^{2\alpha},
$$
where $M= N^{\frac{ \frac{1}{1+\alpha}}{2}-\frac{r_k}{2}}$. Applying Lemma \ref{lemma: donsker class prob bound} with $v= VC(\Pi)+1$, $M$ and $\sigma^2 = \zeta(N)$, we can show that with probability $1-1/N$,
$$
\sup_{\pi \in \Pi, \beta \in \Pi} \bigdkh{\abs{(\PN - P)( \widehat \phi^\pi_\beta  - \phi^\pi_\beta) }}  \leqconst \sqrt{VC(\Pi)+1}~o(\frac{1}{\sqrt{N}}),
$$
for $k\geq 2$ and $N$ sufficiently large. Summarizing together, we can show that with probability at least $1 - (2k+10)/N$
$$
\Regret(\widehat \pi_N)\leqconst \sqrt{\frac{\Sigma(VC(\Pi)+1)}{N}}\log(N).
$$

\subsection{Statistical Efficiency of the Proposed Estimator}
\textbf{Proof of Theorem \ref{thm: efficiency bound}}
As we have shown in the proof of Theorem \ref{thm: policy learning}, for any $\pi \in \Pi$ and $|\beta| \leq R_{\max}$, 
$$
\widehat M_N(\beta, \pi) - M_N(\beta, \pi) = \Rem_N(\beta, \pi) =  o_p(\frac{1}{\sqrt{N}}).
$$
Denote $V^2 = \EE\left[\psi^2(Z;U^\pi, \omega^\pi)\right]$. We can show that
$$
\sqrt{N}\frac{ M_N(\beta, \pi) - M(\beta, \pi)}{V} \xrightarrow{d} \N(0, 1). 
$$
Recall that
$$
M_N(\beta, \pi) - M(\beta, \pi) = \PN \psi(Z; U^\pi, \omega^\pi),
$$
which is sum of martingale differences. We apply Corollary 2 in \cite{jones2004markov}. By Assumption \ref{ass: stationarity} and $\psi$ is uniformly bounded, we have
$$
\sqrt{N}\PN \phi(Z; U^\pi, \omega^\pi) \xrightarrow{d}  \N(0, V^2).
$$
The remaining is to show $EB(N) = V^2$.
By Lemma A.1 in \cite{liao2020batch}, under some regularity condition, we can show that
$$
\triangledown M(\varpi_0) = \EE[\triangledown L_{\varpi_0}(\{D_i\}_{i=1}^n)\PN\psi(Z, U^\pi, \omega^\pi)].
$$
Then by the Cauchy-Schwarz inequality, we have
\begin{align*}
EB(N) &= N\sup \left\{\triangledown^T M(\varpi_0) \left\{\EE\left[ \triangledown L_{\varpi_0}(\{D_i\}_{i=1}^n)\triangledown^T L_{\varpi_0}(\{D_i\}_{i=1}^n)\right]\right\}^{-1} \triangledown M(\varpi_0)\right\}\\[0.1in]
& \leq N\EE\left[\PN\psi(Z, U^\pi, \omega^\pi) \left(\PN\psi(Z, U^\pi, \omega^\pi)\right)^T\right]\\[0.1in]
& = \EE\left[\PN\psi^2(Z, U^\pi, \omega^\pi)\right]\\[0.1in]
& =V^2,
\end{align*}
where the second equality uses $\EE[\psi(Z_i, U^\pi, \omega^\pi)\psi(Z_j, U^\pi, \omega^\pi)] = 0$ for $i\neq j$ and the last equality is based on the stationarity property given in Assumption \ref{ass: stationarity}. We conclude our proof by using a similar argument in the proof of Theorem 2 in \cite{kallus2019efficiently} to show that the upper bound $V^2$ is the supremum over all regular parametric models.

\textbf{Proof of Theorem \ref{thm: efficiency bound for CVaR}}
As we have shown in the proof of Theorem \ref{thm: efficiency bound}, for any $\pi \in \Pi$ and $|\beta| \leq R_{\max}$, 
$$
\sqrt{N}\frac{ M_N(\beta, \pi) - M(\beta, \pi)}{V(\beta)} \xrightarrow{d} \N(0, 1),
$$
where $V^2(\beta) = \EE\left[\psi^2(Z;U^{\pi, \beta}, \omega^\pi)\right]$. Then by functional delta theorem (See Theorem 5.7 of \cite{shapiro2021lectures}) and the assumption that $\max_{\beta \in \mathbb{R}} M(\beta, \pi)$ is unique, we have
$$
\sqrt{N}\frac{ \max_{\beta \in \mathbb{R}}M_N(\beta, \pi) - \max_{\beta \in \mathbb{R}} M(\beta, \pi)}{V(\beta^\ast(\pi))} \xrightarrow{d} \N(0, 1).
$$
The remaining is to show $V(\beta^\ast(\pi))$ is the efficiency bound. By Danskin theorem (e.g., \cite{danskin2012theory}) and some regularity condition, one can show that
$$
\triangledown M(\pi, \beta^\ast(\pi); \varpi_0) = \EE[\triangledown L_{\varpi_0}(\{D_i\}_{i=1}^n)\PN\phi(Z, U^{\pi,\beta^\ast(\pi)}, \omega^\pi)].
$$
Then by the Cauchy-Schwarz inequality, we have
\begin{align*}
EB(N) &= N\sup \left\{\triangledown^T M(\pi, \beta^\ast(\pi); \varpi_0)  \left\{\EE\left[ \triangledown L_{\varpi_0}(\{D_i\}_{i=1}^n)\triangledown^T L_{\varpi_0}(\{D_i\}_{i=1}^n)\right]\right\}^{-1} \triangledown M(\pi, \beta^\ast(\pi); \varpi_0) \right\}\\[0.1in]
& \leq N\EE\left[\PN\phi(Z, U^\pi, \omega^\pi) \left(\PN\phi(Z, U^{\pi,\beta^\ast(\pi)}, \omega^\pi)\right)^T\right]\\[0.1in]
& = \EE\left[\PN\phi^2(Z, U^{\pi,\beta^\ast(\pi)}, \omega^\pi)\right]\\[0.1in]
& =V^2(\beta^\ast(\pi)).
\end{align*}
We conclude our proof by using a similar argument in the proof of Theorem 2 in \cite{kallus2019efficiently} to show that the upper bound $V^2(\beta^\ast(\pi))$ is the supremum over all regular parametric models.

\end{appendices}

\end{document}